\documentclass{amsart}
 \usepackage{amsbsy,amssymb,amscd,amsfonts,latexsym,amstext,delarray,
 amsmath,graphicx}
 \usepackage[dvips]{epsfig}
\usepackage{hyperref}

\def\cutint{{\int \!\!\!\!\!\! -}}

\newtheorem{thm}{Theorem}[section]
\newtheorem{prop}[thm]{Proposition}
\newtheorem{cor}[thm]{Corollary}
\newtheorem{lem}[thm]{Lemma}

\newtheorem{defn}[thm]{Definition}
\newtheorem{rem}[thm]{Remark}

\def\qqq{\,,\quad \forall}

\def\Aut{{\rm Aut}}

\def\Det{{\rm Det}}

\def\Dom{{\rm Dom}}
\def\End{{\rm End}}

\def\id{{\rm id}}

\def\Sp{{\rm Spec}}

\def\Tr{{\rm Tr}}

\def\C{{\mathbb C}}

\def\N{{\mathbb N}}

\def\R{{\mathbb R}}
\def\Z{{\mathbb Z}}

\def\Tr{{\rm Tr}}

\def\cA{{\mathcal A}}
\def\cB{{\mathcal B}}

\def\cE{{\mathcal E}}

\def\cH{{\mathcal H}}

\def\cK{{\mathcal K}}
\def\cL{{\mathcal L}}

\def\cR{{\mathcal R}}
\def\cS{{\mathcal S}}

\newcommand{\ie}{{\it i.e.\/}\ }
\newcommand{\aee}{{\it a.e.\/}\ }
\newcommand{\eg}{{\it e.g.\/}\ }
\newcommand{\cf}{{\it cf.\/}\ }

\def\id{{\mbox{Id}}}

\def\dim{{\mbox{dim}}}

\def\Der{{\mbox{Der}}}

\def\End{{\mbox{End}}}

\def\Rank{{\mbox{Rank}}\,}

\def\weyl{Weyl norms }
\def\limom{{\rm Lim}_{\omega}}
\def\limo{{\rm Lim}}
\def\limtwo{{\rm Lim}^2}

\newcommand{\nil}[1]{}

\def\smooth{expable}
\def\smoothness{expability}
\def\axiom{condition }
\def\axioms{conditions }
\def\axiomss{conditions}
\title
{On the spectral characterization of manifolds}
\author[Connes]{Alain Connes}
\address{A.~Connes: Coll\`ege de France \\
3, rue d'Ulm \\ Paris, F-75005 France
\\ I.H.E.S. and Vanderbilt
University} \email{alain\@@connes.org}

\date{}
\begin{document}
\maketitle \vspace{2cm}

\begin{abstract} We show that the first five of the axioms we had formulated
on spectral triples suffice (in a slightly stronger form) to characterize the
spectral triples associated to smooth compact manifolds. The algebra, which is
assumed to be commutative, is shown to be isomorphic to the algebra of all
smooth functions on a unique smooth oriented compact manifold, while the
operator is shown to be of Dirac type and the metric to be Riemannian.
\end{abstract}

\tableofcontents

\section{Introduction}

  The problem of spectral characterization of
manifolds was initially formulated as an open question in
\cite{CoSM}. The issue is to show that under the simple \axioms of
\cite{CoSM} on a spectral triple $(\cA,\cH,D)$, with  $\cA$ {\em
commutative},  the algebra $\cA$ is the algebra $C^\infty(X)$ of
smooth functions on a (unique) smooth compact manifold $X$. The five
\axioms (\cite{CoSM}), in dimension $p$, are

\begin{enumerate}
  \item  The  $n$-th characteristic value of the resolvent of $D$ is
$O(n^{-1/p})$.
  \item $\left[ [D,a],b\right] = 0 \qquad  \forall \,
a,b \in \cA$.
  \item  For any $a\in \cA$ both $a$
and $[D,a]$ belong to the domain of $\delta^m$, for any integer $m$
where $\delta$ is the derivation: $\delta(T)=[|D|,T]$.
  \item  There exists a
Hochschild cycle $c \in Z_p (\cA ,\cA)$ such that $\pi_D (c)=1$ for
$p$ odd, while for $p$ even, $\pi_D (c) =\gamma$ is a $\Z/2$
grading.
  \item Viewed as an $\cA$-module the space $\cH_{ \infty} =
  \cap \Dom D^m$ is
finite and projective. Moreover the following equality defines a
hermitian structure $( \ |\ )$ on this module: $
    \langle  \xi ,a \,\eta  \rangle = \cutint \, a (\xi |\eta) \, |D|^{-p}
\,, \ \forall   a \in \cA   ,   \forall  \xi ,\eta \in \cH_{ \infty}
$.
\end{enumerate}
The notations are recalled at the beginning of \S \ref{prelem}
below. The strategy of proof was outlined briefly in \cite{CoSM}. It
consists in using the components $a_\alpha^j$ ($j>0$) of the cycle
$c=\sum a_\alpha^0\otimes a_\alpha^1\otimes \cdots \otimes
a_\alpha^p$ as tentative local charts. There are three  basic
difficulties:
\begin{enumerate}
  \item[a)] Show that the spectrum $X$ of $\cA$ is {\em large enough} so
  that the range of ``local charts" $a_\alpha$ contains an open set in $\R^p$.
  \item[b)] Show that the joint spectral measure of the components $a_\alpha^j$ ($j>0$)
  of a ``local chart" is the Lebesgue measure.
  \item[c)] Apply the basic inequality
(\cite{Co-action1}, \cite{Co-book}, Proposition IV.3.14)
  giving an upper bound on the Voiculescu obstruction \cite{Voic1} and use
  \cite{Voic1} Theorem 4.5  to show that the
  ``local charts" are locally injective.
\end{enumerate}
 In a recent paper \cite{ReVa}, Rennie and Varilly considered the above challenging
problem. The paper \cite{ReVa} is a courageous attempt which contains a number
of interesting ideas and a useful smooth calculus but also, unfortunately,
several gaps, each being enough to invalidate the proof of the claimed result.

I will show in this paper how to prove a), b), c). I have tried to be very
careful and give detailed proofs.  The way to prove a) uses a new ingredient:
the implicit function theorem (whose presence is not a real surprise). We shall
first assume that continuous $*$-derivations of $\cA$ exponentiate, \ie are
generators of one-parameter groups (of automorphisms of $\cA$). Then  most of
the work, done in \S\S \ref{dissipsect}, \ref{sectselfadj}, is to show that
this hypothesis can be removed. In this very technical part of the paper we
show that enough self-adjoint derivations of $\cA$ exponentiate. We first prove
in \S \ref{dissipsect} that enough derivations are {\em dissipative} for the
$C^*$-algebra norm. We then proceed in \S \ref{sectselfadj} and use the
self-adjointness of $D$ and the third \axiom (regularity) in the strong form,
to show the surjectivity of the resolvent, and apply the Hille-Yosida Theorem
to integrate these derivations into one-parameter groups of automorphisms   of
the $C^*$-algebra. We then show that they are continuous for the Sobolev norms
and define automorphisms of $\cA$.

\medskip

To prove b) one needs a key   result which is the analogue in our
context of the quasi-invariance under diffeomorphisms of the smooth
measure class on a manifold, whose replacement in our case is given
by the Dixmier trace. This is shown in Proposition
\ref{abscontdirect} at the end of \S \ref{sectselfadj}. We then
prove in \S \ref{sectabscont} the required absolute continuity of
the spectral measure using a smearing argument. In \S
\ref{sectspecmul} we show the required inequality between the
multiplicity of the map $s_\alpha$ and the spectral multiplicity of
the $a_\alpha^j$.

\medskip

 To prove c) a new
strategy is required. Roughly one needs to know that the
multiplicity function of a tentative local coordinate system is
locally bounded while the information one obtains just by applying
the strategy outlined in \cite{CoSM} (and pursued in \cite{ReVa}) is
that it is a lower semicontinuous\footnote{the inverse image of
$]a,\infty]$ is open} integrable function. Typical examples of
Lebesgue negligible dense $G_\delta$ sets\footnote{countable
intersection of open sets} show that, as such, the situation is
hopeless. In order to solve this problem, one needs a local form of
the basic inequality (\cite{Co-action1}, \cite{Co-book}, Proposition
IV.3.14)
  giving an upper bound on the Voiculescu obstruction.
We prove this result in \S \ref{sectvoic}. This key result is
 combined with Voiculescu's Theorem (Theorem 4.5 of \cite{Voic1})
 and with the initial implicit function technique to
conclude the proof in \S \ref{sectloccoord}. Our main result can be
stated as follows (\cf Theorem \ref{regthm}):

\begin{thm}  Let $(\cA,\cH,D)$ be a  spectral triple, with  $\cA$ {\em commutative}, fulfilling the
first five \axioms of \cite{CoSM} (\cf \S \ref{prelem}) in a slightly
stronger form \ie we assume that:
\begin{itemize}
  \item The regularity holds for all $\cA$-endomorphisms of $\cap
  \Dom D^m$.
  \item The Hochschild cycle $c$ is  antisymmetric.
\end{itemize}
Then there exists a compact   oriented smooth manifold $X$ such that
$\cA$ is the algebra $C^\infty(X)$ of smooth functions on $X$.
\end{thm}

Moreover every compact   oriented smooth manifold  appears in this
spectral manner. Our next result is  the following variant (Theorem
\ref{fundthmspinc}):
\begin{thm}\label{withdimhyp}  Let $(\cA,\cH,D)$ be a   spectral triple with $\cA$
commutative, fulfilling the first five \axioms of \cite{CoSM} (\cf \S
\ref{prelem}) with the cycle $c$ antisymmetric. Assume that the
multiplicity of the action of $\cA''$ in $\cH$ is $2^{p/2}$. Then
there exists a smooth oriented compact (spin$^c$) manifold $X$ such
that $\cA=C^\infty(X)$.
\end{thm}
This multiplicity hypothesis is a weak form of the Poincar\'e
duality \axiom 6 of \cite{CoSM} and thus the above theorem can be
seen as the solution of the original problem formulated in
\cite{CoSM} and gives a characterization of spin$^c$ manifolds. It
follows from \cite{CoSM} (\cf \cite{FGV} for the proof) that the
operator $D$ is then a Dirac operator. The reality \axiom selects
spin manifolds among spin$^c$, and the spectral action (\cite{cc2})
selects the Levi-Civita connection.

 Finally we make a few remarks in \S
\ref{sectfinal}.  The first describes a different perspective on our main
result. As  explained many times, it is only because one drops commutativity
that variables with continuous range can coexist with infinitesimal variables
(which  only affect finitely many values larger than a given $\epsilon$). In
the classical formulation of variables, as maps from a set $X$ to the real
numbers, infinitesimal variables cannot coexist with continuous variables. The
formalism of quantum mechanics and the uniqueness of the separable infinite
dimensional Hilbert space cure this problem.  Using this formalism, variables
with continuous range (\ie self-adjoint operators with continuous spectrum)
coexist, in the same operator theoretic framework, with variables with
countable range, such as the infinitesimal ones (\ie compact operators). The
only new fact is that they do not commute. The content of Theorem
\ref{withdimhyp} can be expressed in a suggestive manner from this coexistence
between the continuum and the discrete. We fix the integer $p$ and
$N=2^{[p/2]}$ where $[p/2]$ is the integral part of $p/2$. The continuum will
only be used through its ``measure theoretic" content. This is captured by a
commutative von-Neumann algebra $M$ and, provided there is no atomic part in
$M$, this algebra is then unique (up to isomorphism). It is uniquely
represented in Hilbert space $\cH$ (which we fix once and for all, as a
universal stage) once the spectral multiplicity is fixed equal to  $N$. Thus
the pair $(M,\cH)$ is unique (up to isomorphism). Let us now consider
(separately first) an infinitesimal $ds$ \ie a self-adjoint compact operator in
$\cH$. Equivalently we can talk about its inverse $D$ which is unbounded and
self-adjoint. We assume that $ds$ is an infinitesimal of finite order
$\alpha=\frac 1p$. The information contained in the operator $ds$ is entirely
captured by a list of real numbers, namely the eigenvalues of $ds$ (with their
multiplicity). This list determines uniquely (up to isomorphism) the pair
$(\cH,D)$. Theorem \ref{withdimhyp} can now be restated as the birth of a
geometry  from the coexistence of $(M,\cH)$ with $(\cH,D)$. This coexistence is
encoded by a unitary isomorphism $F$ between the Hilbert space of the canonical
  pair $(M,\cH)$ and the Hilbert space of the canonical pair
  $(\cH,D)$. Thus the full information on the geometric
space is subdivided in two pieces
\begin{enumerate}
  \item The list of eigenvalues of $D$.
  \item The unitary $F$.
\end{enumerate}
We point out in  \S \ref{sectfinal} the analogy between these
parameters for geometry and the parameters of the Yukawa coupling of
the Standard Model (\cite{mc2}) which are encoded similarly by:
\begin{enumerate}
  \item The list of masses.
  \item The CKM matrix $C$.
\end{enumerate}
This analogy as well as the precise definition of the  corresponding unitary
invariant of Riemannian geometry will be dealt with in details in the companion
paper \cite{ckm1}.

\smallskip
 The second remark recalls a result of M. Hilsum on finite propagation
(\cf \cite{hilsum}). We then discuss briefly in \S \ref{sectfinal} the
variations dealing with real analytic manifolds, non-integral dimensions and
the non-commutative case.

\smallskip
We end with two appendices. In the first, \S \ref{regularity}, we
discuss equivalent formulations of the regularity condition. In the
second, \S \ref{appendix}, we recall the basic properties of the
Dixmier trace and its relation with the heat expansion.

\medskip

{\bf Acknowledgment.} The author is grateful to G. Landi, H. Moscovici,  G.
Skandalis and R. Wulkenhaar for useful comments.

\bigskip \section{Preliminaries}\label{prelem}

 Let us recall the \axioms
for {\it commutative} geometry as formulated in \cite{CoSM}.  We shall only use
the first five \axiomss.

 We let $(\cA,\cH,D)$ be a spectral triple, thus $\cH$ is a Hilbert space, $\cA$ an
involutive   algebra represented in $\cH$ and $D $ is a selfadjoint
operator in $\cH$. We assume that $\cA$ is {\em commutative}.
  We are given an integer $p$ which controls the
dimension of our space. The \axioms are:

  1) {\bf  Dimension}: {\it The  $n$-th characteristic value of the resolvent of $D$ is
$O(n^{-1/p})$.}

  2) {\bf Order one}: {\it $\left[ [D,f],g\right] = 0 \qquad  \forall \,
f,g \in \cA$.}

  We let
$\delta (T) = [ \vert D \vert ,T]$ be the commutator\footnote{The
domain of $\delta$ is the set of bounded operators $T$ with $T\Dom
|D| \subset \Dom |D|$ and $\delta(T)$ bounded.}  with the absolute
value of $D$:

  3) {\bf  Regularity}: {\it  For any $a\in \cA$ both $a$
and $[D,a]$ belong to the domain of $\delta^m$, for any integer
$m$.}

  We let $\pi_D $ be the linear map   given by
\begin{equation}\label{piD}
   \pi_D (a^0  \otimes a^1  \otimes \ldots  \otimes a^p) = a^0 [D,a^1]
\ldots [D,a^p] \qqq a^j \in \cA\,.
\end{equation}

  4)  {\bf  Orientability}: {\it  There exists a
Hochschild cycle $c \in Z_p (\cA ,\cA)$ such that $\pi_D (c)=1$ for
$p$ odd, while for $p$ even, $\pi_D (c) =\gamma$ satisfies}
$$
\gamma = \gamma^* \ , \ \gamma^2 = 1 \ , \ \gamma D =-D\gamma \, .
$$

  5) {\bf  Finiteness and absolute continuity}:{\it
Viewed as an $\cA$-module the space $\cH_{ \infty} = \,
  \cap_{m}^{} \, \hbox{\rm Domain} \, D^m$ is}
finite and projective. {\it Moreover the following equality defines
a hermitian structure $( \ |\ )$ on this module,}
\begin{equation}\label{absocont}
    \langle  \xi ,a \,\eta  \rangle = \cutint \, a (\xi |\eta) \, |D|^{-p}
\qquad  \forall \, a \in \cA \ , \  \forall \, \xi ,\eta \in \cH_{
\infty} \, .
\end{equation}

\bigskip

In other words the module can be written as $\cH_{ \infty}=e\cA^n$
with $e=e^*\in M_n(\cA)$ defining the Hermitian structure so that
\begin{equation}\label{hermstruc}
   (\xi|\eta)=\sum \xi_i^*\eta_i\in \cA\qqq \xi,\eta\in e\cA^n
\end{equation}
It follows from \axiom 4) and from\footnote{We shall not use this
result in an essential manner since one can just fix a choice of
Dixmier trace $\Tr_\omega$ throughout the proof.} (\cite{Co-book},
Theorem 8, IV.2.$\gamma$, and \cite{FGV}) that the operators $a \,
|D|^{-p} $, $a \in \cA$ are {\it measurable} (\cite{Co-book},
Definition 7, IV.2.$\beta$)
  so that the coefficient $\cutint \, a \, |D|^{-p}
$ of the logarithmic divergence of their trace is unambiguously
defined.

  It follows from \axiom 5) that the algebra $\cA$
is uniquely determined inside its weak closure $\cA''$ (which is
also the bicommutant of $\cA$ in $\cH$) by the equality
$$
\cA = \{ T\in \cA'' \ ; \ T\in \,   \cap_{m>0}^{} \, \Dom \delta^m
\} \, .
$$
This was stated without proof in \cite{CoSM} and we give the
argument below:

\begin{lem} \label{insideweak} The following conditions are equivalent for $T\in \cA''$:
\begin{enumerate}
  \item $T\in \cA$
   \item $[D,T]$ is bounded and both $T$ and $[D,T]$
   belong to the domain of $\delta^m$, for any integer
$m$
  \item $T$ belongs to the domain of $\delta^m$, for any integer
$m$
  \item $T \cH_{ \infty}\subset \cH_{ \infty}$
\end{enumerate}
\end{lem}

\proof Let us assume the fourth property. Then $T$ defines an
endomorphism of the finite projective module $\cH_{ \infty}=e\cA^n$
over $\cA$.
 As any endomorphism $T$ is of the form,
\begin{equation}\label{Top}
   T=e[a_{ij}]e \,, \  a_{ij}\in \cA
\end{equation}
\ie it is the compression of a matrix $a=[a_{ij}]\in M_n(\cA)$.

 Let us show that
since $T$ belongs to the weak closure of $\cA$ one can choose
$a_{ij}=x\delta_{ij}$ for some element $x$ of $\cA$. The norm
closure $A$ of $\cA$ in $\cL(\cH)$ is a commutative $C^*$-algebra,
$A=C(X)$ for some compact space $X$, and since $\cA$ is a subalgebra
of $\cL(\cH)$ it injects in $A$. The following equality
 \begin{equation}\label{measurelambdadef}
   \lambda(f)=\cutint f |D|^{-p}
 \qqq f\in A\,,
 \end{equation}
 defines a positive measure $\lambda$ on $X$. We let $\cE=eA^n$ be
 the induced finite projective module over $A$, which is
 intrinsically
 defined as $\cE=\cH_\infty\otimes_\cA A$. We let $S$ be the
 hermitian vector bundle on $X$ such that $\cE=C(X,S)$.
 By the absolute continuity relation
\eqref{absocont}, the representation of $A=C(X)$   in $\cH$ is
obtained from its action in $L^2(X,\lambda)$ by the tensor
product
\begin{equation}\label{induced}
   \cH=\cE\otimes_A L^2(X,\lambda)= e L^2(X,\lambda)^n=L^2(X,\lambda,S)\,.
\end{equation}
 This shows that the weak closure $\cA''=A''$ of $\cA$ in $\cH$ is
given by the diagonal action of $L^\infty(X,\lambda)$ in $ e
L^2(X,\lambda)^n$. Thus, since $T\in \cA''$, there exists $f\in
L^\infty(X,\lambda)$ such that $ T= ef $. It follows that the matrix
$eae$ belongs to the center of $eM_n(\cA)e$. This center is $e
(1\otimes \cA)$ and thus $T$ agrees with an element of $\cA$ which
proves the implication $4)\Rightarrow 1)$. To be more specific, and
for later use, let us give  a formula for an element $x\in \cA$ such
that $T=ex$ in term of the matrix elements $t_{ij}\in \cA$ of
$T=e[a_{ij}]e$. First the fact that $T$ belongs to the center of the
algebra $\cB=eM_n(\cA)e$ of endomorphisms of $\cH_\infty$ can be
seen directly since any such endomorphism $S$ is automatically
continuous in $\cH$ using \eqref{induced}. Thus since $T\in \cA''$
one has $ST=TS$.  Since $e$ is a self-adjoint idempotent and $\cA$
injects in $C(X)$ the element $\tau=\Tr(e)=\sum e_{jj}\in \cA$ is
determined by its image in $A$ which is just the function $\chi\in X
\mapsto \dim S_\chi\in \{0,1,\ldots,n\}$. This determines $n+1$
self-adjoint idempotents $p_j\in A$ by
\begin{equation}\label{projpj1}
   \tau=\Tr(e)=\sum j\,p_j\,\ \ \sum p_j=1 \,.
\end{equation}
To check that $p_j\in \cA$ it is enough to show that $p_j=P_j(\tau)$
where $P_j$ is a polynomial with
\begin{equation}\label{projpj2}
P_j(k)=0\qqq k\neq j\, ,\ 0\leq k\leq n\,, \ \ P_j(j)=1 \,.
\end{equation}
One then has the following formula\footnote{Note that $p_0=0$ because of the
faithfulness of the action of $\cA$ in Hilbert space together with \axiom 5.}
for $x$:
\begin{equation}\label{projpj3}
x=(\sum t_{ii})\sum_{j>0}\frac 1j \, p_j \in \cA\,.
\end{equation}
As $T$ belongs to the center of $eM_n(C(X))e$ one gets an equality $ T= ef $
for $f\in C(X)$ and working at every point $\chi\in X$ one then shows that
$T=ex$.

 The
implication $1)\Rightarrow 2)$ follows from the regularity, and $2)\Rightarrow
3)$ is immediate. To show the implication $3)\Rightarrow 4)$ one uses  the
definition of $\cH_{ \infty}$ as the intersection of domains of powers of $|D|$
and the implication $$T\in \Dom \,\delta^m\,, \ \xi\in \Dom \,|D|^m\Rightarrow
T\xi \in \Dom \,|D|^m$$ with the formula
\begin{equation}\label{commudel}
    |D|^m T\xi=\sum_{k=0}^m \left(
               \begin{array}{c}
                 m \\
                 k \\
               \end{array}
             \right)
\delta^k(T)\,|D|^{m-k}\xi \qqq \xi \in \Dom |D|^m
\end{equation}
which is proved by induction on $m$. More precisely this gives an estimate of
the norms but one has to care for the domains and proceed as follows. By
definition any $T\in \Dom \delta$ preserves the domain $\Dom |D|$ thus one gets
\eqref{commudel} for $m=1$. Let now $T\in \Dom \delta^2$  \ie $T\in \Dom
\delta$ and $\delta(T)\in \Dom \delta$. Let $\xi \in \Dom |D|^2$. Then since
$T\in \Dom \,\delta$ and $|D|\xi \in \Dom |D|$, one has $T|D|\xi \in \Dom |D|$.
One has $\delta(T)\xi=|D|T\xi-T|D|\xi$ where both terms make sense separately.
Since $\delta(T)\in \Dom \delta$ one has $\delta(T)\Dom |D|\subset \Dom |D|$.
Thus $\delta(T)\xi \in \Dom |D|$. Hence $|D|T\xi=\delta(T)\xi+T|D|\xi \in \Dom
|D|$ so that $T$ preserves $\Dom |D|^2$. Moreover one gets \eqref{commudel} for
$m=2$ as an equality valid on any vector $\xi \in \Dom |D|^2$. One can now
proceed by induction on $m$. We assume to have shown,
\begin{itemize}
  \item For $q\leq m$, $S\in \Dom \delta^q\Rightarrow S\Dom |D|^q \subset \Dom |D|^q
  $
  \item \eqref{commudel} holds for all $n\leq m$.
  \end{itemize}

For $T\in \Dom \delta^{m+1}$ and $\xi\in \Dom |D|^{m+1}$, one has
$\xi\in \Dom |D|^{m}$ and one can use the induction hypothesis to
get
$$
|D|^m T\xi=\sum_{k=0}^m \left(
               \begin{array}{c}
                 m \\
                 k \\
               \end{array}
             \right)
\delta^k(T)\,|D|^{m-k}\xi
$$
Let us show that $\delta^k(T)\,|D|^{m-k}\xi\in \Dom |D|$. One has
$\delta^k(T)\in \Dom \delta^{m+1-k}\subset \Dom \delta$ and
$|D|^{m-k}\xi\in \Dom |D|^{1+k}\subset \Dom |D| $ which gives the
result. Thus each term of the sum belongs to $\Dom |D| $ and one has
$$
|D|^{m+1} T\xi=\sum_{k=0}^m \left(
               \begin{array}{c}
                 m \\
                 k \\
               \end{array}
             \right)
|D|\delta^k(T)\,|D|^{m-k}\xi
$$
Moreover, as $\delta^k(T)\in  \Dom \delta$ and $|D|^{m-k}\xi\in \Dom
|D| $ one has
$$
|D|\delta^k(T)\,|D|^{m-k}\xi=\delta^{k+1}(T)\,|D|^{m-k}\xi+\delta^k(T)\,|D|^{m-k+1}\xi
$$
which gives \eqref{commudel} for $n+1$.
\endproof

\medskip

 This shows that the whole geometric data $(\cA ,\cH ,D)$ is in
fact uniquely determined by the triple $(\cA'' ,\cH ,D)$ where
$\cA''$ is a commutative von Neumann algebra.

  This also shows that $\cA$ is a pre-$C^*$-algebra, i.e.
  that it is stable under the holomorphic functional
calculus in the $C^*$-algebra norm closure of $\cA$, $A=
\overline{\cA}$. Since we assumed that $\cA$ was commutative, so is
$A$ and by Gelfand's theorem $A = C(X)$ is the algebra of continuous
complex valued functions on $X = \Sp  (A)$. We note finally that
characters $\chi$ of $\cA$ are automatically self-adjoint:
$\chi(a^*)=\overline \chi(a)$ since the spectrum of self-adjoint
elements of $\cA$ is real. Also they are automatically continuous
since  the $C^*$-norm is uniquely determined algebraically by
$$
\Vert a\Vert= \sup \{|\lambda|\,|\, a^*a-\lambda^2\notin \cA^{-1}\}
$$
thus they extend automatically to $A$ by continuity so that
$$
\Sp  A = \Sp  \cA \, .
$$

\medskip
 We shall now show that
$\cA$ is a Frechet  algebra \ie a complete locally convex algebra
whose topology is defined by the submultiplicative norms
\begin{equation}\label{submul}
p_k(xy)\leq p_k(x)p_k(y)\qqq x,y\in \cA
\end{equation}
 associated to the regularity condition, for instance by
 \begin{equation}\label{pk}
 p_k(x)=\Vert \rho_k(x)\Vert \,, \
 \rho_k(x)=
 \left(
   \begin{array}{cccc}
     x & \delta(x) & \ldots & \delta^k(x)/k! \\
     0 & x & \ldots & \ldots \\
     \ldots & \ldots & x & \delta(x) \\
     0 & \ldots & 0 & x \\
   \end{array}
 \right)
 \end{equation}
since $\rho_k$ is a representation of $\cA$.

\begin{prop} \label{generaltop0}
\begin{enumerate}
                                 \item The unbounded derivation $\delta$ is a
closed operator in $\cL(\cH)$.
                                 \item The algebra $\cA$ endowed with the norms $p_k$ is a Frechet algebra.
                                 \item The semi-norms $p_k([D,a])=p'_k(a)$ are continuous.
                               \end{enumerate}

\end{prop}

\proof (1) Let $G(|D|)$ be the graph of $|D|$. The graph of $\delta$ is
$$
G(\delta)=\{(T,S)\in \cL(\cH)^2\;|\; (T\xi,T\eta+S\xi)\in G(|D|) \qqq
(\xi,\eta)\in G(|D|)\}\,.
$$
It is therefore closed.

(2) Let us show that $\cA$ is complete. Let $a_n\in \cA$ be a
sequence which is a Cauchy sequence in any of the norms $p_k$. Then
$a_n\to T$ in norm, so that $T\in A\subset\cA''$. Since $\delta$ is
a closed operator one has $T\in \Dom \delta$ and
$\delta(a_n)\to\delta(T)$ in norm. By induction one gets, using the
closedness of $\delta$ that $T\in \Dom \delta^k$ and
$\delta^k(T)=\lim \delta^k(a_n)$. Thus $T\in \cap \Dom \delta^m$ and
by Lemma \ref{insideweak}, we get $T\in \cA$. Furthermore we also
have the norm convergence $\delta^k(T)=\lim \delta^k(a_n)$. This
shows that the $a_n$ converge to $T$ in the topology of the norms
$p_k$ and hence that $\cA$ is a Frechet space.

(3) Let us show that if we adjoin the semi-norms $p'_k$ to the
topology of $\cA$ we still get a complete space. The argument of the
proof of (1) only uses the closedness of the operator $|D|$ and thus
we get in the same way that the derivation $T\to d(T)=[D,T]$ with
domain $\Dom d=\{T\in \cL(\cH)\,|\, T\Dom D\subset \Dom D\,, \ \Vert
[D,T]\Vert<\infty\}$ is closed for the norm topology of $\cL(\cH)$.
Thus the above proof of completeness applies. The result then
follows from the open mapping Theorem (\cite{Rudin2} Corollary 2.12)
applied to the identity map from $\cA$ endowed with the topology of
the $p_k$, $p'_k$ to $\cA$ endowed with the topology of the
$p_k$.\endproof

In fact Lemma \ref{insideweak} shows that one has Sobolev estimates,
using finitely many generators $\eta_\mu$ of the $\cA$-module
$\cH_\infty$ to define the Sobolev norms on $\cA$ as
\begin{equation}\label{sobolevnorms}
    \Vert a\Vert_s^{\rm sobolev}=(\sum_\mu
    \Vert(1+D^2)^{s/2}a\eta_\mu \Vert^2)^{1/2} \qqq a\in \cA
\end{equation}
One has

\begin{prop} \label{generaltop}
\begin{enumerate}
  \item When endowed with the norms \eqref{sobolevnorms},
$\cA$ is a Frechet separable nuclear space.
  \item One has Sobolev estimates of the form
\begin{equation}\label{sobolevest}
    p_k(a)\leq c_k \Vert a\Vert_{s_k}^{\rm sobolev}\,, \ \
    p_k([D,a])\leq c'_k \Vert a\Vert_{s'_k}^{\rm sobolev}
    \qqq a\in \cA
\end{equation}
with $c_k<\infty$, $c'_k<\infty$ and  suitable sequences $s_k>0$,
$s'_k>0$.
  \item The spectrum $X=\Sp(\cA)$ is metrizable.
  \item  Any $T\in \End_{\cA}\cH_\infty$ is continuous in $\cH_\infty$ and
extends continuously to a bounded operator in $\cH$.
\item The algebraic isomorphism $\cH_\infty=e\cA^n$ is topological.
\item The map $(a,\xi)\mapsto a\xi$ and the $\cA$-valued inner
product are jointly continuous $\cA\times \cH_\infty\to \cH_\infty$
and $\cH_\infty \times \cH_\infty\to \cA$.
\end{enumerate}
\end{prop}

\proof 1) By construction the family \eqref{sobolevnorms} is an
increasing sequence of norms. Let us show that $\cA$ is complete.
Let $a_n$ be a sequence of elements of $\cA$ such that the vectors
$(1+D^2)^{s/2}a_n\eta_\mu$ converge for all $s$ (and all $\mu$). We
then obtain vectors
$$
\zeta_\mu = \lim a_n\eta_\mu\in \cH_\infty \qqq \mu\,,
$$
where the convergence holds in the topology of  $\cH_\infty$. Let
then $T$ be the operator given by
\begin{equation}\label{defofT}
   T\xi=\lim a_n\xi \qqq \xi \in \cH_\infty\,.
\end{equation}
It is well defined since one can write $\xi=\sum b^\mu\eta_\mu$ with
$b^\mu \in \cA$ which gives $a_n\xi=\sum b^\mu a_n\eta_\mu$ which
converges, in the topology of  $\cH_\infty$, to $\sum b^\mu
\zeta_\mu$ since the $b^\mu$ are continuous linear maps on
$\cH_\infty$ using \eqref{commudel} and regularity. Thus $T$ is a
linear map on $\cH_\infty$ and it commutes with $\cA$ \ie it is an
endomorphism of this finite projective module. Thus $T$ is of the
form \eqref{Top} and in particular it is bounded in $\cH$. Also
since endomorphisms of the finite projective module are
automatically continuous in $\cH$ they commute with $T$ using
\eqref{defofT}. Thus the argument of Lemma \ref{insideweak} shows
that $T\in \cA$. Moreover, since the convergence \eqref{defofT}
holds in the topology of  $\cH_\infty$, one has $a_n\to T$ in the
Sobolev topology and $\cA$ is complete in that topology. Thus $\cA$
is a Frechet space. It is by construction a closed subspace of the
sum of finitely many spaces $\cH_\infty$ each being a separable
nuclear space (of sequences of rapid decay). Thus it is a separable
nuclear space.

2) The identity map from the Frechet algebra $\cA$ with the norms
$p_k$ to the Frechet space $\cA$ with the Sobolev topology, is
continuous (using \eqref{commudel}) and surjective. Hence the open
mapping Theorem (\cite{Rudin2} Corollary 2.12) asserts that it is an
open mapping. This shows that the inverse map is continuous which
gives the required estimates for the norms $p_k$. The   result for
the semi-norms $p_k([D,a])$ follows from Proposition
\ref{generaltop0}.

3) Since $\cA$ is a Frechet separable nuclear space there is a
sequence $x_n\in \cA$ which is dense in any of the continuous norms
and in particular  using 2) in the $p_0$ norm. This shows that the
$C^*$-algebra $A$ is norm separable and hence that its spectrum is
metrizable.

4) By hypothesis $T$ being an endomorphism is of the form
\eqref{Top}. Using the inclusion $\cA\subset A=C(X)$ of $\cA$ in its
norm closure, we can   view $T$ as an endomorphism of the induced
$C^*$-module $\cE$ over $A$.   By \eqref{induced}, any element of
$\End_A(\cE)$ defines a bounded operator in $\cH$. This shows that
the graph of the operator $T$ in $\cH_\infty\times \cH_\infty$ is
closed and hence by the closed graph theorem that $T$ is continuous
in the Frechet topology of $\cH_\infty$.

5) The product $\cA\times \cA\to \cA$ is jointly continuous using
the submultiplicative norms $p_k$ of \eqref{pk}. This shows that
$e\cA^n$ is a closed subspace of $\cA^n$ and hence is complete.
Moreover the map $(a_j)\mapsto \sum a_j\xi_j$ for given $\xi_j\in
\cH_\infty$ is continuous from $\cA^n$ to $\cH_\infty$ using
\eqref{commudel}. Thus the open mapping theorem gives the result.

6) Follows from 5) and the joint continuity of the product
$\cA\times \cA\to \cA$.
\endproof

We end this section with the stability of $\cA$ under the smooth
functional calculus as shown in \cite{ReVa} Proposition 2.8. We
repeat the proof for convenience.

\begin{prop}\label{functcalc} Let $a_j=a_j^*$ be $n$ self-adjoint
elements of $\cA$ and $f\,:\,\R^n\mapsto \C$ be a smooth function
defined on a neighborhood of the joint spectrum of the $a_j$. Then
the element $f(a_1,\ldots ,a_n)\in A$ belongs to $\cA\subset A$.
\end{prop}

\proof Let us first show that for $a=a^*\in \cA$ one has for any
$k\in \N$,
\begin{equation}\label{deltakesti}
    \Vert \delta^k(e^{isa})\Vert =O(|s|^k) \,, \  |s|\to \infty\,.
\end{equation}
For $k=1$ one has
\begin{equation}\label{derdeltaforexp}
    \delta(e^{isa})=is\int_0^1 e^{itsa}\delta(a)e^{i(1-t)sa}dt
\end{equation}
which proves \eqref{deltakesti} for $k=1$. In general one has, with
$\beta_u(T)=e^{iusa}Te^{-iusa}$,
\begin{equation}\label{derdeltaforexp1}
   \frac{1}{n!} \delta^n(e^{isa}) e^{-isa}=\sum_{k_j> 0,\,\sum k_j=n} i^\ell s^\ell
    \int_{S_\ell}
    \beta_{u_1}(\frac{\delta^{k_1}(a)}{k_1!})\cdots \beta_{u_\ell}(\frac{\delta^{k_\ell}(a)}{k_\ell!})
    du
\end{equation}
where $S_\ell=\{(u_j)\,|\,0\leq u_1\leq\ldots\leq u_\ell\leq 1\}$ is
the standard simplex. This gives \eqref{deltakesti}. Now the joint
spectrum $K\subset \R^n$ of the $a_j$ is a compact subset and one
can extend $f$ to a smooth function with compact support $f\in
C_c^\infty(\R^n)$. The element $f(a_1,\ldots ,a_n)\in A$ is then
given by
\begin{equation}\label{eletfa}
    f(a_1,\ldots ,a_n)=(2\pi)^{-n}\int \hat f(s_1,\ldots ,s_n)\prod
    e^{is_ja_j}\prod ds_j
\end{equation}
where $\hat f$ is the Fourier transform of $f$ and is a Schwartz
function $\hat f\in \cS(\R^n)$. By \eqref{deltakesti} the integral
\eqref{eletfa} is convergent in any of the norms $p_k$ which define
the topology of $\cA$ and one gets $f(a_1,\ldots ,a_n)\in \cA$.
\endproof

\medskip

\medskip \section{Openness Lemma}\label{sectopen}

In this section, we use the standard implicit function theorem for
smooth maps $\R^p\to \R^p$ to obtain the openness of the tentative
local charts. We formulate the Lemma in a rather abstract manner
below and use it concretely in \S \ref{sectabscont} for the local
charts.

As above and in \cite{ReVa}, we let $\cA$ be a Frechet
pre-$C^*$-algebra.   We recall, for involutive algebras, the reality
condition which defines a $*$-derivation:
\begin{equation}\label{realder}
    \delta_0(a^*)=\delta_0(a)^* \qqq a\in \cA\,.
\end{equation}
We let $\Der \cA$ be the Lie algebra of continuous $*$-derivations
of $\cA$.

\begin{defn}\label{beingtame} Let $\cA$ be a Frechet pre-$C^*$-algebra.
A  continuous $*$-derivation $\delta_0\in \Der \cA$ {\em
exponentiates} iff one has a
  unique solution, depending continuously on $(t,a)\in \R\times \cA$, of
  the differential equation:
  \begin{equation}\label{diffequa}
    \partial_t\, y(t,a)=\delta_0(y(t,a))\,, \ y(0,a)=a\,.
\end{equation}
We say that $\cA$ is {\em \smooth} when any
  continuous $*$-derivation $\delta_0\in \Der \cA$  exponentiates.
\end{defn}

 We shall show in \S\S \ref{dissipsect},
\ref{sectselfadj} that in our context enough derivations
exponentiate but for clarity of the argument we shall first assume
that the algebra $\cA$ is \smooth. We refer to \cite{Hamilton} \S
I.3 for the discussion of differentiability in the context of
Frechet spaces. We just recall that a map $y:F\to G$ of Frechet
spaces is of class $C^1$ when the following directional derivative
exists and is a jointly continuous function of $(x,h)\in F\times F$,
\begin{equation}\label{direder}
    Dy(x,h)=\lim_{\epsilon\to 0}\frac 1\epsilon (y(x+\epsilon
    h)-y(x))
\end{equation}
The map is of class $C^n$ when the higher derivatives
$D^ky(x,h_1,\ldots,h_k)$ which are defined by iteration of
\eqref{direder} exist and are jointly continuous functions for
$k\leq n$. The map is smooth (or of class $C^\infty$) iff it is of
class $C^n$ for all $n$.

\begin{prop} One has for any $a,b\in \cA$,
\begin{equation}\label{diffequa1}
y(t,ab)=y(t,a)y(t,b)\,,\ y(t,a^*)=y(t,a)^*\,, \
y(t,a+b)=y(t,a)+y(t,b)\,.
\end{equation}
\begin{equation}\label{diffequa1.5}
y(t_1+t_2,a)=y(t_1,y(t_2,a))\,, \
y(t,\delta_0(a))=\delta_0(y(t,a))\,.
\end{equation}
Moreover $y(t,a)$ is a smooth function of $(t,a)$ with $n$-th
derivative given by
\begin{equation}\label{diffequan}
D^ny(t,a,s_1,h_1,\ldots,s_n,h_n)= \delta_0^n(y(t,a))\prod s_j+\sum_i
\delta_0^{n-1}y(t,h_i)\prod_{j\neq i} s_j\,.
\end{equation}
\end{prop}

\proof The two equalities \eqref{diffequa1} and \eqref{diffequa1.5}
follow from the uniqueness of the solution. To prove
\eqref{diffequan} we consider the Frechet spaces $F=\R\times \cA$
and $G=\cA$ and compute the first derivative $Dy$. One has
$$
y(t+\epsilon s, a +\epsilon h)-y(t,a)=y(t+\epsilon s, a
)-y(t,a)+\epsilon y(t+\epsilon s,  h)
$$
so that
$$
Dy(t,a,s,h)=s \delta_0(y(t,a))+y(t,h)
$$
Since by \eqref{diffequa1.5} one has
$\delta_0^k(y(t,a))=y(t,\delta_0^k(a))$ for all $k$, one gets
\eqref{diffequan} by induction on $n$.\endproof

 The Taylor expansion
at $(t,a)$ is thus of the form
$$
y(t+s,a+h)\sim \sum (\delta_0^k(y(t,a))s^k+
\delta_0^k(y(t,h))s^k)/k!
$$

\begin{lem} \label{openlem} Let $\cA$ be   commutative,
and $a=(a^j)$ be
  $p$ self-adjoint elements of $\cA$. Let $\chi$ be a
character of $\cA$. Assume that there exits $p$ derivations
$\delta_j\in \Der \cA$   such that
\begin{itemize}
  \item Each $\delta_j$ exponentiates.
  \item The determinant of the matrix $
\chi(\delta_j(a^k))$ does not vanish.
\end{itemize}

Then the image under $a$ of any neighborhood of $\chi$ in the
spectrum $\Sp(\cA)$ of $\cA$ contains a neighborhood of $a(\chi)$ in
$\R^p$.
\end{lem}

\proof By hypothesis the derivations $\delta_j\in \Der \cA$ can be
exponentiated to the corresponding one parameter groups $F^j(t)\in
\Aut(\cA)$ of automorphisms of $\cA$. Note that the flows $F^j$ do
not commute pairwise in general. We then define a map $h$ from
$\R^p$ to $\Sp(\cA)$ by
$$
h=\chi\circ \sigma\,, \ \sigma_{(t_1, \ldots,t_p)}=F_{t_1}^1\circ
\ldots \circ F_{t_p}^p
$$
which defines a character since $F^j(t)\in \Aut(\cA)$ by
\eqref{diffequa1}. The map $h$ is continuous since the topology of
$\Sp(\cA)$ is the weak topology and for any $a\in \cA$ the map
$(t_1, \ldots,t_p)\in \R^p\mapsto \sigma_{(t_1, \ldots,t_p)}(a)$ is
continuous using Definition \ref{beingtame}. The coordinates of the
map $\phi=a\circ h$, from $\R^p$ to $\R^p$, are given by
$$
\phi^k(t_1, \ldots,t_p)=h(t_1, \ldots,t_p)(a^k)=\chi\circ
F_{t_1}^1\circ \ldots \circ F_{t_p}^p(a^k)
$$
The map
$$
(t_1, \ldots,t_p)\in \R^p\mapsto F_{t_1}^1\circ \ldots \circ
F_{t_p}^p(a^k)
$$
is a smooth map from $\R^p$ to $\cA$. Indeed the maps $(t,a)\mapsto
F_t^j(a)$ are smooth, and composition of smooth maps are smooth (\cf
\cite{Hamilton} Theorem 3.6.4), while the above map is the
composition:
\begin{equation}\label{smoothcomp}
F^1\circ(\id\times F^2)\circ \cdots \circ (\id^{p-2}\times
F^{p-1})\circ (\id^{p-1}\times F^p(a^k))
\end{equation}
$$
\R^p \stackrel{F^p(a^k)}
     {\longrightarrow} \R^{p-1}\times \cA
     \stackrel{F^{p-1}}
     {\longrightarrow} \R^{p-2}\times \cA\longrightarrow \ldots \longrightarrow
     \R\times \cA  \stackrel{F^{1}}
     {\longrightarrow} \cA
$$

 Thus the map $\phi=a\circ h$, obtained by composition with the
 character $\chi$ which is linear and continuous and hence smooth,
  is a smooth map from $\R^p$ to $\R^p$.
 The image of $0\in \R^p$ is $a(\chi)$. The partial derivatives at
$0$ are
$$
(\partial_j \phi^k)(0)=\chi(\delta_j(a^k))
$$
thus we know from the hypothesis of the lemma, that the Jacobian
does not vanish at $0$. It then follows from the implicit function
theorem that the mapping $\phi=a\circ h$ maps by  a diffeomorphism a
suitable neighborhood of $0$ to a neighborhood of $a(\chi)$. In
particular the image under $a$ of a neighborhood $W$ of $\chi$
contains the image under $\phi$ of $h^{-1}(W)$ which, since $h$ is
continuous, is a neighborhood of $0\in\R^p$. This shows that the
image under $a$ of any neighborhood of $\chi$ in the spectrum
$\Sp(\cA)$ of $\cA$ contains a neighborhood of $a(\chi)$ in $\R^p$.
\endproof

The above proof yields the following more precise statement:

\begin{lem} \label{openlembis} Under the hypothesis of Lemma
\ref{openlem}, there exists a smooth family $\sigma_t\in \Aut(\cA)$,
$t\in \R^p$, a neighborhood $Z$ of $\chi$ in $X=\Sp(\cA)$ and a
neighborhood $W$ of $0\in \R^p$ such that, for any $\kappa\in Z$,
the map $t\mapsto a(\kappa\circ\sigma_t)$ is a  diffeomorphism,
depending continuously on $\kappa$, of $W$ with a neighborhood of
$a(\kappa)$ in $\R^p$.
\end{lem}

\proof Let as above
\begin{equation}\label{sigmatdef}
    \sigma_{(t_1, \ldots,t_p)}=F_{t_1}^1\circ
\ldots \circ F_{t_p}^p
\end{equation}
The map which to $\kappa\in X$ associates the map $\psi_\kappa$ from
$\R^p$ to $\R^p$ given by $\psi_\kappa(t)=a(\kappa\circ \sigma_t)$
yields by restriction  a continuous map $X\to C^\infty(K,\R^p)$
where $K$ is a closed ball centered at $0\in \R^p$. Indeed for each
$j$ the map $t\in K\mapsto \sigma_t(a^j)\in \cA$ is smooth by
\eqref{smoothcomp} and thus its partial derivatives
$\partial_t^\alpha\sigma_t(a^j)$ are elements of $\cA$ which depend
continuously of $t$. One has
$$
\partial_t^\alpha\psi^j_\kappa(t)=\kappa(\partial_t^\alpha\sigma_t(a^j))
$$
and thus the partial derivatives of $\psi_\kappa(t)$ are continuous
functions of $(\kappa,t)$. Since the determinant of the jacobian
$\chi(\delta_j(a^k))$ does not vanish, the result follows from the
implicit function Theorem (see \eg \cite{Hamilton} Theorem
5.2.3).\endproof

\medskip \section{Jacobian and openness of local charts}\label{sectcomm}

We first briefly recall first the well known properties of multiple
commutators
 which we need later.

\begin{defn} Let $T_j\in \cB$ be elements of a noncommutative
algebra $\cB$, one lets
$$
[T_1,T_2,\ldots,T_n]=\sum_\sigma
\epsilon(\sigma)\,T_{\sigma(1)}T_{\sigma(2)}\cdots T_{\sigma(n)}
$$
where $\sigma$ varies through all permutations of $\{1,\ldots,n\}$
and $\epsilon(\sigma)$ is its signature.
\end{defn}

We mention the following general properties

\begin{prop} \label{multicom} Let $T_j\in \cB$ be elements of a noncommutative
algebra $\cB$.
\begin{enumerate}
  \item[a)] For any permutation $\alpha$ of $\{1,\ldots,n\}$, one has
  $$
[T_{\alpha(1)},T_{\alpha(2)},\ldots,T_{\alpha(n)}]=\epsilon(\alpha)\,[T_1,T_2,\ldots,T_n]
  $$
  \item[b)] If two of the $T_j$ are equal one has
  $$
[T_1,T_2,\ldots,T_n]=0
$$
  \item[c)] Let $\cA\subset\cB$ be a commutative subalgebra and $\cA'\subset \cB$ its relative
  commutant in $\cB$. Let $a_k^j\in \cA$, $\gamma_j\in \cA'$. Then,
  with $T_k=\sum a_k^j \gamma_j$, one has
  \begin{equation}\label{multiple}
    [T_1,T_2,\ldots,T_n]=\Det((a_k^j))\,[\gamma_1,\gamma_2,\ldots,\gamma_n]
\end{equation}
\item[d)] The equality \eqref{multiple} extends to the case
of a rectangular matrix $a_k^j\in \cA$ as follows
\begin{equation}\label{multiple1}
    [T_1,T_2,\ldots,T_n]=\sum_F\Det((a_k^j(F)))\,[\gamma_1(F),\gamma_2(F),\ldots,\gamma_n(F)]
\end{equation}
where the sum is over all subsets $F\subset \{1,\ldots,m\}$ with $\#
F=n$, the matrix $a_k^j(F)$ is the restriction of $a_k^j$ to $j\in
F$ and the $\gamma_j(F)$ are the $\gamma_i$, $i\in F$, ordered with
increasing index in $F$.
\end{enumerate}
\end{prop}

\proof a) This follows from
$\epsilon(\sigma\circ\alpha)=\epsilon(\sigma)\epsilon(\alpha)$.

b) The permutation of the two indices is odd, but does not affect
the expression which must vanish.

c) One has
 \begin{equation}\label{summ}
    [T_1,T_2,\ldots,T_n]=\sum_{(j_k)} \prod_{k=1}^n
a_k^{j_k}\,[\gamma_{j_1},\gamma_{j_2},\ldots,\gamma_{j_n}]
 \end{equation}
 where, a priori, the $(j_k)$ is an arbitrary map from
$\{1,\ldots,n\}$ to $\{1,\ldots,n\}$. By the second statement of the
lemma, these terms vanish when two of the indices $j_k$ are equal.
Thus one can take the sum over permutations $(j_k)$ and one can use
the first statement of the lemma to rewrite the corresponding term
as
$$
[\gamma_{j_1},\gamma_{j_2},\ldots,\gamma_{j_n}]=\epsilon(j)\,[\gamma_1,\gamma_2,\ldots,\gamma_n]
$$
It follows that
$$
[T_1,T_2,\ldots,T_n]=\Det((a_k^j))\,[\gamma_1,\gamma_2,\ldots,\gamma_n]
  $$
  d) One decomposes the sum \eqref{summ} according to the range $F$
  of the injection $j$ from $\{1,\ldots,n\}$ to $\{1,\ldots,m\}$.
\endproof

Let us now go back to spectral triples $(\cA,\cH,D)$ fulfilling the
five \axioms of \S \ref{prelem}.

\begin{lem} \label{invariantintegration0} Let $\cB$ be the algebra of endomorphisms of $\cH_\infty$. One has a finite decomposition
\begin{equation}\label{dirac}
   [D,a]=\sum \delta_j(a) \gamma_j \qqq a\in \cA\,,
\end{equation}
where $\gamma_j\in \cB$ and the $\delta_j$ are derivations of the
form
\begin{equation}\label{formforder}
    \delta_j(a)=i(\xi_j|[D,a]\xi_j)\qqq a \in \cA
\end{equation}
for some $\xi_j\in \cH_\infty$.
\end{lem}

\proof First $[D,a]\cH_\infty\subset \cH_\infty$ using regularity
and \eqref{commudel}. Thus the order one condition shows that
$[D,a]\in \cB$. One has $\cH_\infty=e\cA^n$, $\cB=eM_n(\cA)e$ for a
self-adjoint idempotent $e\in M_n(\cA)$. Thus every element $T\in
\cB$ can be written uniquely, as any element of $M_n(\cA)$ in the
form
$$
T=\sum a_{k\ell}\varepsilon_{k\ell}\,, \ \ a_{k\ell}\in \cA
$$
in terms of the matrix units $\varepsilon_{ij}$. The coefficients
$a_{k\ell}\in \cA$ are uniquely determined, using the elements
$\eta_k=e\zeta_k\in \cH_\infty$ where $\zeta_k\in \cA^n$ is the
element all of whose components vanish except the $k$-th one which
is equal to $1$. Using the $\cA$-valued inner product, one has
\begin{equation}\label{formforder0}
a_{k\ell}=(\eta_k|T\eta_\ell)=L_{k\ell}(T)\qqq k,\ell\,.
\end{equation}
Moreover one has, since $T=eTe$ and the $a_{k\ell}$ commute with
$e$:
\begin{equation}\label{sumofder}
T=\sum a_{k\ell}e\varepsilon_{k\ell}e\,.
\end{equation}
 One has
$L_{k\ell}(aT)=a L_{k\ell}(T)$ for any $a\in \cA$. Applying this to
$T=[D,b]$, the maps $a \mapsto L_{ij}([D,a])$ give derivations of
$\cA$. They are not self-adjoint but can be decomposed as linear
combinations of self-adjoint derivations, which, using
\eqref{sumofder}, gives the required formula \eqref{dirac}. More
precisely, the derivations $\delta_j$ can be written using the
$\cA$-valued inner product on $\cH_\infty$ in the form
\eqref{formforder} for some $\xi_j\in \cH_\infty$ (with $i$ to
ensure self-adjointness). Indeed  one obtains \eqref{formforder}
applying to \eqref{formforder0} the polarization identity:
\begin{equation}\label{polarization}
   2   (\xi|T\eta)=((\xi+\eta)|T(\xi+\eta))-(\xi|T\xi)-(\eta|T\eta)
\end{equation}
$$
-i\left(
   ((\xi+i\eta)|T(\xi+i\eta))-(\xi|T\xi)-(i\eta|Ti\eta)\right)
   $$
 In particular, using
Proposition \ref{generaltop}, the $\delta_j$ are continuous.
\endproof

 By hypothesis the cycle $c$ is of the form:
\begin{equation}\label{cyclec}
 c=\sum_\alpha a_\alpha^0\,\omega_\alpha \,, \ \ \omega_\alpha
=\sum_\beta \epsilon(\beta) 1\otimes a_\alpha^{\beta(1)}\otimes
\cdots \otimes a_\alpha^{\beta(p)}
\end{equation}
where one can assume that the $a_\alpha^\mu$ are self-adjoint for
$\mu>0$.  We define the conditional expectation $E_{\cA}\,:\,
\End_{\cA}(\cH_\infty)\to \cA$, using the projections $p_j$ of
\eqref{projpj1},
\begin{equation} \label{tracemap}
 E_{\cA}(T)=\sum_{j>0}\frac 1j p_j\sum T_{kk}\qqq T=(T_{k\ell})\in
eM_n(\cA)e
\end{equation}
using the identification $\cH_\infty=e\cA^n$. We   obtain a
self-adjoint $\rho_\alpha\in \cA$ given\footnote{there is no
$\gamma$ in the odd case} by
 \begin{equation} \label{omegaalpha}
\rho_\alpha =i^{\frac{p(p+1)}{2}}\,E_{\cA}(\gamma \sum_\beta
\epsilon(\beta) [D,a_\alpha^{\beta(1)}] \cdots
[D,a_\alpha^{\beta(p)}])\,.
\end{equation}
 One lets
\begin{equation} \label{calpha}
C_\alpha=\{x\in X\,|\,\rho_\alpha(x)=0\}
\end{equation}
and $U_\alpha=C_\alpha^c$ be its complement \ie the open set where
$\rho_\alpha$ does not vanish.

\begin{lem} \label{coverlem} The $U_\alpha$ form an open cover of
$X=\Sp(\cA)$.

Each $U_\alpha$ is  the disjoint union of the two open subsets
$U_\alpha^\pm$ corresponding to the sign of $\rho_\alpha$,
$$
\pm\rho_\alpha(x) >0 \qqq x\in U_\alpha^\pm\, .
$$
\end{lem}

\proof It is enough to show that any $x\in X$ belongs to some
$U_\alpha$. One has $\pi_D(c)=\gamma$, so that by \eqref{cyclec}
$$
\gamma \sum_\alpha a_\alpha^0\,\sum_\beta \epsilon(\beta)
[D,a_\alpha^{\beta(1)}] \cdots [D,a_\alpha^{\beta(p)}]=1\,.
$$
By \eqref{omegaalpha} and the conditional expectation module
property $E_{\cA}(aT)=aE_{\cA}(T)$,
$$
i^{-\frac{p(p+1)}{2}}\,\sum_\alpha a_\alpha^0\,\rho_\alpha=1
$$
and $\rho_\alpha(x)\neq 0$ for some $\alpha$. The second statement
follows since $\rho_\alpha$ is a non-vanishing real valued function
on $U_\alpha$.
\endproof
  We let  $s_\alpha$ be the natural
continuous  map from $X$ to $\R^p$ given by
\begin{equation}\label{defnsalpha}
\chi\in \Sp(\cA)\to(\chi(a_\alpha^j))\in \R^p\,.
\end{equation}

\begin{lem} \label{coropen} Assume  that derivations
of the form \eqref{formforder} exponentiate. Let $\chi \in
U_\alpha$.
\begin{itemize}
  \item  There exists $p$ derivations
  $\delta_j\in \Der(\cA)$ such that $\chi(\Det((\delta_j(a_\alpha^k))))\neq
  0$
  \item The map $a_\alpha$ from $U_\alpha$ to $\R^p$ is open.
  \item There exists a smooth family $\sigma_t\in \Aut(\cA)$,
$t\in \R^p$, a neighborhood $Z$ of $\chi$ in $X=\Sp(\cA)$ and a
neighborhood $W$ of $0\in \R^p$ such that, for any $\kappa\in Z$,
the map $t\mapsto s_\alpha(\kappa\circ\sigma_t)$ is a
diffeomorphism, depending continuously on $\kappa$, of $W$ with a
neighborhood of $a(\kappa)$ in $\R^p$.
\end{itemize}
\end{lem}

\proof We let, as above,  $\cB$ be the algebra of endomorphisms of
the $\cA$-module $\cH_\infty$. It contains $\cA\subset\cB$ as a
subalgebra of its center. By Lemma \ref{invariantintegration0}, one
has derivations $\delta_j\in \Der(\cA)$ of the form
\eqref{formforder}, such that the formula \eqref{dirac} holds:
$$
[D,a]=\sum_1^m \delta_j(a) \gamma_j\qqq a \in \cA\,.
$$
By hypothesis we have  $  \rho_\alpha(\chi) \neq 0 $. Thus,  the
following endomorphism of the $\cA$-module $\cH_\infty$ does not
vanish,
$$
[[D,a_\alpha^1],[D,a_\alpha^2],\ldots,[D,a_\alpha^p]](\chi)\neq 0
\qqq \chi \in U_\alpha\,.
$$
 It thus
follows, from \eqref{multiple1} of proposition \ref{multicom}, that
for $\chi \in U_\alpha$ one can find $p$ elements $\delta_j\in
\Der(\cA)$ among the above $\delta_j$ such that:
\begin{equation}\label{nonzerojac}
\chi(\Det((\delta_j(a_\alpha^k))))\neq 0\,.
\end{equation}
Now let $V\subset U_\alpha$ be open. To show that $s_\alpha(V)$ is
open one needs to show that, for any $\chi\in V$, $s_\alpha(V)$
contains a neighborhood of $s_\alpha(\chi)$. But $V$ is a
neighborhood of $\chi$ in $\Sp(\cA)$ and the hypothesis of Lemma
\ref{openlem} is fulfilled so that this lemma shows that
$s_\alpha(V)$ contains a neighborhood of $s_\alpha(\chi)$. The third
statement follows from Lemma \ref{openlembis}.
\endproof

\begin{figure}
\begin{center}
\includegraphics[scale=0.7]{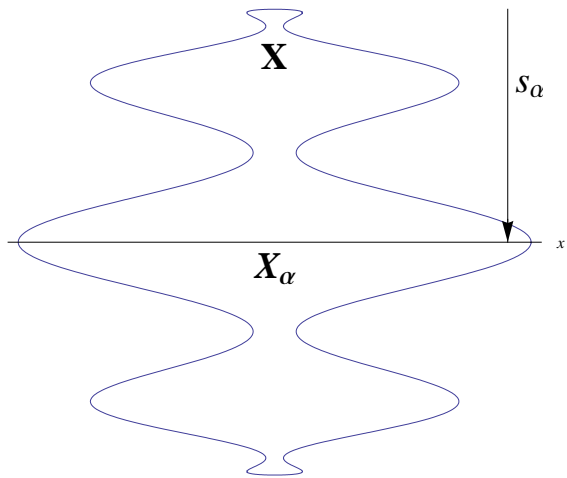}
\end{center}
\caption{The map $s_\alpha$ from $X$ to $X_\alpha$ \label{mapproj}}
\end{figure}

\medskip
  \section{Dissipative
derivations}\label{dissipsect}

We assumed in the above discussion that the algebra $\cA$ is
\smooth. It is of course desirable to remove this hypothesis, and
this will be done in this section and the next one. We need a form
of existence and uniqueness for solutions of linear differential
equations with values in a Frechet space $E$. Simple examples show
that in that generality one has neither existence nor uniqueness.
For failure of existence just let $\cA=C^\infty([0,1])$,
$\delta_0=\partial_x$. For failure of uniqueness,   let $E$ be the
space of sequences $x_n\in \R$ for $n\geq 1$ and take the shift
operator $S$, then the equation
$$
\partial_t x=Sx\,, \ \ (Sx)_n=x_{n+1}\,,
$$
has no uniqueness of solutions. Indeed we can take for $x_1(t)$ any
smooth function which is flat at $t=0$, \ie $\partial_t^n x_1(0)=0$,
and then define by induction $x_{n+1}(t)=\partial_t x_n(t)$ so that
$\partial_t x=Sx$ holds and the initial condition $x(0)=0$ does not
imply uniqueness.

In our case we need to know that any derivation $\delta\in \Der \cA$
can be exponentiated, \ie that one has existence and uniqueness for
the differential equation
$$
\frac{dy(t)}{dt}=\delta(y(t)) \,, \ \ y(t)\in \cA\,.
$$
It is only the compactness of $X$ that ensures this, and also the
fact that one is dealing with a {\em real} vector field.  This means
that we first need to make sure that the derivation exponentiates
  at the level of the $C^*$-algebra as discussed in \cite{Bratteli}.

One step towards this would be to show directly the following
corollary of \smoothness:
\begin{lem} \label{diss} Assume that the derivations of the form
the form \eqref{formforder} exponentiate, then, for any $h=h^*\in \cA$, the
commutator $[D,h]$ vanishes where\footnote{This makes sense since $[D,h]$
commutes with $h$} $h$ reaches its maximum. Conversely if this property holds
the derivations $\pm\delta_j$ of the form \eqref{formforder} are dissipative
(\cf \cite{Bratteli}, Definition 1.4.6), \ie
\begin{equation}\label{dissdefn}
  \Vert x+\lambda\delta_j(x)\Vert \geq \Vert x\Vert \qqq x\in
\cA\,, \ \lambda \in \R\,.
\end{equation}
\end{lem}

\proof One has, by \eqref{dirac}, $[D,h]=\sum \delta_j(h)\gamma_j$
where $\delta_j\in \Der\cA$. Thus it is enough to show that
$\delta_j(h)(\chi)=0$ where $h=h^*\in \cA$ reaches its maximum at
$\chi$. This follows from the existence of $e^{t\delta_j}\in
\Aut(\cA)$ using the differentiable function
$f(t)=\chi(e^{t\delta_j}(h))$ which has a maximum at $t=0$ and hence
vanishing derivative.

Conversely, the derivations $\delta_j$ are of the form
\eqref{formforder} \ie $\delta_j(h)=i(\xi|[D,h]\xi)$. Thus the
vanishing of $[D,h](\chi)$ where $h=h^*\in \cA$ reaches its maximum,
ensures that
$\delta_j(h)(\chi)=i\langle\xi(\chi),[D,h](\chi)\xi(\chi)\rangle=0$
also vanishes. Thus one has
$$
\Vert h+\lambda\delta_j(h)\Vert \geq \Vert h\Vert \qqq h=h^*\in
\cA\,, \ \lambda \in \R\,,
$$
since for a character $\chi$ of $\cA$ with $\chi(\pm h)=\Vert h\Vert$ one has
$\chi(\pm (h +\lambda\delta_j(h)))=\Vert h\Vert$. In the complex case \ie for
an arbitrary $x\in \cA$, let $\psi$ be a state on $A\supset \cA$, such that
$|\psi(x)|=\Vert x\Vert$. Replacing $x\mapsto u x$ for $u\in \C$, $|u|=1$, one
can assume that $\psi(x)>0$. Then writing $x=h+ik$ with $h=h^*$ and $k=k^*$,
one has $\psi(x)=\psi(h)=\Vert h\Vert$ so that $\psi(\delta_j(h))=0$ from the
above discussion. Then one has, for $\lambda\in \R$,
$$
\psi(x+\lambda\delta_j(x))=\psi(h)+i\lambda\psi(\delta_j(k))
$$
and $|\psi(x+\lambda\delta_j(x))|\geq \psi(h)=\Vert x\Vert$.
\endproof

Note that the commutativity of $[D,h]$  with $h$ and the
self-adjointness of $D$ do not suffice to entail the conclusion of
Lemma \ref{diss}. This can be seen with the following spectral
triple:
\begin{equation} \label{bdcond0}
\cA=C^\infty([0,1])\,, \ \ \cH=L^2([0,1])\otimes \C^2\,, \ \
D=\left(
  \begin{array}{cc}
    0 & \partial_x\\
    -\partial_x & 0 \\
  \end{array}
\right)
\end{equation}
with the boundary condition
\begin{equation} \label{bdcond}
\Dom D=\{\xi= \left(
  \begin{array}{c}
    \xi_1 \\
    \xi_2 \\
  \end{array}
\right)\;|\; \xi_1(0)=0\,,\ \ \xi_2(1)=0\}\,.
\end{equation}
For any $h\in \cA$ one has $[D,h]=\partial_x h\gamma_1$,
$$\gamma_1= \left(
  \begin{array}{cc}
    0 & 1 \\
    -1 & 0 \\
  \end{array}
\right)$$ so that $[D,h]$ commutes with $h$. For $h(x)=x$ the
maximum is at $x=1$ and $[D,h]$ does not vanish there.
 This example shows that the hypothesis of \smoothness
$\,$of the algebra $\cA$ appears at first sight as essential.
However, in this example,   \axiom 5  fails since the boundary
condition \eqref{bdcond} does not yield a finite projective
submodule of $C^\infty([0,1])\otimes \C^2$ over
$\cA=C^\infty([0,1])$. Also $[D,h]=\partial_x h\gamma_1$ does  not
preserve the domain of $D$ which is the same as the domain of $|D|$
thus {\em regularity fails}. Orientability also fails in this
example. We shall now show that   regularity allows in fact to
obtain the required dissipativity.

Let us consider the one parameter group of automorphisms of
$\cL(\cH)$ given by
\begin{equation}\label{alphat}
    \alpha_t(T)=e^{itD}Te^{-itD}\qqq t\in \R
\end{equation}

\begin{lem} \label{tentative} Let $T$ preserve $\Dom D$ and $[D,T]$
be bounded. Then the function $t\mapsto \alpha_t(T)$ is norm
continuous,
\begin{equation}\label{lip}
    \Vert \alpha_s(T)-\alpha_t(T)\Vert \leq |s-t|\, \Vert [D,T]\Vert
\end{equation}
and when $s\to 0$ the difference quotient
\begin{equation}\label{lip1}
\frac{\alpha_s(T)-T}{s}=\frac is \int_0^s \alpha_t([D,T])dt
\end{equation}
converges to $i[D,T]$ in the strong topology.
\end{lem}

\proof Let $\xi \in \Dom D$. Then $\frac 1s(e^{isD}-1)\xi\to i D\xi$
(in norm) when $s\to 0$. Thus using
$$\frac
1s(e^{isD}Te^{-isD}-T)\xi =\frac 1s e^{isD}T(e^{-isD}-1)\xi+\frac
1s(e^{isD}-1)T \xi$$ one gets that (in norm)  when $s\to 0$
$$
\frac 1s(e^{isD}Te^{-isD}-T)\xi\to i[D,T]\xi
$$
 Thus $t\mapsto \alpha_t(T)\xi$ is of class $C^1$. Its derivative is $t\mapsto i\alpha_t([D,T])\xi$. Thus
\begin{equation}\label{lip2}
(\alpha_s(T)-\alpha_t(T))\xi=i\int_t^s \alpha_u([D,T])\xi du
\end{equation}
holds for all $\xi \in \Dom D$ and hence all $\xi\in \cH$ since the
map $u\mapsto \alpha_u([D,T])\xi$ is continuous, as follows from the
continuity of $s\mapsto e^{isD}\eta$ for any $\eta\in \cH$. Both
statements follow.
\endproof

We can now consider the $C^*$-algebra $C$ generated by the
$\alpha_s(h)$ for $h=h^*\in \cA$ as above. It is norm separable and
the $\alpha_s\in\Aut(C)$ form a norm continuous one parameter group.
To try and prove that  $[D,h]$ vanishes where $h=h^*\geq 0$ reaches
its maximum, one considers a state $\phi$ on $C$ such that
$\phi(h)=\Vert h\Vert$. It is obtained by extension using the
inclusion $C^*(h)\subset C$. The function
$$
f(s)=\phi(\alpha_s(h))
$$
is a Lipschitz function and reaches its maximum: $\Vert h\Vert$ at
$s=0$. Thus if one could assert that the derivative at $s=0$ is
given by $\phi([D,h])$, one would get the vanishing $\phi([D,h])=0$.
The problem is that $\alpha_u([D,h])$ is not in general a norm
continuous function of $u$ and thus the differentiability only holds
in the strong topology but not in the norm topology.

 Things are easier with $|D|$ since the regularity
conditions ensures that the map
\begin{equation}\label{geod}
t\to \gamma_t(a)=e^{it|D|}ae^{-it|D|}
\end{equation}
is in fact of class $C^\infty$ in the norm topology (\cf Lemma \ref{geodcinfty}
of \S \ref{regularity}). Moreover the following Lemma shows that it is enough
to show the vanishing of $[[D^2,a],b]$ at $\chi$ for all $b\in \cA$ to get the
vanishing of $[D,a]$ at $\chi$.

\begin{lem} \label{fromdsquared} Let $h=h^*\in
\cA$ and $\chi\in \Sp(\cA)$. If $[[D^2,h],h]$ vanishes at $\chi$
then $[D,h]$ vanishes at $\chi$.
\end{lem}

\proof
\begin{equation}\label{double}
[[D^2,h],h]=2[D,h]^2
\end{equation}
 using the order one condition.
\endproof

Note that \eqref{double} shows that $[D,h]^2$ and hence $|[D,h]|$
only depends upon $D^2$ and hence $|D|$ and not upon the phase of
the polar decomposition of $D$. This comes from the order one
condition.  Moreover one has the following vanishing of $[|D|,h]$
where $h=h^*\in \cA$ reaches its maximum.

\begin{lem} \label{diss1} For any $h=h^*\in
\cA$, $h\geq 0$, reaching its maximum at $\chi\in \Sp(\cA)$ and any
sequence $b_n\in \cA$, $\Vert b_n\Vert\leq 1$, with support tending
to $\{\chi\}$, one has
\begin{equation}\label{conver}
   \Vert b_n^* [|D|,h] b_n\Vert\to 0
\end{equation}
\end{lem}

\proof Let $\xi_n\in \cH$ be unit vectors with support tending to
$\{\chi\}$. Then consider any limit state on $\cL(\cH)$:
\begin{equation}\label{conver0}
\eta(T)=\lim_\omega \langle \xi_n,T\xi_n\rangle
\end{equation}
One has $\eta(h)=h(\chi)$ since $h$ is a continuous function on
$X=\Sp(\cA)$. Thus $\eta(h)=\Vert h\Vert$. When applied to $|D|$
instead of $D$, Lemma \ref{lip} shows that both $\gamma_s(h)$ and
$\gamma_s(\delta(h))$ are Lipschitz functions of $s$, while
\begin{equation}\label{lip1bis}
\frac{\gamma_s(h)-h}{s}=\frac 1s \int_0^s \gamma_t(i\delta(h))dt\to
i\delta(h)
\end{equation}
so that $\gamma_s(h)$ is of class $C^1$ in norm. It follows that the
function $s\mapsto \eta(\gamma_s(h))$ is of class $C^1$. It is
maximal for $s=0$ and hence its derivative vanishes so that
$\eta(\delta(h))=0$. Thus
\begin{equation}\label{conver1}
\lim_\omega \langle \xi_n,\delta(h)\xi_n\rangle=0
\end{equation}
and this continues to hold for any bounded sequence $\xi_n\in \cH$
 with support tending to $\{\chi\}$. Now, let $b_n$
be as in the Lemma, then if \eqref{conver} does not hold, one can
find a subsequence $n_k$ with $\Vert b_{n_k}^* [|D|,h] b_{n_k}\Vert
\geq \epsilon >0$ for all $k$. Using polarization
\eqref{polarization}, one gets
 unit vectors $\xi'_k\in \cH$ such that
$$
 \vert\langle\xi'_k,b_{n_k}^* [|D|,h] b_{n_k}\xi'_k\rangle\vert\geq \epsilon' >0
$$
which contradicts \eqref{conver1}, for $\xi_k=b_{n_k}\xi'_k$.
\endproof

 We shall use the analogue in our context of the notion of symbol for
pseudodifferential operators. The symbol of $T$ can be viewed as a
weak limit of the conjugate operators of the form
$$
\tau^{-k}\,e^{i\tau\phi}Te^{-i\tau\phi}\,, \ \tau\to \infty
$$
where the integer $k$ is the order of $T$. For instance the symbol
of $D$ is given by $-i[D,\phi]$ since the order one condition gives
\begin{equation}\label{symofD}
\tau^{-1}\,e^{i\tau\phi}De^{-i\tau\phi}=-i[D,\phi]+\tau^{-1}\,D
\end{equation}
 One
expects the symbol of $D^2$ to be of the form
\begin{equation}\label{symofD2}
\lim_{\tau\to\infty}\frac{1}{\tau^2}e^{i \tau \phi}D^2e^{-i \tau
\phi}\xi=-[D,\phi]^2\xi \qqq \xi \in \Dom D^2
\end{equation}
 This is obtained by squaring \eqref{symofD} but one needs to know
 that
  $\Dom D$ is invariant under $[D,\phi]$ to control the
term $D\,[D,\phi]$. This is insured by regularity.

\begin{rem}\label{noregexample}{\rm   In the example \eqref{bdcond0} considered above, $[D,\phi]$
does not map $\Dom D^2$ to $\Dom D$ so that $D\,[D,\phi]\xi$ does
not make sense in that case. In fact  regularity fails, and $\Dom
D^2$ is not invariant under $\phi$, unless $[D,\phi]$ vanishes on
the boundary. To see this note that the boundary condition for $D^2$
is
\begin{equation} \label{bdcondbis}
\xi= \left(
  \begin{array}{c}
    \xi_1 \\
    \xi_2 \\
  \end{array}
\right)\in \Dom D^2 \Leftrightarrow \xi_1(0)=0\,,\
\partial_x\xi_2(0)=0  \,,\                     \xi_2(1)=0\,,
\partial_x\xi_1(1)=0  \,.
\end{equation}
which contains the Neumann condition $\partial_x\xi_2(0)=0$ while
$\xi_2(0)$ is arbitrary. Thus
$\partial_x\phi\xi_2(0)=\partial_x\phi(0)\xi_2(0)$ vanishes only
when $\partial_x\phi(0)=0$.}\end{rem}

In the case of an operator of order $0$ there is no power of $\tau$
and one deals with a bounded family so that one can expect the limit
to be a weak limit. We need to guess the symbol of $[|D|,h]$. We
expect that if we choose $\phi=h$, this symbol will just be
$i|[D,h]|\in {\rm End}\,S$. The symbol of $[D^2,h]$ is  (using
$[D^2,h]=D[D,h]+[D,h]D$)
\begin{equation}\label{symbd2}
    \tau^{-1}\,e^{i\tau\phi}[D^2,h]e^{-i\tau\phi}
    =-i([D,\phi][D,h]+[D,h][D,\phi])+\tau^{-1}\,[D^2,h]
\end{equation}
To see why we should expect the symbol of $[|D|,h]$ for $\phi=h$ to
just be $i|[D,h]|$ we have:
\begin{lem}\label{symbroot}
Assume that when $\tau\to \infty$ the following limit holds in the
strong topology:
$$
\lim_{\tau\to \infty} e^{i\tau h}[|D|,h]e^{-i\tau h}=T \,.
$$
Then one has, with strong convergence on $\Dom D$:
$$
\lim_{\tau\to \infty} \tau^{-1}\,e^{i\tau h}|D|e^{-i\tau h}=-iT \,.
$$
\end{lem}

\proof One defines a one parameter group $\beta_u$ of automorphisms:
\begin{equation}\label{betau}
   \beta_u(Y)=e^{iu h}Ye^{-iu h}
\end{equation}
One has, at the formal level, $\frac{d}{du}\beta_u(Y)=-i\beta_u([Y,h])$. Taking
$Y=|D|$ one gets, with the notations of the Lemma, and using regularity
\begin{equation}\label{cesaro0}
e^{iu h}[|D|,h]e^{-iu h}\xi=i \frac{d}{du}\beta_u(|D|)\xi \qqq \xi
\in \Dom |D|
\end{equation}
which gives
\begin{equation}\label{cesaro0bis}
\int_0^\tau e^{iuh}[|D|,h]e^{-iuh}\xi du=i (\beta_\tau(|D|)-|D|)\xi
\qqq \xi \in \Dom |D|
\end{equation}
Note that this equality continues to hold for any $\xi \in \cH$
since $\beta_\tau(|D|)-|D|$ is a bounded operator.
 Now $e^{iuh}[|D|,h]e^{-iuh}$ is uniformly
bounded and converges strongly by hypothesis to $T$. Thus one has,
for the Cesaro mean:
$$
\lim_{\tau\to \infty} \tau^{-1}\,\int_0^\tau
e^{iuh}[|D|,h]e^{-iuh}\xi \,du= T\xi \qqq \xi \in \cH
$$
which gives the result since one controls $ \tau^{-1}\,|D|\xi\to 0$
for $\xi \in\Dom D$.\endproof

Moreover we expect the symbol map to be a morphism so that the
symbol of $|D|$ is given by the absolute value of the symbol of $D$
\ie by $|[D,h]|$. In fact we do not need to prove the converse of
Lemma \ref{symbroot} since we can use the regularization by Cesaro
mean to compose the states   $\eta$, with weak limits of
$C_\tau(T)$,
\begin{equation}\label{cesaro}
C_\tau(T)= \tau^{-1}\,\int_0^\tau  \beta_u(T)du
\end{equation}
\begin{lem}\label{symbroot1} With $h$ as above one has
\begin{enumerate}
  \item $C_\tau$ is a completely positive map from $\cL(\cH)$ to itself and $C_\tau(1)=1$
  \item $C_\tau(aTb)=aC_\tau(T)b$ for all $a,b\in \cA$.
  \item $C_\tau([|D|,h])=\frac i\tau\,(e^{i\tau h}|D|e^{-i\tau h}-|D|) $
\end{enumerate}
\end{lem}

\proof The first two statements follow from \eqref{cesaro} using the
commutativity of $\cA$ to get $\beta_u(aTb)=a\beta_u(T)b$. The last
statement follows from \eqref{cesaro0bis}.\endproof

We can then compose the vector states $\langle \xi_n,\bullet\,
\xi_n\rangle$ used in the construction of $\eta$ \eqref{conver0}
with $C_{\tau_n}$ to replace $[|D|,h]$ by $i|[D,h]|$.

Thus we need to determine the principal symbol of $|D|$. The
intuitive idea is as follows: one has
\begin{equation}\label{conjD}
    \beta_\tau(D)=e^{i \tau h}De^{-i \tau h}=D-i \tau[D,h]
\end{equation}
since $h$ commutes with $[D,h]$ so that $[D,f(h)]=f'(h)[D,h]$ for $f$ smooth
(\cf \cite{hilsum}). Thus, by homogeneity of the absolute value,
\begin{equation}\label{conjD1}
    \frac 1 \tau\beta_\tau(|D|)=\frac 1 \tau e^{i \tau h}|D|e^{-i \tau h}=|\frac
    D \tau-i[D,h]|\qqq \tau>0
\end{equation}
We need the weak limit in $\cH$ for $\tau\to \infty$ of $\frac 1
\tau\beta_\tau(|D|)\xi$ for $\xi\in \Dom D$. These vectors are
bounded in norm as follows from
\begin{equation} \label{conjD2}
\Vert \frac 1 \tau\beta_\tau(|D|)\xi\Vert=\Vert |\epsilon D
-i[D,h]|\xi\Vert=\Vert (\epsilon D -i[D,h])\xi\Vert\,, \ \ \epsilon
=1/\tau
\end{equation}
which is bounded since $\xi\in \Dom D$ so that $\Vert \epsilon D
\xi\Vert \to 0$. Note also that $\frac 1 \tau\beta_\tau(|D|)$ is a
positive operator so that any weak limit $\eta$ of $\frac 1
\tau\beta_\tau(|D|)\xi$ fulfills $\langle \xi,\eta\rangle \geq 0$.

Let us now show how to use regularity to obtain the strong
convergence of
\begin{equation} \label{conjD3}
\frac 1 \tau\beta_\tau(|D|)\xi=|\epsilon D -i[D,h]|\xi
\end{equation}
when $\epsilon \to 0$ and $\xi \in \Dom D$. By \eqref{conjD2} we can
assume that $\xi \in \cH_\infty$. We let $X(\epsilon)=\epsilon D
-i[D,h]$. By \eqref{conjD} it is  a self-adjoint operator with
$\cH_\infty$ as a core since $\cH_\infty$ is invariant under $e^{i
\tau h}$. The same holds for $|X(\epsilon)|$.

\begin{lem} \label{modofx} One has, with $X(\epsilon)=\epsilon D
-i[D,h]$,
\begin{equation}\label{modofx1}
   |X(\epsilon)|=Y(\epsilon)+f_0(X(\epsilon))
\end{equation}
where
\begin{equation}\label{modofx2}
    Y(\epsilon)\xi=\frac 2\pi\,
    \int_0^\infty\frac{X(\epsilon)^2}{1+u^2+X(\epsilon)^2}\xi\,du\qqq
    \xi \in \Dom D
\end{equation}
and
\begin{equation}\label{modofx3}
   f_0(x)=|x|-x^2(1+x^2)^{-1/2}\qqq
   x\in \R
\end{equation}
\end{lem}

\proof For any self-adjoint operator $T$ one has $\Vert
(1+u^2+T^2)^{-1}\Vert\leq (1+u^2)^{-1}$ and the norm convergent
expression
$$
(1+T^2)^{-1/2}=\frac 2\pi\, \int_0^\infty\frac{1}{1+u^2+T^2}\,du
$$
which gives for any $\xi\in \Dom T$,
$$
T^2(1+T^2)^{-1/2}\xi=\frac 2\pi\,
\int_0^\infty\frac{T^2}{1+u^2+T^2}\xi\,du
$$
Note that the partial sums
$$
\int_0^v\frac{T}{1+u^2+T^2}\,du
$$
are uniformly bounded but do not converge in norm to
$T(1+T^2)^{-1/2}$ since the function $x(1+x^2)^{-1/2}$ does not
vanish at $\infty$. Thus we get strong convergence on $\Dom T$. This
applies to $X(\epsilon)$ which is, up to a scale factor, conjugate
to $D$ by an automorphism of $\Dom D$ so that \eqref{modofx2} holds
with $Y(\epsilon)=f(X(\epsilon))$, $f(x)=x^2(1+x^2)^{-1/2}$. Finally
one has $f_0\in C_0(\R)$ and $f(x)+f_0(x)=|x|$.
\endproof

For each $\lambda \geq 0$ we define a transformation on operators
acting in $\cH_\infty$ by
\begin{equation}\label{thetalambda}
    \theta_\lambda(T)=(D^2+\lambda)T(D^2+\lambda)^{-1}
\end{equation}

\begin{lem} \label{dermodofx} Let $h=h^*\in \cA$, there exists
$\lambda<\infty$  such that
\begin{equation}\label{dermodofx1}
    \Vert \theta_\lambda((1+u^2+X(0)^2)^{-1})\Vert \leq
    (\frac 12+u^2)^{-1} \qqq u\,.
\end{equation}
\end{lem}

\proof Let $\xi\in \cH_\infty$ and let us give a lower bound for $
\Vert\theta_\lambda(1+u^2+X(0)^2)\xi\Vert$. Using
\begin{equation}\label{thetader}
    \theta_\lambda(T)=T+[D^2,T](D^2+\lambda)^{-1}
\end{equation}
we get
\begin{equation}\label{dermodofx2}
\theta_\lambda(X(0)^2)=X(0)^2-[D^2,[D,h]^2](D^2+\lambda)^{-1}
\end{equation}
Now the regularity shows (\cf \S \ref{regularity}) that
$[D^2,[D,h]^2](D^2+\lambda)^{-1}$ is compact so that for $\lambda\to
\infty$   its norm goes to $0$ (in fact it is of the form $B_\lambda
(D^2+\lambda)^{-1/2}$ with the norm of $B_\lambda$ bounded, so its
norm decays like $\lambda^{-1/2}$). Thus we can choose $\lambda$
large enough   so that $Z =\theta_\lambda(X(0)^2)-X(0)^2$ fulfills
$\Vert Z\Vert\leq 1/2$. We then get
$$
\langle \xi,\theta_\lambda(1+u^2+X(0)^2)\xi\rangle\geq \langle
\xi,1+u^2+X(0)^2\xi \rangle -|\langle \xi,Z\xi\rangle|\geq (\frac
12+ u^2)\Vert \xi\Vert^2
$$
(using $X(0)^2=-[D,h]^2\geq 0$) so that
\begin{equation}\label{dermodofx3}
\Vert \theta_\lambda(1+u^2+X(0)^2)\xi\Vert \geq (\frac 12+ u^2)\Vert
\xi\Vert \qqq \xi \in \cH_\infty\,.
\end{equation}
It remains to show that $\theta_\lambda(1+u^2+X(0)^2)$ is invertible
as an operator acting in $\cH_\infty$. Since $D^2+\lambda$ is an
automorphism of $\cH_\infty$, it is enough to show that
$1+u^2+X(0)^2$  is invertible as an operator acting in $\cH_\infty$.
One has $X(0)^2=-[D,h]^2\geq 0$ so that $1+u^2+X(0)^2$ is invertible
as an operator in $\cH$. Its invertibility in $\cH^\infty$ follows
from the stability under smooth functional calculus (Proposition
\ref{functcalc}) of the algebra
$$
\{T\in \cL(\cH)\,|\,T\cH_\infty \subset\cH_\infty\,,\
\Vert\delta^m(T)\Vert  <\infty \,,\ \forall m\}
$$
and the fact that, by regularity, $[D,h]$ belongs to this algebra.
\endproof

\begin{lem} \label{dermodofxbis} Let $h=h^*\in \cA$, then when
$\epsilon\to 0$,
\begin{equation}\label{dermodofx1ter}
    Y(\epsilon)\xi\to Y(0)\xi \qqq \xi \in \cH_\infty
\end{equation}
\end{lem}

\proof  One has for the action on $\cH_\infty$,
$$X(\epsilon)^2=(\epsilon D -i[D,h])^2= \epsilon^2D^2-i\epsilon(
D[D,h]+[D,h]D)-[D,h]^2$$
$$=\epsilon^2D^2-i\epsilon [D^2,h]-[D,h]^2$$ We first estimate, for $\xi \in \cH_\infty$,
$$
\eta(u,\epsilon)=(\frac{X(\epsilon)^2}{1+u^2+X(\epsilon)^2}-\frac{X(0)^2}{1+u^2+X(0)^2})\xi
$$
One has
$$
\eta(u,\epsilon)=(1+u^2)(\frac{1}{1+u^2+X(0)^2}-\frac{1}{1+u^2+X(\epsilon)^2})\xi=
$$
$$
\frac{1+u^2}{1+u^2+X(\epsilon)^2}(X(\epsilon)^2-X(0)^2)\frac{1}{1+u^2+X(0)^2}\xi
$$
$$
= \frac{1+u^2}{1+u^2+X(\epsilon)^2}(\epsilon^2D^2-i\epsilon [D^2,h]
)\frac{1}{1+u^2+X(0)^2}\xi
$$
$$
=\frac{1+u^2}{1+u^2+X(\epsilon)^2}(\epsilon^2D^2-i\epsilon [D^2,h]
)(D^2+\lambda)^{-1}\theta_\lambda((1+u^2+X(0)^2)^{-1})(D^2+\lambda)\xi
$$
Now one has, using regularity,
$$
\Vert (\epsilon^2D^2-i\epsilon [D^2,h]
)(D^2+\lambda)^{-1}\Vert=k(\epsilon)=O(\epsilon)
$$
while, since $X(\epsilon)$ is self-adjoint,
$$
\Vert \frac{1+u^2}{1+u^2+X(\epsilon)^2}\Vert \leq 1\,.
$$
Moreover $(D^2+\lambda)\xi\in \cH_\infty\subset \cH$. By Lemma \ref{dermodofx},
for $\lambda$ large enough, one thus gets
$$
\Vert\theta_\lambda((1+u^2+X(0)^2)^{-1})(D^2+\lambda)\xi\Vert\leq
(\frac 12+u^2)^{-1} \Vert (D^2+\lambda)\xi\Vert\qqq u\,.
$$
Thus after integrating in $u$ we get the following estimate
$$
\Vert(Y(\epsilon)- Y(0))\xi\Vert\leq \frac 2\pi\, \int_0^\infty
k(\epsilon)(\frac 12+u^2)^{-1} \Vert (D^2+\lambda)\xi\Vert
du=O(\epsilon)
$$
which gives the required result.
\endproof

It remains to estimate the continuity for $\epsilon\to 0$ of
$f_0(X(\epsilon))\xi$. The above proof shows that for
$g_a(x)=(a+x^2)^{-1}$ and any $a>0$ one has the norm continuity of
$g_a(X(\epsilon))\xi$ when $\epsilon \to 0$ (we showed convergence
only for $\xi\in \cH_\infty$ but it holds in general using the
boundedness of the functions $g_a$). The even functions in $f\in
C_0(\R)$ for which the following holds
\begin{equation}\label{continu}
   \Vert f(X(\epsilon))\xi-f(X(0))\xi\Vert \to 0 \qqq \xi \in \cH
\end{equation}
form a norm closed subalgebra of $C_0(\R)^{\rm even}$. This algebra
contains the functions $g_a$, thus the Stone-Weierstrass Theorem
shows that \eqref{continu} holds for all $f\in C_0(\R)^{\rm even}$
and in particular for $f_0$. We thus get:

\begin{prop}\label{symbabsd} Let $h=h^*\in \cA$, then
one has, with norm convergence:
\begin{equation}\label{symbabsd1}
\lim_{\tau\to \infty} \tau^{-1}\,e^{i\tau h}|D|e^{-i\tau
h}\xi=|[D,h]|\xi \qqq \xi \in \Dom D\,.
\end{equation}
\end{prop}

\proof By \eqref{conjD3} we just need to show that
$|X(\epsilon)|\xi\to |X(0)|\xi$ when $\epsilon \to 0$ for any $\xi
\in \cH_\infty$. By \eqref{modofx1}
$|X(\epsilon)|=Y(\epsilon)+f_0(X(\epsilon))$. By Lemma
\ref{dermodofx} we have $Y(\epsilon)\xi\to Y(0)\xi$ for $\xi \in
\cH_\infty$, and by the above discussion $f_0(X(\epsilon))\xi$ is
continuous at $\epsilon=0$. Thus we get the required result for $\xi
\in \cH_\infty$. The general case $\xi \in \Dom D$ follows using
\eqref{conjD2}.\endproof

\begin{rem} \label{modda}{\rm Proposition \ref{symbabsd} shows that,
under the regularity  hypothesis,
\begin{equation}\label{commsymbd}
    [|[D,h]|,[D,a]]=0\qqq h=h^*\,,\ a\in \cA\,.
\end{equation}
Indeed one has $$[e^{i\tau h}|D|e^{-i\tau h}, [D,a]] =e^{i\tau
h}[|D|,[D,a]] e^{-i\tau h}$$ and the norm of $[|D|,[D,a]]$ is finite
so that $\tau^{-1}\Vert [e^{i\tau h}|D|e^{-i\tau h}, [D,a]]\Vert\to
0$ for $\tau\to \infty$. Thus one has
$$
\lim_{\tau\to \infty} \tau^{-1}\,(e^{i\tau h}|D|e^{-i\tau
h}[D,a]\xi-[D,a]e^{i\tau h}|D|e^{-i\tau h}\xi)=0 \qqq \xi \in \Dom
D\,,
$$
and, since $[D,a]$ preserves $\Dom D$,
$$
[|[D,h]|,[D,a]]\xi=0 \qqq \xi \in \Dom D\,.
$$
Note also that, by the same argument, under the {\em strong
regularity} hypothesis of Definition \ref{stronregdefn} below, this
shows that
\begin{equation}\label{strongregcom}
   [D,h]^2\in \cA \qqq h=h^*  \in \cA\,.
\end{equation}
Indeed $|[D,h]|$ then commutes with all endomorphisms of
$\cH_\infty$. Its square $[D,h]^2$, being itself an endomorphism,
 belongs to the center of $\End_\cA(\cH_\infty)$ and is, by
\eqref{projpj3}, an element of $\cA$.

}\end{rem}

We can now show that regularity suffices to ensure the dissipative
property of Lemma \ref{diss}.

\begin{thm}\label{dissipthm} Let $(\cA,\cH,D)$ be a regular spectral
triple with $\cA$ commutative fulfilling the order one condition. Then for any
$h=h^*\in \cA$, the commutator $[D,h]$ vanishes where $h$ reaches its maximum,
\ie for any sequence $b_n\in \cA$, $\Vert b_n\Vert\leq 1$, with support tending
to $\{\chi\}$, where $\chi$ is a character such that $|\chi(h)|$ is maximum,
one has
$$
\Vert [D,h]b_n\Vert \to 0\,.
$$
\end{thm}

\proof By Proposition \ref{symbabsd} combined with the third
statement of Lemma \ref{symbroot1} one has, first for $\xi \in \Dom
D$ and then by uniformity for all $\xi \in \cH$,
$$
\lim_{\tau\to \infty} C_{\tau}([|D|,h])\xi=\lim_{\tau\to \infty}
\frac i\tau\,(e^{i\tau h}|D|e^{-i\tau h}\xi-|D|\xi)=i|[D,h]|\xi
$$
and thus the $C_{\tau}([|D|,h])$ converge strongly to $i|[D,h]|$
when $\tau\to \infty$.
 By the second statement of Lemma \ref{symbroot1}, one has
\begin{equation}\label{ctaunbis}
 C_{\tau}(b_n^*[|D|,h]b_n)=b_n^*C_{\tau}([|D|,h])b_n
 \end{equation}
Thus fixing $n$ and taking the limit for $\tau\to \infty$ one gets
\begin{equation}\label{ctaun3}
 C_{\tau}(b_n^*[|D|,h]b_n)\to i b_n^*|[D,h]|b_n
 \end{equation}
 One thus gets
 \begin{equation}\label{converend}
   \Vert b_n^*|[D,h]|b_n\Vert \leq \Vert b_n^* [|D|,h] b_n\Vert
\end{equation}
 But Lemma \ref{diss1} shows that $
   \Vert b_n^* [|D|,h] b_n\Vert\to 0
$ when $n\to\infty$ which gives the required result. Moreover, since
$|[D,h]|$ commutes with the $b_n$ this can be formulated as $ \Vert
[D,h]b_n\Vert \to 0\,. $
\endproof

\begin{cor} \label{dissipcor} Let $(\cA,\cH,D)$ be a spectral
triple with $\cA$ commutative fulfilling the five \axioms of \S
\ref{prelem}. The derivations $\pm\delta_j$ of Lemma
\ref{invariantintegration0} are dissipative.
\end{cor}

\proof This follows from Theorem \ref{dissipthm} and Lemma
\ref{diss}.\endproof

\begin{cor} \label{gradient} Let $h=h^*\in
\cA$. The principal symbol of the operator
\begin{equation}\label{nabla}
  {\rm Grad}  (h)=[D^2,h]
\end{equation}
vanishes where  $h$ reaches its maximum.
\end{cor}

\proof One has $[D^2,h]=D[D,h]+[D,h]D$ and since $[D,h]$ commutes
with $\cA$ one gets the principal symbol of $[D^2,h]$ from that of
$D$ which gives
\begin{equation}\label{nabla1}
\lim_{\tau\to\infty}\frac 1\tau
e^{i\tau\phi}[D^2,h]e^{i\tau\phi}=-i([D,\phi][D,h]+[D,h][D,\phi])
\end{equation}
Thus the result follows from Theorem \ref{dissipthm}.\endproof

\medskip \section{Self-adjointness and derivations}\label{sectselfadj}

We now introduce a technical hypothesis which will play an important
role.

\begin{defn} \label{stronregdefn} A spectral triple is {\em strongly regular} when all
endomorphisms of the $\cA$-module $\cH_\infty$ are regular.
\end{defn}

Our goal is to obtain self-adjoint operators from the operator $D$,
in the form $A^*DA$ where $A$ is regular \ie belongs to the domains
of $\delta^m$ for all $m$.

\begin{lem} \label{lemmaopda} Let $A$ be regular, then $A\Dom D\subset \Dom D$ and
the adjoint of $A^*D$ is the closure of the densely defined operator
$T$
\begin{equation}\label{opda}
    \Dom T=\Dom D\,,\ \ T\xi= D(A\xi)\qqq \xi \in \Dom D
\end{equation}
\end{lem}

\proof By regularity both $A$ and $A^*$ preserve the domain $\Dom
|D|=\Dom D$ so that \eqref{opda} makes sense. The domain of $A^*D$
is the domain of $D$. An $\eta \in \cH$ belongs to the domain of the
adjoint $S=(A^*D)^*$ when there exists a constant $C<\infty$ such
that
\begin{equation}\label{domadj}
    |\langle A^*D\xi,\eta\rangle|\leq C\Vert \xi\Vert \qqq \xi \in
    \Dom D
\end{equation}
One has $\langle A^*D\xi,\eta\rangle=\langle  D\xi,A\eta\rangle$
and, since $D$ is self-adjoint, the above condition means that
$A\eta\in \Dom D$. Moreover one then has  $S\eta=D A\eta$. In other
words $S=DA$ with domain
\begin{equation}\label{opT}
    \Dom S=\{\eta\,|\, A\eta\in\Dom D\}\,,\ \ S\xi= D(A\xi)\qqq \xi \in \Dom
    S
\end{equation}
To prove the Lemma we need to show that $S$ is the closure of the
operator $T$ of \eqref{opda}. Let $\eta \in \Dom S$. We construct a
sequence $\eta_n\in \Dom D$ such that
\begin{equation}\label{sequetan}
    \eta_n \to \eta \,, \  \ DA\eta_n\to DA\eta
\end{equation}
In fact we let
\begin{equation}\label{sequetan1}
    \eta(\epsilon)=(1+\epsilon |D|)^{-1} \eta \qqq \epsilon >0\,.
\end{equation}
It belongs to $\Dom D$ by construction and $\eta(\epsilon) \to \eta$
when $\epsilon \to 0$. One has
$$
DA\eta(\epsilon)=D (1+\epsilon |D|)^{-1} A\eta +\epsilon D
(1+\epsilon |D|)^{-1}[|D|,A](1+\epsilon |D|)^{-1}\eta
$$
Since $ A\eta\in\Dom D$ one has $D (1+\epsilon |D|)^{-1} A\eta=
(1+\epsilon |D|)^{-1}D A\eta\to DA\eta$. The remainder is of the
form $B(\epsilon)[|D|,A]\eta(\epsilon)$ where $B(\epsilon)=\epsilon
D (1+\epsilon |D|)^{-1}$ is of norm less than $1$, $[|D|,A]$ is
bounded and $\eta(\epsilon) \to \eta$. Thus it behaves like
$B(\epsilon)[|D|,A]\eta$ and hence tends to $0$ when $\epsilon \to
0$, since $B(\epsilon)\zeta\to 0$ for any $\zeta\in \cH$. This shows
that $DA\eta(\epsilon)\to DA\eta$ and $S$ is the closure of $T$.
\endproof

\begin{cor} \label{coropda} Let $\varphi=\varphi^*\in \cA$ then the operator
$H=\varphi D\varphi$ with domain $\Dom D$ is essentially
self-adjoint.
\end{cor}

\proof One has $H=\varphi^2 D+\varphi[ D,\varphi]$ on $\Dom D$. The
bounded perturbation $P=\varphi[ D,\varphi]$ does not alter the
domain of the adjoint $H^*$ which is thus the same as the domain of
$H_0^*$, $H_0=\varphi^2 D$. By Lemma \ref{lemmaopda}, the adjoint of
$H_0$ is the closure of $D \varphi^2$ with domain $\Dom D$. This is
the same as the closure of $\varphi D\varphi+[ D,\varphi]\varphi$
with domain $\Dom D$. Since $[ D,\varphi]\varphi$ is bounded we thus
get that the adjoint $H_0^*$ of $H_0$ is the sum of the closure of
$\varphi D\varphi$ with domain $\Dom D$ with the bounded operator $[
D,\varphi]\varphi$. Thus when adding $P^*=-[ D,\varphi]\varphi$ to
$H_0^*$ we obtain the closure of $\varphi D\varphi$ with domain
$\Dom D$, \ie the operator $H$. \endproof

\begin{lem}\label{lemmaopda2} Let $A$ be regular. Then $A^*D$ is
closable and
\begin{itemize}
  \item For any $\xi$ in the domain of the closure $\overline{A^*D}$ of
$A^*D$, one has, for $\epsilon >0$,
\begin{equation}\label{domclos}
    (1+\epsilon |D|)^{-1}\overline{A^*D}\xi=A^*D(1+\epsilon
    |D|)^{-1}\xi- (1+\epsilon |D|)^{-1}[|D|,A^*]\epsilon D(1+\epsilon
    |D|)^{-1}\xi
\end{equation}
  \item The domain of $\overline{A^*D}$ is the set of $\xi\in \cH$
for which the $A^*D(1+\epsilon |D|)^{-1}\xi$ converge in norm for
$\epsilon \to 0$.
  \item The limit of the $A^*D(1+\epsilon |D|)^{-1}\xi$ gives
  $\overline{A^*D}\xi$.
\end{itemize}
\end{lem}

\proof The operator $A^*D$ is closable since its adjoint is densely
defined by Lemma \ref{lemmaopda}. The right hand side of
\eqref{domclos} is a bounded operator, thus it is enough to prove
the equality for $\xi \in \Dom D$ since $\overline{A^*D}$ is the
closure of its restriction to $\Dom D$. For $\xi \in \Dom D$
\eqref{domclos} follows from
$$
[(1+\epsilon |D|)^{-1},A^*]= -(1+\epsilon |D|)^{-1}[\epsilon|D|,A^*]
(1+\epsilon |D|)^{-1}\,.
$$
Let then $\xi$ be in the domain of the closure $\overline{A^*D}$. By
\eqref{domclos}, $A^*D(1+\epsilon |D|)^{-1}\xi$ is the sum of
$(1+\epsilon |D|)^{-1}\overline{A^*D}\xi\to\overline{A^*D}\xi$, and
of $(1+\epsilon |D|)^{-1}[|D|,A^*]\epsilon D(1+\epsilon
    |D|)^{-1}\xi$ which converges to $0$ in norm since $(1+\epsilon
    |D|)^{-1}[|D|,A^*]$
    is uniformly bounded while $\epsilon
D(1+\epsilon
    |D|)^{-1}\xi$ converges to $0$ in norm. Thus $A^*D(1+\epsilon |D|)^{-1}\xi$
    is convergent when $\epsilon\to 0$. Conversely, if the $A^*D(1+\epsilon |D|)^{-1}\xi$ converge in norm for
$\epsilon \to 0$, then since $(1+\epsilon |D|)^{-1}\xi\to \xi$ and
$(1+\epsilon |D|)^{-1}\xi \in \Dom D$, one gets that $\xi$ is in the
domain of the closure $\overline{A^*D}$ of $A^*D$ and that moreover
the limit of the $A^*D(1+\epsilon |D|)^{-1}\xi$ gives
  $\overline{A^*D}\xi$.
\endproof

\begin{prop}\label{selfsqueeze} Let $A$ be regular then the operator $H=A^*DA$ with
domain $\Dom D$ is essentially self-adjoint. The domain of the
closure of $H$ is the set  of $\xi\in \cH$ for which the
$A^*DA(1+\epsilon |D|)^{-1}\xi$ converge in norm for $\epsilon \to
0$. The limit of the $A^*DA(1+\epsilon |D|)^{-1}\xi$ gives
  $\overline{H}\xi$.
\end{prop}

\proof Let us first check that $H$ is symmetric. One has for $\xi$
and $\eta$ in $\Dom D$,
$$
\langle H\xi,\eta\rangle=  \langle A^*DA\xi,\eta\rangle=\langle
DA\xi,A\eta\rangle=\langle A\xi,DA\eta\rangle=\langle
\xi,A^*DA\eta\rangle=\langle
 \xi,H\eta\rangle
$$
Let us now show that $H^*$ is the closure of $H$. Let $\eta \in \Dom
H^*$. Then there exists $C<\infty$ with
$$
 |\langle A^*DA\xi,\eta\rangle| \leq C \Vert\xi \Vert\qqq \xi\in \Dom
 D\,.
$$
Since $\langle A^*DA\xi,\eta\rangle=\langle  DA\xi,A\eta\rangle$,
this means that $A\eta$ is in the domain of the adjoint of $DA$ with
domain $\Dom D$ \ie $$ A\eta \in \Dom T^*\,,\ \ H^*\eta=T^*A\eta
$$ where $T$ is defined in \eqref{opda}. By Lemma
\ref{lemmaopda} the adjoint  of $A^*D$ is the closure of $T$:
$(A^*D)^*=\bar T$. The adjoint  $T^*$ of $T$ is the same as the
adjoint of  the closure $\bar T$, and is the closure
$\overline{A^*D}=(A^*D)^{**}$ of $A^*D$. Thus by Lemma
\ref{lemmaopda2} we have, since $A\eta$ is in the domain of
$\overline{A^*D}$, the convergence of $A^*D(1+\epsilon
|D|)^{-1}A\eta$ to $\overline{A^*D}A\eta=H^*\eta$. Moreover, as
above, we have
$$
A^*D(1+\epsilon |D|)^{-1}A\eta-A^*DA(1+\epsilon
|D|)^{-1}\eta=-A^*\epsilon D(1+\epsilon |D|)^{-1}[|D|,A] (1+\epsilon
|D|)^{-1}\eta
$$
and the right hand side converges to $0$ in norm when $\epsilon\to
0$. Thus we have shown that for any $\eta \in \Dom H^*$ one gets the
convergence of $A^*DA(1+\epsilon |D|)^{-1}\eta$ to $H^*\eta$. This
shows, since $(1+\epsilon |D|)^{-1}\eta\in \Dom H$, that $H^*$ is
the closure of $H$ and hence that $H$ is essentially self-adjoint.
It also gives a characterization of the domain of the closure of $H$
as required.\endproof

\smallskip

We now want to apply this result using endomorphisms of the
$\cA$-module $\cH_\infty$ which are of rank one, in order to obtain
an operator on $\cA$ itself.

\begin{lem}\label{rankonecstar} Let $\xi,\eta\in \cH_\infty$
then the following gives an endomorphism of the $\cA$-module
$\cH_\infty$:
\begin{equation}\label{rankonecstar1}
    T_{\xi,\eta}(\zeta)=(\eta|\zeta )\xi \qqq
    \zeta \in \cH_\infty
\end{equation}
where $(\eta|\zeta )$ is the $\cA$-valued inner product. One has
\begin{equation}\label{rankonecstar2}
    T_{a\xi,\,b\eta}=ab^*T_{\xi,\eta}\qqq a, b\in \cA\,,\ \
    T_{\xi,\eta}^*=T_{\eta,\xi}
\end{equation}
\end{lem}

\proof This follows from the $\cA$-linearity of the inner product,
which is linear in the second variable and antilinear in the first.
The equality $T_{\xi,\eta}^*=T_{\eta,\xi}$ follows from
$$
\langle T_{\eta,\xi}\alpha,\beta\rangle=\langle
(\xi|\alpha)\eta,\beta\rangle=\int
(\xi|\alpha)^*(\eta|\beta)d\lambda
$$
$$
\langle \alpha,T_{\xi,\eta}\beta\rangle=\langle
\alpha,(\eta|\beta)\xi\rangle=\int (\alpha|\xi)(\eta|\beta)d\lambda
$$
\endproof

By Proposition \ref{generaltop} (4), the $T_{\xi,\eta}$ are bounded operators
in $\cH$. Let us now assume that all endomorphisms of the $\cA$-module
$\cH_\infty$ are regular as in Definition \ref{stronregdefn}. We can then apply
Proposition \ref{selfsqueeze} and get that the following defines an essentially
self-adjoint operator with domain $\Dom D$,
\begin{equation}\label{opdxieta}
    D_{\xi,\,\eta}=T_{\eta,\,\xi}\,D\,T_{\xi,\,\eta}
\end{equation}
We need to relate this operator with the derivation of $\cA$ given
by \eqref{formforder} \ie
\begin{equation}\label{formforderbis}
    \delta_0(a)=i(\xi|[D,a]\xi)\qqq a \in \cA
\end{equation}

\begin{lem}\label{rankonecstar3} One has
\begin{equation}\label{rankonecstar4}
    D_{\xi,\,\eta}\zeta=-i \,\delta_0((\eta|\zeta))\,\eta +
    (\xi|D\xi)T_{\eta,\,\eta}\zeta \qqq \xi, \eta,\zeta\in
    \cH_\infty
\end{equation}
The operator $V_\eta(a)=a\eta $, $\forall a \in \cA$ extends to a
bounded linear map $V_\eta$ from $L^2(X,d\lambda)$ to $\cH$, and one
has
\begin{equation}\label{veta}
    V_\eta^*V_\eta=(\eta|\eta)\, ,\ V_\eta V_\eta^*=T_{\eta,\,\eta}  \, ,\  V_\eta^*(\zeta)=(\eta|\zeta)
    \qqq \zeta \in \cH_\infty
\end{equation}
\end{lem}

\proof One has
$D_{\xi,\,\eta}\zeta=T_{\eta,\,\xi}\,D\,T_{\xi,\,\eta}\zeta=T_{\eta,\,\xi}\,D((\eta|\zeta
)\xi) $. Thus using $(\xi,[D, a]\xi)=-i\delta_0(a)$ and
$(\xi|aD\xi)\eta=(\xi|D\xi)a\eta$ for $a=(\eta|\zeta )$ one gets
$$
D_{\xi,\,\eta}\zeta=(\xi|Da\xi)\eta=(\xi|[D,a]\xi)\eta +
(\xi|D\xi)a\eta=-i \delta_0((\eta|\zeta))\eta +
    (\xi|D\xi)T_{\eta,\,\eta}\zeta
$$
which gives \eqref{rankonecstar4}. To show that $V_\eta$ is bounded,
note that
$$
\langle V_\eta(a), V_\eta(a)\rangle=\langle a\eta , a\eta
\rangle=\cutint a^*a (\eta| \eta)|D|^{-p}=\int a^*a (\eta|
\eta)d\lambda
$$
which also shows that $V_\eta^*V_\eta=(\eta|\eta)$. Let us check
that $V_\eta^*(\zeta)=(\eta|\zeta)$. One has
$$
\langle \zeta,V_\eta(a)\rangle=\langle \zeta,a\eta\rangle=\int a
(\zeta| \eta)d\lambda=\int (\eta| \zeta)^* a d\lambda=\langle
V_\eta^*(\zeta),a\rangle
$$
The equality $V_\eta V_\eta^*=T_{\eta,\,\eta}$ follows from
\eqref{rankonecstar1}.\endproof

The strategy now is to use the self-adjointness of $D_{\xi,\,\eta}$
and the fact that $\delta_0$ can be compared to $iD_{\xi,\,\eta}$
plus a bounded perturbation to show that the resolvent problem
$(1+\epsilon\delta_0)\xi=\eta$ can be solved first in $L^2$. Then
one wants to use the regularity to show that this problem can also
be solved in the Sobolev spaces. Finally one wants to use the
Sobolev estimates to show that it can be solved in the $C^*$ norm.
Then together with the results on dissipative derivations of \S
\ref{dissipsect} one gets the existence of the resolvent for the
action on the $C^*$-algebra. One notes that it is enough to solve
the resolvent problem for $\epsilon$ small enough. One then applies
the Hille-Yoshida Theorem.

\smallskip

More specifically we consider the equation
\begin{equation}\label{resolvequ1}
   (1+i\epsilon(D_{\xi,\,\eta}-(\xi|D\xi)T_{\eta,\,\eta}))\zeta=a\eta
\end{equation}
where $a\in \cA$ is given and $\epsilon$ can be taken as small as
needed. Given a solution $\zeta$ of \eqref{resolvequ1}, one can
under suitable regularity conditions on $\zeta$ take the inner
product $(\eta|\zeta)=b$. One then has, at the formal level
\begin{equation}\label{resolvequ2}
   b+\epsilon \delta_0(b)(\eta|\eta)=a (\eta|\eta)
\end{equation}
We can assume that the support of $\xi$ is small enough so that we
can find $\eta$ such that $(\eta|\eta)=1$ in a neighborhood of the
support $K$ of $\xi$. Then by \eqref{formforderbis} one gets that
since $\delta_0$ vanishes outside $K$, one can replace
$\delta_0(b)(\eta|\eta)$ in \eqref{resolvequ2} by $\delta_0(b)$.
Moreover one then gets
\begin{equation}\label{resolvequ3}
   c+\epsilon \delta_0(c)=a \,,\  \ c=b+(1-(\eta|\eta))a
\end{equation}
since $(1-(\eta|\eta))a$ belongs to the kernel of $\delta_0$ since
its support is disjoint from $K$. We need to know that $c\in A$
where $A=C(X)$ is the norm closure of $\cA$ and in fact also that
$[D,c]$ is bounded, just to formulate the result. Thus we need to
control the Sobolev norms of the solution of \eqref{resolvequ1}. To
do that we use the transformation $\theta_\lambda$ of
\eqref{thetalambda}. One has, as in \eqref{thetader},
\begin{equation}\label{thetaderbis}
    \theta_\lambda(T)=T+\cE_\lambda(T)\,, \ \ \cE_\lambda(T)=[D^2,T](D^2+\lambda)^{-1}
\end{equation}
so that the binomial formula expresses $\theta^N_\lambda(T)$ in
terms of the $\cE^k_\lambda(T)$ for $k\leq N$, as
\begin{equation}\label{thetaderbin}
    \theta_\lambda^N(T)=T+\sum_{k\geq 1}
\left(
       \begin{array}{c}
         N \\
         k \\
       \end{array}
     \right)
\cE^k_\lambda(T)
\end{equation}
Note also that, for $T$ regular, and on $\Dom D$ one has
\begin{equation}\label{thetanddelta}
[D^2,T]=2\delta(T)|D|+\delta^2(T)
\end{equation}
as follows from
$[D^2,T]=[|D|^2,T]=\delta(T)|D|+|D|\delta(T)=2\delta(T)|D|+\delta^2(T)$.

\begin{lem} \label{resolvequlem1prel}  Let $T$ be regular.

1) The $\cE^k_\lambda(T)$ are compact operators for $k>0$ and
converge in norm to $0$ when $\lambda \to \infty$.

2) One has (with $\lambda \geq 1$)
\begin{equation}\label{leftrightbound}
   \Vert\cE_\lambda(T)D\Vert\leq 2\Vert \delta(T)\Vert+\Vert \delta^2(T)\Vert\,, \ \
    \Vert D\cE_\lambda(T)\Vert\leq 2\Vert \delta(T)\Vert+3\Vert \delta^2(T)\Vert+\Vert \delta^3(T)\Vert
\end{equation}

3) For $k>1$, the operators $D\cE^k_\lambda(T)$ and
$\cE^k_\lambda(T)D$ are compact operators which converge in norm to
$0$ when $\lambda \to \infty$.
\end{lem}

\proof 1) One has, using \eqref{thetaderbis} and
\eqref{thetanddelta}, that
\begin{equation}\label{leftrightbound0}
\cE_\lambda(T)=(2\delta(T)|D|+\delta^2(T))(D^2+\lambda)^{-1}
\end{equation}
 Thus
the answer follows for $k=1$ since both $(D^2+\lambda)^{-1}$ and
$|D|(D^2+\lambda)^{-1}$ are compact operators which converge in norm
to $0$ when $\lambda \to \infty$. Since $\cE_\lambda(T)$ is also
regular it follows also for $k>1$.

2) The first inequality of \eqref{leftrightbound} follows from
\eqref{leftrightbound0}.
 For the second, one has
$$
[|D|,\cE_\lambda(T)]=\cE_\lambda(\delta(T))
$$
which gives \eqref{leftrightbound} using,
$$
\Vert D\cE_\lambda(T)\Vert \leq \Vert \cE_\lambda(T)|D|\Vert + \Vert
\cE_\lambda(\delta(T))\Vert
$$
and the second inequality follows using the first and
\eqref{leftrightbound0}.

3) The statement is immediate for $\cE^k_\lambda(T)D$ since
$|D|(D^2+\lambda)^{-1}$ is compact. For the second one uses
$[|D|,\cE^k_\lambda(T)]=\cE_\lambda^k(\delta(T))$ as in the proof of
2).
\endproof

\begin{lem} \label{resolvequlem1}
1) For any integer $N\in \N$, there exists $\lambda <\infty$ and
$\epsilon_0>0$ such that the following operator with domain $\Dom D$
is closable and invertible in $\cH$ for any $\epsilon\leq
\epsilon_0$,
\begin{equation}\label{opinH}
   \theta_\lambda^N(1+i\epsilon S_{\xi,\,\eta})  \,, \
   S_{\xi,\,\eta}=D_{\xi,\,\eta}-(\xi|D\xi)T_{\eta,\,\eta}
\end{equation}
and the norm of its inverse fulfills
\begin{equation}\label{opinHbis}
  \Vert
  (\theta_\lambda^N(1+i\epsilon S_{\xi,\,\eta}))^{-1}\Vert\leq
  1+N c_{\xi,\,\eta}\epsilon\,,
\end{equation}
where $c_{\xi,\,\eta}<\infty$ only depends on $\xi$ and $\eta$.

2)  For any integer $N$ there exists $\epsilon_N>0$ such that
\eqref{resolvequ1} can be solved in $\cH_N=\Dom |D|^N$.
\end{lem}

\proof  The operator $P=(\xi|D\xi)T_{\eta,\,\eta}$ is bounded and regular since
it is an endomorphism of the $\cA$-module $\cH_\infty$. Thus it preserves the
domain of $(D^2+\lambda)^{N}$ and the $\cE^k_\lambda(P)$ are compact operators
for $k>0$ and converge in norm to $0$ when $\lambda \to \infty$ by Lemma
\ref{resolvequlem1prel}. By \eqref{opdxieta}, one has
$D_{\xi,\,\eta}=T_{\eta,\,\xi}\,D\,T_{\xi,\,\eta}$ thus, by regularity of the
$T_{\xi,\,\eta}$, the operator
$$
\theta^N_\lambda(D_{\xi,\,\eta})=(D^2+\lambda)^{N}D_{\xi,\,\eta}(D^2+\lambda)^{-N}
$$
is well defined on $\Dom D$. Moreover one has, by
\eqref{thetaderbin}, and on $\Dom D$,
$$
\theta^N_\lambda(D_{\xi,\,\eta})=\theta^N_\lambda(T_{\eta,\,\xi})
\theta^N_\lambda(D)\theta^N_\lambda(T_{\xi,\,\eta})= \sum_{k,m}
\left(
       \begin{array}{c}
         N \\
         k \\
       \end{array}
     \right)
\cE^k_\lambda(T_{\eta,\,\xi})D \left(
       \begin{array}{c}
         N \\
         m \\
       \end{array}
     \right)
     \cE^m_\lambda(T_{\xi,\,\eta})
$$
so that one gets
\begin{equation}\label{epsilonlambdad}
    \theta^N_\lambda(D_{\xi,\,\eta})=D_{\xi,\,\eta}+N
    \cE_\lambda(T_{\eta,\,\xi})D T_{\xi,\,\eta}+N
    T_{\eta,\,\xi}D\cE_\lambda(T_{\xi,\,\eta})+Q(N,\lambda)
\end{equation}
where the remainder $Q(N,\lambda)$ is a sum of terms proportional to
$\cE^k_\lambda(T_{\eta,\,\xi})D \cE^m_\lambda(T_{\xi,\,\eta})$ for
$k+m>1$. By Lemma \ref{resolvequlem1prel} we get that $Q(N,\lambda)$
is a compact operator and $\Vert Q(N,\lambda)\Vert \to 0$ when
$\lambda \to \infty$. Thus for $\lambda\to \infty$,  we get the
following estimate: there exists $C_{\xi,\,\eta}<\infty$ only
depending on $\xi$ and $\eta$ such that
\begin{equation}\label{epsilonlambdad1}
   \liminf_{\lambda\to \infty} \Vert \theta^N_\lambda(D_{\xi,\,\eta})-D_{\xi,\,\eta}\Vert \leq N
    C_{\xi,\,\eta}\,.
\end{equation}
Since the $\cE^k_\lambda(P)$ are compact operators for $k>0$ and
converge in norm to $0$ when $\lambda \to \infty$, one gets
similarly
\begin{equation}\label{epsilonlambdad2}
   \liminf_{\lambda\to \infty} \Vert \theta^N_\lambda(S_{\xi,\,\eta})-S_{\xi,\,\eta}\Vert \leq N
    C_{\xi,\,\eta}\,.
\end{equation}
Let $\lambda$ be large enough so that $\Vert
\theta^N_\lambda(S_{\xi,\,\eta})-S_{\xi,\,\eta}\Vert \leq 2 N
    C_{\xi,\,\eta}$. For $\epsilon$ small enough,
    \begin{equation}\label{epsilonlambdad3}
\Vert \theta^N_\lambda(1+i\epsilon S_{\xi,\,\eta})-(1+i\epsilon
D_{\xi,\,\eta})\Vert \leq 2 N
    C_{\xi,\,\eta}\epsilon + \epsilon \Vert
    (\xi|D\xi)T_{\eta,\,\eta}\Vert < 1
    \end{equation}
     Since $D_{\xi,\,\eta}$ is essentially self-adjoint, the
operator $K=1+i\epsilon D_{\xi,\,\eta}$ is closable and invertible
for any $\epsilon >0$ and the norm of its inverse is  $\leq 1$.
$\theta^N_\lambda(1+i\epsilon S_{\xi,\,\eta})$ is closable since it
is a bounded perturbation $K-B$ of $K=1+i\epsilon D_{\xi,\,\eta}$.
Moreover by \eqref{epsilonlambdad3} it is invertible, and the norm
of its inverse, which is given by the Neumann series
$(K-B)^{-1}=\sum (K^{-1}B)^m K^{-1}$, fulfills \eqref{opinHbis}.
    This proves the first statement.

2) Let
    $\xi=a\eta\in \cH_{2N}$, and consider
    $\xi'=(D^2+\lambda)^{N}\xi\in \cH$. Then, by the first statement, one
    can find a sequence $\zeta_n\in \Dom D$, $\zeta_n\to \zeta' \in \cH$ such that
    $$
\theta_\lambda^N(1+i\epsilon(D_{\xi,\,\eta}-(\xi|D\xi)T_{\eta,\,\eta}))\zeta_n\to\xi'
    $$
    where the convergence is in $\cH$.
Applying the bounded operator $(D^2+\lambda)^{-N}$  gives
$$
(1+i\epsilon(D_{\xi,\,\eta}-(\xi|D\xi)T_{\eta,\,\eta}))(D^2+\lambda)^{-N}\zeta_n\to
\xi\,.
$$
One has $\zeta =(D^2+\lambda)^{-N}\zeta'\in \cH_{2N}\subset \Dom D$
and $(D^2+\lambda)^{-N}\zeta_n\to (D^2+\lambda)^{-N}\zeta'$ in the
topology of $\Dom D$. Thus $$
(1+i\epsilon(D_{\xi,\,\eta}-(\xi|D\xi)T_{\eta,\,\eta}))\zeta=\xi $$
and $\zeta\in \cH_{2N}$ gives the required solution.\endproof

We now need to show that if $\eta\in \cH_\infty$ and $\zeta \in
\cH_{N}$ for $N$ large enough, the inner product $(\eta|\zeta)$
gives an element of $A=C(X)$ and in fact in the domain of
$\delta^k$. To see this we use Proposition \ref{generaltop}. We
recall that the Sobolev norms on $\cA$ are defined using generators
$\eta_\mu$ of the $\cA$-module $\cH_\infty$ by \eqref{sobolevnorms}
\ie
$$
 \Vert a\Vert_s^{\rm sobolev}=(\sum_\mu
    \Vert(1+D^2)^{s/2}a\eta_\mu \Vert^2)^{1/2} \qqq a\in \cA
$$
Thus when we want to control the Sobolev norms of $(\eta|\zeta)$ we
need to control the norms
$$
\Vert(1+D^2)^{s/2}(\eta|\zeta)\eta_\mu \Vert
$$
The point then is that
$(\eta,\zeta)\eta_\mu=T_{\eta_\mu,\,\eta}\zeta$ while the
endomorphism $T_{\eta_\mu,\,\eta}$ is regular by hypothesis so that
$\theta_\lambda^N(T_{\eta_\mu,\,\eta})$ is bounded and (with
$\lambda=1$) one gets:

\begin{lem} \label{strongreglem}
Assuming strong regularity, one has for $\eta\in \cH_\infty$
\begin{equation}\label{sobolevinnerbound}
    \Vert (\eta|\zeta)\Vert_s^{\rm sobolev}\leq C_s
    \Vert(1+D^2)^{s/2}\zeta\Vert
\end{equation}
\end{lem}

\proof It is enough to prove the estimate when $s/2=N$ is an
integer. For each $\mu$ one has
$$
\Vert(1+D^2)^{N}(\eta|\zeta)\eta_\mu \Vert=\Vert
\theta_1^N(T_{\eta_\mu,\,\eta})(1+D^2)^{N}\zeta\Vert\leq \Vert
\theta_1^N(T_{\eta_\mu,\,\eta}) \Vert \Vert (1+D^2)^{N}\zeta\Vert
$$
\endproof

\begin{thm} \label{hileyosidathm} Let $(\cA,\cH,D)$ be a strongly regular spectral triple with $\cA$
commutative, fulfilling the five \axioms of \S \ref{prelem}. Then any
derivation of $\cA$ of the form \eqref{formforder} \ie
$\delta_0(a)=i(\xi|[D,a]\xi)$, $ \forall a \in \cA$, is closable for
the $C^*$-norm of $A$, and its closure is the generator of a
one-parameter group of automorphisms $U(t)$ of $A=C(X)$,
$X=\Sp(\cA)$.
\end{thm}

\proof By Corollary \ref{dissipcor} the derivation $\delta_0$, with
domain $\cA\subset A$, is dissipative for the $C^*$-norm of $A$.
Thus it is closable (\cite{Bratteli} Proposition 1.4.7) and we let
$D(\delta_0)$ be the domain of its closure. To apply the
Hille-Yosida-Lumer-Phillips Theorem we need to show that for
sufficiently small $\epsilon$ one has
\begin{equation}\label{surjectc0}
    (1+\epsilon \delta_0)D(\delta_0)=A
\end{equation}
By Corollary \ref{dissipcor}, we have
$$
\Vert (1+\epsilon \delta_0)(x)\Vert \geq \Vert x\Vert \qqq x\in
D(\delta_0)
$$
Thus $(1+\epsilon \delta_0)D(\delta_0)$ is closed in norm and it is
enough to show that $(1+\epsilon \delta_0)\cA$ is norm dense in $A$.
Let then $\eta \in \cH_\infty$ be such that $(\eta|\eta)=1$ in a
neighborhood of the support of $\xi$ (with
$\delta_0(a)=i(\xi|[D,a]\xi)$). Let then $N\in \N$ be such that the
Sobolev estimate holds (Proposition \ref{generaltop})
\begin{equation}\label{sobestiexpected}
    \Vert a\Vert_{C^*}   \leq C \Vert a\Vert_N^{\rm
    sobolev} \qqq a \in \cA
\end{equation}
Let $a\in \cA$ one has $a\eta\in \cH_\infty$. By Lemma
\ref{resolvequlem1} there exists $\epsilon_{N+1}
>0$ such that for any $\epsilon \leq \epsilon_{N+1}$ one can find a
solution in $\zeta\in\cH_{N+1}$ of the equation \eqref{resolvequ1}.
Since $\cH_\infty $ is dense in $\cH_{N+1}$ we thus get a sequence
$\zeta_n \in \cH_\infty$ such that $\zeta_n \to \zeta$ in
$\cH_{N+1}$. The operator
$S_{\xi,\,\eta}=D_{\xi,\,\eta}-(\xi|D\xi)T_{\eta,\,\eta})$ is
continuous from $\cH_{N+1}$ to $\cH_N$. One thus has, with
convergence in $\cH_N$:
$$
(1+i\epsilon S_{\xi,\,\eta})\zeta_n\to (1+i\epsilon
S_{\xi,\,\eta})\zeta=a\eta
$$
Combining Lemma \ref{strongreglem} with \eqref{sobestiexpected}, one gets that
the $b_n=(\eta|\zeta_n)\in \cA$ converge in the $C^*$-norm $\Vert x\Vert$.
Moreover, by \eqref{rankonecstar4} and \eqref{opinH},
$$
 S_{\xi,\,\eta}\zeta_n=-i \,\delta_0((\eta|\zeta_n))\,\eta\,,\
 (1+i\epsilon S_{\xi,\,\eta})\zeta_n=\zeta_n+\epsilon
 \delta_0(b_n)\,\eta\to a\eta
$$
with convergence in $\cH_N$. Thus applying $(\eta|\bullet)$ and
using \eqref{sobolevinnerbound} and \eqref{sobestiexpected},
$$
b_n+\epsilon \delta_0(b_n)(\eta|\eta)\to a (\eta|\eta)
$$
in the $C^*$-norm, as in \eqref{resolvequ2}. Since $(\eta|\eta)=1$
in a neighborhood of the support of $\xi$ one has
$\delta_0(b_n)(\eta|\eta)=\delta_0(b_n)$. Moreover one has, since
$(1-(\eta|\eta))a$ vanishes in a neighborhood of the support of
$\xi$ that $\delta_0((1-(\eta|\eta))a)=0\,.$ Thus we have the norm
convergence
$$
c_n+\epsilon \delta_0(c_n)\to a \,, \ \ c_n=b_n+ (1-(\eta|\eta))a
$$
and this shows that $(1+\epsilon \delta_0)\cA$ is norm dense in $A$.
Since $(1+\epsilon \delta_0)D(\delta_0)$ is norm closed it is equal
to $A$. Thus, we have shown that for sufficiently small $\epsilon$
one has \eqref{surjectc0}. Thus the Hille-Yosida-Lumer-Phillips
Theorem (\cite{Bratteli} Theorem 1.5.2, \cite{Simon} Theorem X.47 a)
shows that $\delta_0$ generates a contraction semi-group of $A$.
Since the same holds for $-\delta_0$, one gets a one parameter group
of isometries $U(t)=e^{t\delta_0}$ of the $C^*$-algebra $A$.
Moreover $U(t)(a)$ is a norm continuous function of $t$ for fixed
$a\in A$. Using the operators of the form
\begin{equation}\label{smearing}
   U(f)=\int f(t)U(t)dt\,:\, A\to A
\end{equation}
for $f$ such that the $L^1$-norms of the derivatives $\Vert
f^{(n)}\Vert_1$ fulfill $\sum \frac{t^n}{n!}\Vert f^{(n)}\Vert_1
<\infty$, one gets a dense domain of analytic elements and one
checks that since $\delta_0$ is  a derivation on $D(\delta_0)$ the
$U(t)$ are automorphisms of $A$.\endproof

It remains to show that the $U(t)\in \Aut(A)$ respect the
smoothness. Let us first show that we need only understand what
happens to $U(t)(a)\eta$ as an element of $\cH$ because $U(t)$ is
the identity in the complement of the support of $\xi$.

\begin{lem} \label{supportlem} Let $x\in A$ have support disjoint from the support of $\xi$.
Then $U(t)(x)=x$ for all $t\in \R$.
\end{lem}

\proof We can assume that $x\in \cA$. Let us show that
$\delta_0(x)=0$. There exists $\phi\in \cA$ with $x= x \phi^2$ and
$\phi \xi=0$. One has $\delta_0(x)=i(\xi|[D,x]\xi)$ and $[D,x]=
[D,x]\phi^2+2x  [D,\phi]\phi$ so that $[D,x]\xi=0$ and
$\delta_0(x)=0$. It follows that for $f$ as in \eqref{smearing}, one
gets $\delta_0(U(f)(x))=0$ since $U(f)$ commutes with $\delta_0$.
With $f$ analytic for $L^1$ one gets that $U(f)(x)$ is an analytic
element such that $\delta_0(U(f)(x))=0$ and hence
$\delta_0^n(U(f)(x))=0$ for all $n\geq 1$. It follows that
$$
U(t)(U(f)(x))=\sum \frac{t^n}{n!}\delta_0^n(U(f)(x))=U(f)(x)\qqq
t\in \R
$$
Thus since $U(f_n)(x)\to x$ in norm for a suitable sequence $f_n$,
one gets that $U(t)(x)=x$ for all $t\in \R$.\endproof

\begin{lem} \label{restrictedop} Let $S$ be the closure of
$S_{\xi,\,\eta}=D_{\xi,\,\eta}-(\xi|D\xi)T_{\eta,\,\eta}$ as an
unbounded operator in $\cH$. Then for any $a\in \Dom \delta_0$ one
has $a\eta\in \Dom S$ and
\begin{equation}\label{restrictedop1}
    S(a\eta)=-i\delta_0(a)\eta
\end{equation}
For any $a\in A$ and $\epsilon >0$, let $b=(1+\epsilon
\delta_0)^{-1}(a)$. Then $b\eta \in \Dom S$ and
\begin{equation}\label{restrictedop2}
    (1+i\epsilon S)(b\eta)=a\eta
\end{equation}
\end{lem}

\proof For $a_n\to a$ and $\delta_0(a_n)\to \delta_0(a)$ in norm one
has $a_n \eta\to a\eta$ and $\delta_0(a_n)\eta\to \delta_0(a)\eta$
in $\cH$. Thus, since $S$ is closed and $\cA$ is a core for
$\delta_0$ it is enough to prove \eqref{restrictedop1} for $a\in
\cA$. In that case one gets
$$
 S(a\eta)=S_{\xi,\,\eta}(a\eta)=\left (
 (\xi|D(a(\eta|\eta)\xi))-(\xi|D\xi)(\eta|a\eta)\right ) \eta
$$
One has $a(\eta|\eta)\xi=a\xi$ since $(\eta|\eta)=1$ on the support
of $\xi$. Similarly
$(\xi|D\xi)(\eta|a\eta)=(\xi|D\xi)a=(\xi|aD\xi)$. Thus we get
$$
S(a\eta)=(\xi|[D ,a]\xi) \qqq a \in \cA\,,
$$
which gives \eqref{restrictedop1}.

To prove \eqref{restrictedop2}, note that by Theorem
\ref{hileyosidathm} the resolvent $(1+\epsilon \delta_0)^{-1}$
exists for any $\epsilon >0$ and maps $A$ to the domain of
$\delta_0$. Thus applying the first part of the Lemma to
$b=(1+\epsilon \delta_0)^{-1}(a)$, one gets
$$
(1+i\epsilon S)(b\eta)=b\eta + i\epsilon
(-i\delta_0(b)\eta)=((1+\epsilon \delta_0)b)\eta=a\eta
$$
which gives \eqref{restrictedop2}.
\endproof

\begin{lem} \label{respectsobolevlem} The one parameter group
$U(t)\in \Aut(A)$ fulfills for each $N$ an estimate of the form
\begin{equation}\label{restrictedop3}
   \Vert U(t)(a)\Vert_N^{\rm
    sobolev}\leq c_1 e^{Nc_{\xi,\,\eta}|t|   } \Vert a\Vert_N^{\rm
    sobolev}
\end{equation}
\end{lem}

\proof We first use  the Sobolev semi-norm   given by
$$
\Vert a\Vert_{N,\,\eta,\,\lambda}^{\rm
    sobolev}=(D^2+\lambda)^{N/2} a \eta
    $$
    with $\lambda>0$ determined by Lemma \ref{resolvequlem1}. We let $\epsilon_0>0$ be as in
    Lemma \ref{resolvequlem1}. By
    \eqref{restrictedop2} one has for any $a\in A$
\begin{equation}\label{coincide}
    (1+\epsilon \delta_0)^{-1}(a)\eta=(1+i\epsilon S)^{-1}a\eta \qqq
    \epsilon\leq \epsilon_0
\end{equation}
Assume that $\Vert a\Vert_{N,\,\eta,\,\lambda}^{\rm
    sobolev}<\infty$. One then has $$a\eta=(D^2+\lambda)^{-N/2}\zeta
    \,, \
    \zeta=(D^2+\lambda)^{N/2} a \eta\in \cH$$
     By Lemma \ref{resolvequlem1}
    one gets, for $\epsilon\leq \epsilon_0$, using \eqref{opinHbis},
    $$
\Vert(1+\epsilon \delta_0)^{-1}(a)\Vert_{N,\,\eta,\,\lambda}^{\rm
    sobolev}=\Vert (D^2+\lambda)^{N/2}(1+\epsilon
    \delta_0)^{-1}(a)\eta\Vert=\Vert (D^2+\lambda)^{N/2}(1+i\epsilon
    S)^{-1}a\eta\Vert$$
    $$
    =\Vert \theta_\lambda^{N/2}((1+i\epsilon
    S)^{-1})\zeta\Vert\leq (1+N c_{\xi,\,\eta}\epsilon)\Vert\zeta\Vert=
    (1+N c_{\xi,\,\eta}\epsilon)\Vert a\Vert_{N,\,\eta,\,\lambda}^{\rm
    sobolev}
    $$
    This shows that $\Vert(1+\epsilon \delta_0)^{-1}(a)\Vert_{N,\,\eta,\,\lambda}^{\rm
    sobolev}<\infty$ and thus one can iterate and get,
    \begin{equation}\label{iterated}
\Vert(1+\epsilon \delta_0)^{-m}(a)\Vert_{N,\,\eta,\,\lambda}^{\rm
    sobolev}\leq (1+N c_{\xi,\,\eta}\epsilon)^m\Vert a\Vert_{N,\,\eta,\,\lambda}^{\rm
    sobolev}  \qqq
    \epsilon\leq \epsilon_0
    \end{equation}
     Now for $t>0$ and with norm convergence in $A$ one has, $$U(-t)a=\lim_{n\to\infty} (1+\frac
    {t\delta_0}{n})^{-n}(a)\,.$$
    This shows that $U(-t)a$ is the norm limit of the sequence $a_n=(1+\frac
    {t\delta_0}{n})^{-n}(a)$ and moreover one has, from
    \eqref{iterated},
    $\Vert a_n\Vert_{N,\,\eta,\,\lambda}^{\rm
    sobolev}\leq (1+ N c_{\xi,\,\eta}t/n)^n\Vert a\Vert_{N,\,\eta,\,\lambda}^{\rm
    sobolev}$. Thus
    \begin{equation}\label{iterated1}
\limsup \Vert a_n\Vert_{N,\,\eta,\,\lambda}^{\rm
    sobolev}\leq e^{N c_{\xi,\,\eta}|t|}\Vert a\Vert_{N,\,\eta,\,\lambda}^{\rm
    sobolev}
    \end{equation}
    Since $a_n\to b=U(-t)a$ in norm, one has $a_n\eta\to
    b\eta$ also in norm in $\cH$. Since the operator
    $(D^2+\lambda)^{N/2}$ is closed, and by \eqref{iterated1} the $(D^2+\lambda)^{N/2}a_n\eta$
    are uniformly bounded, it follows that $b\eta\in \Dom
    (D^2+\lambda)^{N/2}$ and thus $\Vert U(-t)a\Vert_{N,\,\eta,\,\lambda}^{\rm
    sobolev}<\infty$. More precisely we get
    $$
\Vert U(-t)a\Vert_{N,\,\eta,\,\lambda}^{\rm
    sobolev}\leq e^{N c_{\xi,\,\eta}|t|}\Vert a\Vert_{N,\,\eta,\,\lambda}^{\rm
    sobolev}
    $$
    Now the semi-norm $\Vert a\Vert_{N,\,\eta,\,\lambda}^{\rm
    sobolev}$ is not equivalent to the Sobolev norm but the latter
    is equivalent to the sum
    $$
    \Vert(D^2+\lambda)^{N/2} a \eta\Vert + \sum \Vert(D^2+\lambda)^{N/2} a \eta_\mu\Vert
    $$
where one can choose the $\eta_\mu$ so that their supports are
disjoint from the support of $\xi$. This can be seen using the
strong regularity. It then follows from Lemma \ref{supportlem}, that
the semi-norm $\sum \Vert(D^2+\lambda)^{N/2} a \eta_\mu\Vert$ is
preserved by $U(t)$ since $U(t)(a)\eta_\mu=a\eta_\mu$ for all $\mu$.
Thus one obtains \eqref{restrictedop3}. \endproof

\begin{thm} \label{expabilitythm} Let $(\cA,\cH,D)$ be a strongly regular spectral triple with $\cA$
commutative, fulfilling the five \axioms of \S \ref{prelem}. Then any
derivation of $\cA$ of the form \eqref{formforder} \ie
$\delta_0(a)=i(\xi|[D,a]\xi)$, $ \forall a \in \cA$, is the
generator of a one-parameter group of automorphisms $\sigma_t\in
\Aut(\cA)$ such that
\begin{itemize}
  \item $\partial_t \sigma_t(a)=\delta_0(\sigma_t(a))$.
  \item The map $(t,a)\in \R\times \cA\mapsto \sigma_t(a)\in \cA$ is
  jointly continuous.
 \end{itemize}
\end{thm}

\proof By Lemma \ref{respectsobolevlem}, the one parameter group
$U(t)\in \Aut(A)$ preserves the subalgebra $\cA\subset A$. We let
$\sigma_t\in \Aut(\cA)$ be the corresponding automorphisms. For
$a\in \cA$ one has $a\in \Dom \delta_0$ and thus
$$
  \sigma_t(a)-a=  \int_0^t\sigma_u(\delta_0(a))du
$$
where $\sigma_u(\delta_0(a))$ is a norm continuous function of $u$.
By \eqref{restrictedop3} applied to $\delta_0(a)$,  this shows that
$\Vert \sigma_t(a)-a\Vert_N^{\rm
    sobolev}=O(|t|)$ when $t\to 0$. One has
$$
\frac 1t (\sigma_t(a)-a)-\delta_0(a)=\frac 1t
\int_0^t(\sigma_u(\delta_0(a))-\delta_0(a))du
$$
which, since $\Vert \sigma_u(\delta_0(a))-\delta_0(a)\Vert_N^{\rm
    sobolev}=O(|u|)$, gives
$$
\Vert \frac 1t (\sigma_t(a)-a)-\delta_0(a)\Vert_N^{\rm
    sobolev}=O(|t|) \,, \ {\rm for }\  \ t\to 0\,.
$$
This shows that $\partial_t \sigma_t(a)=\delta_0(\sigma_t(a)$ in the
Frechet space $\cA$. Let us check the joint continuity of
$(t,a)\mapsto \sigma_t(a)$. Let $(t_n,a_n)\to (t,a)\in \R\times
\cA$. One has
$\sigma_{t_n}(a_n)-\sigma_t(a)=\sigma_{t_n}(a_n-a)+\sigma_{t_n}(a)-\sigma_t(a)$.
The norm $\Vert \sigma_{t_n}(a)-\sigma_t(a)\Vert_N^{\rm
    sobolev}$  converges to $0$ from the above discussion. Moreover Lemma \ref{respectsobolevlem}
    shows that one controls the Sobolev norms of
    $\sigma_{t_n}(a_n-a)$ from those of $(a_n-a)$ which gives the required continuity. \endproof

We can now also prove directly the absolute continuity of the
transformed measure $\sigma_t^*(\lambda)$ with respect to $\lambda$.

\begin{prop}\label{abscontdirect} Let $(\cA,\cH,D)$, $\delta_0$ and
$\sigma_t$ be as in Theorem \ref{expabilitythm}. Then for each $t\in
\R$ the measure $\lambda$ of \eqref{measurelambdadef} is
strongly\footnote{$\mu$ is strongly equivalent to $\nu$ iff there is
$c>0$ with $c\nu\leq\mu\leq c^{-1}\nu$} equivalent to its transform
under $\sigma_t$.
\end{prop}

\proof Let $\delta_0(a)=i(\xi|[D,a]\xi)$, $ \forall a \in \cA$. By
Lemma \ref{supportlem} the measure $\sigma_t^*(\lambda)$, given by
$\sigma_t^*(\lambda)(f)=\lambda(\sigma_t(f))$ agrees with
$\lambda(f)$ whenever the support of $f$ is disjoint from the
support of $\xi$. With $\eta\in \cH_\infty$ as above one has
$(\eta|\eta)=1$ in a neighborhood $V$ of the support of $\xi$. To
obtain the required strong equivalence, it is enough to compare
$\lambda(\sigma_t(f))$ and $\lambda(f)$ for $f$  and $\sigma_t(f)$
with support contained in $V$. One then has, using \eqref{absocont}
$$
\lambda(\sigma_t(f))=\cutint \sigma_t(f)|D|^{-p}=\cutint
\sigma_t(f)(\eta|\eta)|D|^{-p}=\langle\eta,\sigma_t(f)\eta\rangle
$$
Let, as above,  $S$ be the closure of
$S_{\xi,\,\eta}=D_{\xi,\,\eta}-(\xi|D\xi)T_{\eta,\,\eta}$. It is by
construction a bounded perturbation of the self-adjoint operator
(closure of) $D_{\xi,\,\eta}$ and one can define $e^{itS}$ for $t\in
\R$ using the expansional formula (\cite{Araki})
\begin{equation}\label{expansional}
e^{A+B}e^{-A}=\left( \sum_n \int_{S_n} \alpha_{u_1}(B)\cdots
  \alpha_{u_n}(B) du \right)\,,\ \ \alpha_u(B)=e^{uA}Be^{-uA}
\end{equation}
with $A=it D_{\xi,\,\eta}$ and $B=-it(\xi|D\xi)T_{\eta,\,\eta}$. Let
us show that
\begin{equation}\label{sigmatandt}
\sigma_t(a)\eta=e^{itS}a\eta \qqq a\in \cA\,.
\end{equation}
By Theorem \ref{expabilitythm} and \eqref{restrictedop1} the
$\cH$-valued function $t\mapsto \eta(t)=\sigma_t(a)\eta$ solves the
differential equation
$$
\frac{d\eta(t)}{dt}=iS\eta(t)\,,\ \eta(0)=a\eta\,, \ \eta(t)\in\Dom
S\qqq t\in \R\,.
$$
This implies that $\frac{d}{dt}(e^{-itS}\eta(t))=0$ and thus
$e^{-itS}\eta(t)=a\eta$ which proves \eqref{sigmatandt}. It follows
from \eqref{sigmatandt} that
\begin{equation}\label{tminuststar0}
\langle\eta,\sigma_t(a)\eta\rangle=\langle\eta,e^{itS}a\eta\rangle=\langle
e^{-itS^*}\eta,a\eta\rangle
\end{equation}

\medskip

Note that $S$ is not self-adjoint in general because of the
additional term $-(\xi|D\xi)T_{\eta,\,\eta}$.  The difference
$S-S^*$ is a bounded operator and an endomorphism of the
$\cA$-module $\cH_\infty$ given by
\begin{equation}\label{tminuststar}
   S-S^*=\rho \,T_{\eta,\,\eta}\,, \ \  \rho=(\xi|D\xi)^*-(\xi|D\xi)
\end{equation}
since $T_{\eta,\,\eta}$ is self-adjoint  by \eqref{rankonecstar2}.
We can now write a formula for $e^{-itS^*}\eta$,
\begin{equation}\label{tminuststar1}
  e^{-itS^*}\eta=\left( \sum_n i^nt^n\int_{S_n} \sigma_{-tu_1}(\rho)\cdots
  \sigma_{-tu_n}(\rho) du\right)\,\eta
\end{equation}
with $S_n=\{(u_j)\,|\,0\leq u_1\leq\ldots\leq u_n\leq 1\}$ the
standard simplex. Indeed one has $-itS^*=-itS+P$ with $P=it\rho
T_{\eta,\,\eta}$ which is bounded which allows one to use the
expansional formula \eqref{expansional}, with $A=-itS$, $B=P$. Now
by \eqref{sigmatandt} one has $e^{itS}\eta=\eta$ for all $t\in \R$
thus the left hand side of \eqref{expansional} applied to $\eta$
gives $e^{-itS^*}\eta$. Let us compute the right hand side. We first
show that
\begin{equation}\label{induct}
e^{isS}\rho T_{\eta,\,\eta} a\eta=\sigma_s(\rho a)\eta \qqq  a\in
\cA\,.
\end{equation}
Indeed one has $T_{\eta,\,\eta}
a\eta=(\eta|a\eta)\eta=a(\eta|\eta)\eta$ and since $(\eta|\eta)=1$
on the support of $\rho$ (using \eqref{tminuststar}) one gets that
$\rho T_{\eta,\,\eta} a\eta=\rho a(\eta|\eta)\eta=\rho a \eta$. Thus
\eqref{induct} follows from \eqref{sigmatandt}. We then get
$$
\alpha_{u_1}(P)\cdots
  \alpha_{u_n}(P)\eta=e^{-itu_1S}Pe^{-it(u_2-u_1)S}P\cdots
  e^{-it(u_n-u_{n-1})S}P\eta=
$$
$$
i^nt^n e^{-itu_1S}\rho T_{\eta,\,\eta} e^{-it(u_2-u_1)S}\rho
T_{\eta,\,\eta}\cdots
  e^{-it(u_n-u_{n-1})S}\rho T_{\eta,\,\eta}\eta=
$$
$$
=i^nt^n \sigma_{-tu_1}(\rho \sigma_{-t(u_2-u_1)}(\rho( \ldots (\rho
\sigma_{-t(u_n-u_{n-1})}(\rho))\ldots )\eta
$$
which yields \eqref{tminuststar1} from \eqref{expansional}. Now the
series
\begin{equation}\label{radonnik}
    h(t)=\sum_n i^nt^n\int_{S_n} \sigma_{-tu_1}(\rho)\cdots
  \sigma_{-tu_n}(\rho) du
\end{equation}
converges in the Frechet algebra $\cA$ since, for each $k$, the
$p_k(\sigma_s(\rho))$ are uniformly bounded on compact sets of $s$,
while the volume of the simplex $S_n$ is $1/n!$. Thus $h(t)\in \cA$
and combining \eqref{tminuststar0} and \eqref{tminuststar1} one has:
\begin{equation}\label{radonnik1}
\langle\eta,\sigma_t(f)\eta\rangle=\langle
e^{-itS^*}\eta,f\eta\rangle=\langle
h(t)\eta,f\eta\rangle=\langle\eta,\bar h(t) f\eta\rangle
\end{equation}
so that we get, for all $f$ with support in $V$,
\begin{equation}\label{radonnik2}
\lambda(\sigma_t(f))=\lambda(\bar h(t) f)
\end{equation}
Since, by construction one has $h=1$ outside the support of $\xi$,
(using Lemma \ref{supportlem}), equality \eqref{radonnik2} holds for
all $f\in \cA$. The norm continuity $\Vert h(t)-1\Vert\to 0$ when
$t\to 0$ (using \eqref{radonnik}) then gives the required strong
equivalence.
\endproof

\medskip \section{Absolute continuity}\label{sectabscont}

The following equality defines a positive measure $\lambda$ on $X$.
\begin{equation} \label{lebesgue}
\int ad\lambda=\cutint  \, a|D|^{-p} \qqq a\in C(X)\,.
\end{equation}
This measure is locally equivalent to the spectral measure of the
representation of $A=C(X)$ in $\cH$. More precisely:

\medskip
\begin{lem} \label{onedlem}  For any open set $V\subset X$
 the following two measures are strongly equivalent:
\begin{itemize}
    \item The restriction $\lambda|_V$ to $V$ of the measure
  $\lambda$ of \eqref{lebesgue}.
  \item The restriction to $V$ of the spectral measure associated to
  a vector $\xi\in \cH^\infty$ whose $\cA$-valued inner product
  $(\xi,\xi)$ is strictly positive on $\overline V$.
\end{itemize}
\end{lem}

\proof By the \axiom of absolute continuity one has a relation of the
form
$$
\langle \xi, a \xi\rangle= \cutint  \, a  (\xi,\xi) |D|^{-p}
$$
and since $(\xi,\xi)\in \cA$ is strictly positive on $\overline V$,
one gets the strong equivalence between  the restriction to $V$ of
the spectral measure associated to
 the vector $\xi\in \cH^\infty$ and the measure $\lambda|_V$ of \eqref{lebesgue}.
\endproof

We let
$$B_\epsilon=\{t\in \R^p\,|\,|t|<\epsilon\}$$
Given an automorphism $\sigma\in \Aut(\cA)$ we use the covariant
notation
\begin{equation}\label{covariant}
    \sigma(\kappa)=\kappa\circ \sigma^{-1}\qqq \kappa\in \Sp(\cA)
\end{equation}
and view $\sigma$ as an homeomorphism of $X=\Sp(\cA)$. We use the notations
$U_\alpha$, $s_\alpha$ of Lemma \ref{coverlem} and of \eqref{defnsalpha}.

\medskip

\begin{lem} \label{specmeas0} Let $V\subset U_\alpha$ be an open set
and $\chi \in V$. There exists a smooth family $\sigma_t\in
\Aut(\cA)$, $t\in \R^p$, a neighborhood $Z$ of $\chi$ in $V$ and
$\epsilon>0$, $\epsilon'>0$ such that:
\begin{enumerate}
  \item For any $\kappa\in Z$, the map $t\mapsto
s_\alpha(\sigma_t(\kappa))=F(\kappa,t)$ is a diffeomorphism,
depending continuously on $\kappa$, of $B_\epsilon$ with a
neighborhood of $s_\alpha(\kappa)$ in $\R^p$ and
\begin{equation}\label{surject}
    s_\alpha(\kappa)+B_{\epsilon'}\subset F(\kappa,B_{\epsilon/2})\qqq \kappa\in \overline Z
\end{equation}
\item For any $t\in B_\epsilon$ one has
\begin{equation}\label{quasi}
 \frac 12\,  \lambda \leq \sigma_t(\lambda)\leq  2 \lambda \,.
\end{equation}
  \item $Z_1=\cap_{B_\epsilon}\sigma_tZ$ is a neighborhood of $\chi$.
  \item $Z_2=\cup_{B_\epsilon}\sigma_tZ$ is contained in $V$.
\end{enumerate}
\end{lem}

\proof Let  $\sigma_t\in \Aut(\cA)$, $t\in \R^p$, $W$ and $Z$ as in
Lemma \ref{coropen}.  We can replace the $Z_0$ of Lemma
\ref{coropen} by any neighborhood of $\chi$ contained in $Z_0$ and
hence by a ball centered at $\chi$ and contained in $V\cap Z_0$. We
use a metric $d$ on $X$ compatible with the topology (Proposition
\ref{generaltop}). Thus
$$
Z=\{\kappa\in X\,|\,d(\kappa,\chi)< r\}
$$
and we can take $r$ small enough so that
\begin{equation}\label{surjectbis1}
\{\kappa\in X\,|\,d(\kappa,\chi)\leq 3/2\,r\}\subset V
\end{equation}
The continuity of the map $(\kappa,t)\mapsto
\sigma_t^{-1}(\kappa)=\kappa\circ \sigma_t$ yields $\epsilon >0$
with $B_\epsilon\subset W$ and
\begin{equation}\label{surjectbis2}
d(\kappa, \sigma_t^{\pm 1}(\kappa))\leq r/2\qqq \kappa\in X\,, \
t\in B_\epsilon
\end{equation}
Then the first statement (1) follows from Lemma \ref{coropen}, with
\eqref{surject} coming from the continuity in $\kappa$. The second
statement follows from \eqref{surjectbis2} since for
$d(\kappa,\chi)< r/2$ one gets $d(\chi, \sigma_t^{-1}(\kappa))< r$
and $\sigma_t^{-1}(\kappa)\in Z$. Similarly the third statement
follows from \eqref{surjectbis2} and \eqref{surjectbis1}. Finally
\eqref{quasi} follows from Proposition \ref{abscontdirect} for
$\epsilon$ small.\endproof

\begin{lem} \label{specmeas} Let $V\subset X$ be an open set with $\overline V\subset U_\alpha$
and $\lambda_V$ (resp. $\lambda_{\overline V}$) be the spectral
measure of the restriction to $V$ (resp. $\overline V$) of the
representation of $C(X)$ in $\cH$. Then $s_\alpha(\lambda_V)$ is
equivalent to the Lebesgue measure on $s_\alpha(V)$ and there exists
$c<\infty$ such that
\begin{equation}\label{specmeasmaj}
   \int_{\overline V} f\circ s_\alpha\,d\lambda_{\overline V}\leq c
   \, \int_{s_\alpha({\overline V})} f(x) dx^p \qqq f\in C_c^+(\R^p)\,.
\end{equation}
\end{lem}

\proof  By Lemma \ref{onedlem}, the spectral measure $\lambda_V$
(resp. $\lambda_{\overline V}$) is equivalent to the measure
$\lambda$ of \eqref{lebesgue} restricted to $V$ (resp. $\overline
V$). We show that
\begin{itemize}
  \item For any $\chi\in \overline V$ one can find a
neighborhood $Z_1$ of $\chi$ in $U_\alpha$ such that
$s_\alpha(\lambda|_{Z_1})\leq c\, dx^p$ for some $c<\infty$.
  \item For any $\chi \in V$ one can find a
neighborhood $Z_2$ of $\chi$ in $V$ such that
$ds_\alpha(\lambda|_{Z_2})/dx^p=\rho(x)>0$  in a neighborhood of
$s_\alpha(\chi)$
\end{itemize}
Let $\chi\in \overline V$. We apply Lemma \ref{specmeas0} relative
to $V=U_\alpha$. We let $\sigma_t$, $Z$, $Z_j$, $\epsilon$ and
$\epsilon'$ be as in Lemma \ref{specmeas0}.
  We can assume that for $|t|<\epsilon$ one
has \eqref{quasi}. Let then $h\in C_c^\infty(B_\epsilon)$, $h(t)\in
[0,1]$, be equal to $1$ on $B_{\epsilon/2}$. By Lemma
\ref{specmeas0}, for any $\kappa\in Z$, the map $t\mapsto
F(\kappa,t)=s_\alpha(\sigma_t(\kappa))$ is a diffeomorphism
$F_\kappa$
 of $B_\epsilon$ with a
neighborhood of $s_\alpha(\kappa)$ in $\R^p$. It follows then
 that for fixed $\kappa$ the image in $\R^p$ of the measure
$h(t)dt^p$, is a smooth multiple $g_\kappa(u)$ of the Lebesgue
measure $du^p$,
\begin{equation}\label{surject1}
\int_{B_\epsilon} f(F(\kappa,t))h(t)dt^p=\int_{\R^p} f(u)
g_\kappa(u) du^p \qqq f\in C_c(\R^p)\,.
\end{equation}
The function $g_\kappa$ vanishes outside $F_\kappa(B_\epsilon)$ and
is given inside by
\begin{equation}\label{jacobian}
g_\kappa(u)=h(\psi(u))|d\psi(u)/du|
\end{equation}
 where $\psi$ is the inverse
of the diffeomorphism $F_\kappa$ and $d\psi(u)/du$ its Jacobian. The
continuity of the map $\kappa\mapsto F_\kappa$ gives a uniform upper
bound
\begin{equation}\label{jacobian1}
g_\kappa(u)\leq c_1 \qqq u\in \R^p \,,\ \kappa \in Z\,.
\end{equation}
Since $h=1$ on $B_{\epsilon/2}$ and
$s_\alpha(\kappa)+B_{\epsilon'}\subset F_\kappa(B_{\epsilon/2})$ by
\eqref{surject}, one has $h(\psi(u))=1$ for $u\in
s_\alpha(y)+B_{\epsilon'}$. The continuity of the map $\kappa\mapsto
F_\kappa$ then yields $\epsilon_1>0$ such that
\begin{equation}\label{surject2}
g_\kappa(u)\geq \epsilon_1 \qqq u\in
s_\alpha(\kappa)+B_{\epsilon'}\qqq \kappa\in Z
\end{equation}
 We
consider the image $d\nu$ under $(\kappa,t)\in Z\times
B_\epsilon\mapsto F(\kappa,t) \in \R^p$ of the finite positive
measure $ d\lambda(\kappa)h(t)dt^p  $ on $Z\times B_\epsilon$. It is
given by
\begin{equation}\label{surjectmeas}
\int_{\R^p} f(x)
d\nu(x)=\int_Z\int_{B_\epsilon}f(F(\kappa,t))h(t)dt^p
d\lambda(\kappa)\qqq f\in C_c(\R^p)\,,
\end{equation}
and is equal, by \eqref{surject1}, to
\begin{equation}\label{dnu}
    \int f(x) d\nu(x)=\int \int f(u) g_\kappa(u) du^p d\lambda(\kappa)=\int
f(u)\rho(u)du^p
\end{equation}
where
\begin{equation}\label{dnu1}
\rho(u)=\int_Z  g_\kappa(u)   d\lambda(\kappa)
\end{equation}
By \eqref{jacobian1} one has
\begin{equation}\label{jacobian2}
\rho(u)\leq c_1 \lambda(Z)<\infty\qqq u\in \R^p \,.
\end{equation}
Moreover \eqref{surject2} shows that
\begin{equation}\label{surject3}
    \rho(u)\geq \epsilon_1 \lambda(\{\kappa\in Z\,|\,|u-s_\alpha(\kappa)|< \epsilon'\})\,.
\end{equation}
We then have
\begin{equation}\label{surject4}
\rho(u)>0\qqq
    u\in s_\alpha(Z)\,.
\end{equation}
This strict positivity follows from the \axiom of absolute continuity
which shows that the support of the measure $\lambda$ is $X$.
Indeed, for $u\in s_\alpha(Z)$,  the open set $\{\kappa\in
Z\,|\,|u-s_\alpha(\kappa)|< \epsilon'\}$ is non-empty and  it has
strictly positive measure. This  shows that the restriction of the
measure $\nu$ to the open set $s_\alpha(Z)$  is equivalent to the
Lebesgue measure.

\smallskip

We now use the quasi-invariance of $d\lambda$ given by \eqref{quasi}
to compare $s_\alpha(\lambda|_{Z_j})$ with $\nu$. Using
\eqref{quasi} (for $\delta=1/2$), one has for $t\in B_\epsilon$, $
\frac 1 2 d\lambda \leq d (\sigma_t(\lambda))\leq 2 d\lambda $ so
that, for any subset $E\subset X$ and any positive $f\in
C_c^+(\R^p)$ one has, with $1_E$ the characteristic function of $E$,
$$
\frac 1 2\int (f\circ s_\alpha)1_E  d\lambda\leq \int  (f\circ
s_\alpha) 1_E d (\sigma_t(\lambda))\leq 2\int  (f\circ s_\alpha) 1_E
d\lambda
$$
The middle term is
$$
\int (f\circ s_\alpha)1_E d (\sigma_t(\lambda))= \int (f\circ
s_\alpha\circ \sigma_t)(1_E\circ \sigma_t) d \lambda
$$
and we thus get
\begin{equation}\label{surject5}
    \frac 1 2\int_E (f\circ s_\alpha)  d\lambda\leq
    \int_{\sigma_t^{-1}E}
(f\circ s_\alpha\circ \sigma_t) d \lambda \leq 2 \int_E (f\circ
s_\alpha)  d\lambda\qqq f\in C_c^+(\R^p)
\end{equation}

We let, as in Lemma \ref{specmeas0},
$$
Z_1=\cap_{B_\epsilon}\sigma_tZ\,, \ \
Z_2=\cup_{B_\epsilon}\sigma_tZ\,.
$$
One has $\sigma_t^{-1}(Z_1)\subset Z$ for $t\in B_\epsilon$ so that,
by the first inequality of \eqref{surject5} for $E=Z_1$,
$$
\frac 1 2\int_{Z_1} (f\circ s_\alpha) d\lambda\leq
\int_{\sigma_t^{-1}Z_1} (f\circ s_\alpha\circ \sigma_t) d
\lambda\leq \int_Z (f\circ s_\alpha\circ \sigma_t)d \lambda \qqq
t\in B_\epsilon\,, \ f\in C_c^+(\R^p)
$$
so that, multiplying by $h(t)dt^p$ and  integrating over $t\in
B_\epsilon$, we get $C<\infty$ with
$$
\int_{Z_1} (f\circ s_\alpha) d\lambda\leq  C \int_Z\int_{B_\epsilon}
f(s_\alpha(\sigma_t(\kappa))h(t)dt^p d\lambda(\kappa)=C\int_{\R^p}
f(x) d\nu(x)\,,
$$
where we used Fubini's theorem, and the equality
$s_\alpha(\sigma_t(\kappa))=F(\kappa,t)$ for $\kappa\in Z$ and $t\in
B_\epsilon$. Thus, using \eqref{dnu} and \eqref{jacobian2}
\begin{equation}\label{surject6}
\int_{Z_1} (f\circ s_\alpha) d\lambda\leq  C \int_{\R^p}
f(u)\rho(u)du^p \leq  C' \int_{\R^p} f(u)du^p\qqq f\in C_c^+(\R^p)
\end{equation}
hence the image $s_\alpha(\lambda|_{Z_1})$ is absolutely continuous
with respect to the Lebesgue measure and is majorized by a constant
multiple of Lebesgue measure. Thus,
 every point of $\overline V$ has a neighborhood $Z_1$ such that
$s_\alpha(\lambda|_{Z_1})\leq c_1 dx^p$. Covering the compact set
$\overline V$ by finitely many such $Z_1$ gives \eqref{specmeasmaj}.

\smallskip

Let us now assume that $\chi\in V$. We  can then assume, by Lemma
\ref{specmeas0} that $Z_2=\cup_{B_\epsilon}\sigma_tZ$ is contained
in $V$. One has $Z\subset \sigma_t^{-1}(Z_2) $ for $t\in B_\epsilon$
so that, by the second inequality of \eqref{surject5} for $E=Z_2$,
$$
\int_Z (f\circ s_\alpha\circ \sigma_t) d \lambda\leq
\int_{\sigma_t^{-1}(Z_2)} (f\circ s_\alpha\circ \sigma_t) d \lambda
\leq 2\int_{Z_2} (f\circ s_\alpha) d\lambda\qqq t\in B_\epsilon
$$
thus, after integration over $t\in B_\epsilon$,
$$
C'\int_{\R^p} f(x) d\nu(x)=C'\int_Z\int_{B_\epsilon}
f(s_\alpha(\sigma_t(y))h(t)dt^p d\lambda(y)\leq \int_{Z_2} (f\circ
s_\alpha) d\lambda
$$
This shows, using \eqref{dnu}, that
\begin{equation}\label{surject7}
   \int_{Z_2} (f\circ
s_\alpha) d\lambda\geq C'\int_{\R^p} f(u)\rho(u)du^p \qqq f\in
C_c^+(\R^p)
\end{equation}
 By \eqref{surject4} one has
$\rho(u)>0$ for all
    $u\in s_\alpha(Z)$ thus $\rho(u)>0$ in a neighborhood of  $s_\alpha(\chi)$, in other words
    $ds_\alpha(\lambda|_{Z_2})/dx^p=\rho_2(x)>0$  in a neighborhood of $s_\alpha(\chi)$ as required.
    This shows that $s_\alpha(\lambda|_V)$ is equivalent to the Lebesgue measure on $s_\alpha(V)$.
\endproof

\medskip

\section{Spectral multiplicity}\label{sectspecmul}

 We want to get an upper bound for the number of
elements in the fiber of the map $s_\alpha:U_\alpha\to \R^p$. We
shall first
  relate the multiplicity of the map $s_\alpha$ with the spectral
multiplicity of the operators $a_\alpha^j$ in the Hilbert space
$\cH$. This is not automatic, indeed the first difficulty is that
for an injective representation $\pi$ of a $C^*$-algebra $B$ with a
subalgebra $A\subset B$, one can have the same double-commutants
 $\pi(A)''=\pi(B)''$ even though
$A\neq B$. Thus for instance one can take the subalgebra
$C[0,1]\subset C(K)$ where $K=\{0,1,\ldots,9\}^\N$ is the Cantor set
of the decimal digits and the inclusion is given by the decimal
expansion. Both act in $L^2[0,1]$ (by multiplication) and the
spectral multiplicity of the function $x\in C[0,1]$ is equal to one
but the number of elements in the fiber is equal to $2$ for numbers
of the form $k \,10^{-n}$. The point  in this example is that the
projection map $s:K\to [0,1]$ is not an open mapping. Thus in
particular the subset of $K$ where the multiplicity of $s$ is two is
not an open subset of $K$ (it is countable).

\begin{lem} \label{semicont} Let $s:X\to Y$ be a continuous open map of compact
spaces. Then the function $n(y)=\#s^{-1}(y)$ is lower semicontinuous on
$Y$.
\end{lem}

\proof Assume that $n(y)\geq m$ and let us show that this inequality
still holds in a neighborhood of $y$. Let $x_j\in X$ be $m$ distinct
points in $s^{-1}(y)$. One can then find disjoint open sets $V_j\ni
x_j$ and let $W=\cap_j\, s(V_j)$ which is an open neighborhood of
$y$. For any $z\in W$ the preimage $s^{-1}(z)$ contains at least $m$
points.
\endproof

Now, let $s:X\to Y$ be a continuous open map of compact spaces. Let
$\mu$ be a positive measure on $X$ with support $X$ and $\pi$ the
corresponding representation of $C(X)$ in $L^2(X,\mu)$. We want to
compare the  spectral multiplicity function $\Sigma(y)$ of the
restriction of $\pi$ to $C(Y)$ with  $n(y)=\#s^{-1}(y)$. Let
$\nu=s(\mu)$ be the image of the measure $\mu$. One can desintegrate
$\mu$ in the form
$$
\mu=\int_Y \rho_y\,d\nu(y)
$$
where the conditional measure $\rho_y$ is supported by the closed
subset $s^{-1}(y)$. The issue is what is the dimension of the
Hilbert space $L^2(X,\rho_y)$. It might seem at first that if the
support of the measure $\mu$ is $X$ one should be able to conclude
that the support of $\rho_y$ is $s^{-1}(y)$ and obtain that the
spectral multiplicity function $\Sigma(y)$ is larger than
$n(y)=\#s^{-1}(y)$.
 However this fails as shown by
  the following example:
$$
X=Y\times \{1,\ldots,m\}\,,\ \ s(y,k)=y
$$
and let $\mu_k$ be the measure on $Y$ corresponding to the
restriction of $\mu$ to $Y\times \{k\}$.

\begin{lem}\label{singular}
If the measures $\mu_k$ are mutually singular then the spectral multiplicity
function $\Sigma(y)$ is equal to $1$ \aee.
\end{lem}

\proof The representations of $C(Y)$ in $L^2(Y,\mu_k)$ are pairwise
disjoint, and each is of multiplicity one. Thus the commutant of
$C(Y)$ in the direct sum of these representations only contains
block diagonal operators and is hence commutative so that the
multiplicity is equal to one.
\endproof

The above example gives the needed condition for the relation
between $\Sigma(y)$ and $n(y)$ and one has:

\begin{lem} \label{lowermulti}Let $X$ be a compact space and $\lambda$ a
finite positive measure on $X$, $\pi$ the
representation\footnote{by multiplication} of $C(X)$ in
$L^2(X,d\lambda)$. Let $a_j=a_j^*\in C(X)$ and $s$ the map from $X$
to $\R^p$ with coordinates $a_j$. We let $U\subset X$ be an open set
and $\nu$   a measure on $\R^p$ and assume that
\begin{itemize}
  \item The restriction of $s$ to $U$ is an open mapping.
  \item For every open subset $V\subset U$ the image $s(\lambda|_V)$ is equivalent to
  the restriction of $\nu$ to $s(V)$.
  \end{itemize}
Let then $V\subset U$ be an open set and consider the operators
$T_j=\pi(a_j)|_V$ obtained by restriction of the $\pi(a_j)$ to the
subspace $L^2(V,d\lambda)\subset L^2(X,d\lambda)$. Then the joint
spectral measure of the $T_j$ is $\nu|_{s(V)}$ and the spectral
multiplicity $\Sigma(y)$ fulfills
\begin{equation}\label{specmult}
   \Sigma(y)\geq n(y)=\#\{s^{-1}(y)\cap V\} \qqq y\in s(V)\,,
\end{equation}
almost everywhere modulo $\nu$.
\end{lem}

\proof Let $W=s(V)$ which is a bounded open set in $\R^p$. One can
desintegrate $\lambda|_V$ in the form
\begin{equation}\label{desinte}
\lambda|_V=\int_W \rho_y d\nu(y)
\end{equation}
where the $\rho_y$ are positive measures carried by
$F_y=s^{-1}(y)\cap V$. Moreover the total mass of $\rho_y$ is $>0$
almost everywhere modulo $\nu$ for $y\in s(V)$. This follows from
the assumed equivalence $s(\lambda|_V)\sim \nu|_{s(V)}$. One then
has
$$
L^2(V,d\lambda)=\int_W^\oplus L^2(F_y,\rho_y)d\nu(y)
$$
For any $\xi\in L^2(V,d\lambda)$ and any $f\in C_c(\R^p)$ one has
$$
\langle \xi, f((a_j))\xi\rangle=\int_W\int_{F_y} |\xi(x)|^2d\rho_y
f(y)d\nu(y)
$$
which shows that the joint spectral measure of the $a_j$ is
absolutely continuous with respect to $\nu|_W$. Its equivalence with
$\nu|_W$ follows from \eqref{desinte} taking $\xi(x)=1$ and using
the assumed equivalence of $s(\lambda|_V)$ with
  the restriction of $\nu$ to $s(V)$.

\smallskip

Let us prove \eqref{specmult}.   Let $y\in W$ with
$n(y)=\#\{s^{-1}(y)\cap V\}\geq m>0$. Let $x_j\in V$ be $m$ distinct
points in $s^{-1}(y)\cap V$. One can then find disjoint open sets
$B_j\ni x_j$ and let $Z=\cap_j\, s(B_j)$ which is an open
neighborhood of $y$. For any $z\in Z$ the preimage $s^{-1}(z)\cap V$
contains at least $m$ points since it contains at least one in each
$B_j$. Moreover one has $s(s^{-1}(Z)\cap B_j)=Z$ for all $j$. Let
$\lambda_j$ be the restriction of $\lambda$ to $V_j=s^{-1}(Z)\cap
B_j$. From the first part of the Lemma,
 for each $j$ the spectral measure of the
action of the $a_j$ in $L^2(V_j,d\lambda_j)$ is the restriction
$\nu|_{s(V_j)}=\nu_Z$. The action of the $a_j$ in
$L^2(V,d\lambda_j)$ contains the direct sum of the actions in the
$L^2(V_j,d\lambda_j)$ and hence a copy of the action of the
coordinates $y_j$ in
$$
\oplus_1^m L^2(Z,\nu_Z)
$$
which shows that the spectral multiplicity fulfills $\Sigma(z)\geq m$ \aee on
the neighborhood $Z$ of $y$. This shows that any element $y$ in the open set $
U_m=\{y\in U\,|\,n(y)\geq m\}$ admits an open neighborhood $Z_y$ where
$\Sigma(z)\geq m$ holds \aee Since $U_m$ is $\sigma$-compact it follows that
$\Sigma(y)\geq m$ almost everywhere modulo $\nu$ on $U_m$ so that
\eqref{specmult} holds.
\endproof

\begin{rem} \label{remlowermulti} {\rm With the hypothesis of Lemma \ref{lowermulti},
let $E$ be a complex hermitian  vector bundle over $X$ with non-zero fiber
dimension everywhere. Then the inequality $ \Sigma_E(y)\geq n(y)$ holds, where
$\Sigma_E$ is the spectral multiplicity of the $T_j$ acting on
$L^2(X,d\lambda,E)$. This follows since, at the measurable level, one can find
a nowhere vanishing section of $E$ which shows that the representation $\pi_E$
of $C(X)$ in $L^2(X,d\lambda,E)$ contains the representation $\pi$ of $C(X)$ in
$L^2(X,d\lambda)$. Since $\pi_E$ is contained in the sum of $N$ copies of
$\pi$, one obtains the conclusion. }\end{rem}

\begin{thm} \label{thmspecmult} Let $V\subset U_\alpha$ be an open set and let
$a_\alpha^j|_V$ be the restriction of $a_\alpha^j\in \cA$ to the
range  $1_V\cH\subset \cH$. Then
\begin{itemize}
  \item The joint spectral measure of the $a_\alpha^j|_V$ is the Lebesgue measure on $s_\alpha(V)$.
  \item The spectral multiplicity
$m_{ac}(y)$ fulfills
\begin{equation}\label{specmultlebesgue}
   m_{ac}(y)\geq n(y)=\#\{s_\alpha^{-1}(y)\cap V\} \qqq y\in s_\alpha(V)\,,
\end{equation}
almost everywhere modulo the Lebesgue measure.
\end{itemize}
\end{thm}

\proof By Lemma \ref{specmeas} the hypothesis of Lemma
\ref{lowermulti} are fulfilled by the compact space $X$ with measure
$\lambda$, the open set $U_\alpha$, the measure $d\nu=dx^p$ and the
elements $a_\alpha^j$. Thus the result follows from Lemma
\ref{lowermulti} and Remark \ref{remlowermulti}.\endproof

\medskip \section{Local form of the $\cL^{(p,1)}$ estimate}\label{sectvoic}

We fix $p\in [1,\infty[$. Our goal is to control the size of the
Lebesgue multiplicity $m_{ac}(y)$ which appears in Theorem
\ref{thmspecmult}. The idea here is to use a local form of the
$\cL^{(p,1)}$ estimate of \cite{Co-book} Proposition IV.3.14, with
the right hand side of the inequality now involving a closed subset
$K\subset X$, by
$$
 \lambda(K)=\inf_{b\in \cA^+,\,b1_K=1_K}
\cutint b |D|^{-p}
$$
It relies on the estimate given in \cite{Co-action1} and on the
crucial results of Voiculescu (\cite{Voic1}, \cite{Voic2},
\cite{Voic3}). The norm $\Vert T\Vert_{(p,1)}$ is
defined\footnote{For $p=1$ it agrees with the $\cL^1$-norm} for a
compact operator $T$ with characteristic values $\mu_n(T)$ in
decreasing order by (\cf \cite{Voic1} p. 5)
\begin{equation}\label{p1norm}
   \Vert T\Vert_{(p,1)}=\sum_1^\infty n^{-1+1/p}\mu_n
\end{equation}
In order to get an upper bound on $\Vert T\Vert_{(p,1)}$ for $T$ an
operator of finite rank, we can use an inequality of the form
\begin{equation}\label{esti1}
    \Vert T\Vert_{(p,1)}\leq C_p\, (\Rank T)^{1/p}\, \Vert T\Vert_\infty
\end{equation}
which follows using the characteristic values $\mu_n(T)$ from
$$
\Vert T\Vert_{(p,1)}=\sum_{n=1}^{N} n^{-1+1/p}\mu_n\leq \Vert
T\Vert_\infty \sum_{n=1}^{N} n^{-1+1/p}\leq C_p\, N^{1/p}\, \Vert
T\Vert_\infty
$$
where $N=\Rank T$.
 Note also that the $\cL^{(p,1)}$ norm fulfills
\begin{equation}\label{bimine}
\Vert ATB \Vert_{(p,1)}\leq \Vert A \Vert_\infty\Vert T
\Vert_{(p,1)}\Vert B \Vert_\infty\,.
\end{equation}

\smallskip

Let $D$ be a self-adjoint unbounded operator such that its resolvent
is an infinitesimal of order $1/p$ \ie is such that the
characteristic values  fulfill $\mu_n(|D|^{-1})=O(n^{-1/p})$. We let
for any $\lambda>0$,
\begin{equation}\label{spectralprojbis}
    P(\lambda)={\bf 1}_{[0,\lambda]}(|D|)\,, \ \ \alpha(\lambda)=
    \Tr P(\lambda)\,.
\end{equation}
By construction $\alpha(\lambda)$ is a non decreasing integer valued
function. The hypothesis $\mu_n(|D|^{-1})=O(n^{-1/p})$ implies
$\mu_n(|D|^{-1})< C n^{-1/p}$ for some $C<\infty$ and it follows
that $\alpha(C^{-1} n^{1/p})< n$, since the $n$-th eigenvalue of
$|D|$ in increasing order is $>C^{-1} n^{1/p}$. Thus, using for $n$
the smallest integer above $C^p\,\lambda^p$, we get
\begin{equation}\label{upperboundalpha}
    \alpha(\lambda)\leq C^p\,\lambda^p \qqq \lambda >0\,.
\end{equation}

 Let us show that
\begin{lem}\label{lemnonzerodix} Let $f\in C_c^\infty(\R)$,
then there is a finite constant $C_f$ such that
\begin{equation}\label{esti0}
    \Vert [f(\epsilon
D),a]\Vert_\infty \leq C_f\epsilon \Vert [D,a]\Vert \qqq a\in \cA\,.
\end{equation}
Under the hypothesis of Theorem \ref{thmspecmult}, one has:
\begin{equation}\label{nonzerodixbis}
\liminf \lambda^{-p} \alpha(\lambda)>0\,.
\end{equation}
\end{lem}

\proof One has
\begin{equation}\label{dercom}
[e^{is\epsilon D},a]=is\epsilon \int_0^1 e^{ius\epsilon
D}[D,a]e^{i(1-u)s\epsilon D}du
\end{equation}
 which gives \eqref{esti0} using the finiteness of
 $\int |s\,\hat f(s)|ds$ and
 \begin{equation}\label{dercom1}
[f(\epsilon D),a]=(2\pi)^{-1}\int \hat f(s)[e^{is\epsilon D},a]ds\,.
\end{equation}
Assume that \eqref{nonzerodixbis} does not hold. Then let
$\lambda_n\to \infty$ be such that $\lim \lambda_n^{-p}
\alpha(\lambda_n)=0$. Let $f\in C_c^\infty(\R)$ be an (even) cutoff
function vanishing outside $[-1,1]$. For $\epsilon_n=\lambda_n^{-1}$
one has
$$
\Rank f(\epsilon_n D)\leq \alpha(\lambda_n)\,, \ \Rank [f(\epsilon_n
D),a]\leq 2\alpha(\lambda_n)
$$
so that by \eqref{esti1} one gets, using \eqref{esti0},
$$
\Vert [f(\epsilon_n D),a]\Vert_{(p,1)}\leq C_p
(2\alpha(\lambda_n))^{1/p}C\epsilon_n \Vert [D,a]\Vert
$$
and since $\lim \lambda_n^{-p} \alpha(\lambda_n)=0$,
\begin{equation}\label{liminfzero}
\lim_{n\to \infty} \Vert [f(\epsilon_n D),a]\Vert_{(p,1)}=0
\end{equation}
The Voiculescu obstruction relative to an ideal $J$ of compact
operators is given by
\begin{equation}\label{voicobstr}
k_J(\{a_j\})=\liminf_{A\in\cR_1^+,A\uparrow 1}\,\max \Vert
[A,a_j]\Vert_J
\end{equation}
where $\cR_1^+$ is the partially ordered set of positive, finite
rank operators of norm less than one, in $\cH$. We take
$A_n=f(\epsilon_n D)$. It is by construction an element of
$\cR_1^+$. Moreover since $f(\epsilon_n D)\to 1$ strongly in $\cH$
this shows that for the ideal $J=\cL^{(p,1)}$ one gets
$k_J(\{a_j\})=0$. This contradicts the existence, shown in Theorem
\ref{thmspecmult}, of $p$ self-adjoint elements $a_j$ of $\cA$ whose
joint spectral measure is the Lebesgue measure, using Theorem 4.5 of
\cite{Voic1} which gives the equality, valid for $p$ self-adjoint
operators
\begin{equation}\label{voicequfirst}
   k_J(\{a_j\})^p=\gamma_p\int_{\R^p} m(y)d^py
\end{equation}
where the function $m(y)$ is the multiplicity of the Lebesgue
spectrum.
\endproof
 The
rank of the operator $T=[f(\epsilon D),a]$ is controlled by twice
the rank of $f(\epsilon D)$. We take $f$ compactly supported and
thus $f\leq g$ where $g$ is equal to one on the support of $f$
yields an inequality of the form
\begin{equation}\label{ranktrace}
\Rank f(\epsilon D)\leq \Tr( g(\epsilon D))
\end{equation}

One has  by Corollary \ref{corheatequdix} (Appendix 2) an estimate
of the form
\begin{equation}\label{liminfestim}
\liminf \epsilon^p\,\Tr (g(\epsilon D))\leq c_g\cutint |D|^{-p}
\end{equation}
This gives
\begin{equation}\label{esti2}
\liminf \epsilon^p\,\Rank f(\epsilon D)\leq c_g\cutint |D|^{-p}
\end{equation}
and:

\begin{lem} \label{booklem} Let $f\in C_c^\infty(\R)$,
then there is a finite constant $c_f$ such that
\begin{equation}\label{esti3}
   \liminf \Vert [f(\epsilon
D),a]\Vert_{(p,1)}\leq c_f (\cutint |D|^{-p} \,)^{1/p}\Vert
[D,a]\Vert \qqq a\in \cA\,.
\end{equation}
\end{lem}

\proof Using  \eqref{esti2}, and $\Rank [f(\epsilon D),a]\leq 2
\Rank f(\epsilon D)$ one gets a sequence $\epsilon_q\to 0$ with
$$ \Rank T_q \leq 3\,\epsilon_q^{-p}c_g\cutint |D|^{-p}\,,\ \
T_q=[f(\epsilon_q D),a]$$ Using \eqref{esti1} and \eqref{esti0} then
gives
$$
\Vert T_q\Vert_{(p,1)}\leq C_p\, (\Rank T_q)^{1/p}\, \Vert
T_q\Vert_\infty \leq C_p\,(3\,\epsilon_q^{-p}c_g\cutint
|D|^{-p})^{1/p}C_f\epsilon_q \Vert [D,a]\Vert
$$
which is the required estimate since
$(\epsilon_q^{-p})^{1/p}\epsilon_q=1$.\endproof

We now let $K\subset X$ be a compact subset and we want to localize
the estimate \eqref{esti3} to $K$ \ie to the range of $K$ in $\cH$.

\begin{lem} \label{derivlem}\footnote{This is Lemma 10.29 in \cite{FGV}, but the
proof given there is not correct, so we give the full details here.}
 Let $h\in C_c^\infty(\R)$ be an (even) cutoff function
and $f=h^2$. Then
\begin{equation}\label{esti4}
   \Vert [f(\epsilon
|D|),a]-\frac 12 \epsilon (f'(\epsilon |D|)\delta(a)+
\delta(a)f'(\epsilon |D|))\Vert_{(p,1)}=O(\epsilon)
\end{equation}
where $\delta(a)=[|D|,a]$ and one assumes that $a\in
\cap_{j=1}^2\Dom \,\delta^j$.
\end{lem}

\proof First one has (\cf Corollary 10.16 of \cite{FGV}),
\begin{equation}\label{derivlembis}
    \Vert [h(\epsilon|D|),a]-\epsilon h'(\epsilon|D|)\delta(a)\Vert
    \leq C_2 \,\epsilon^2\Vert \delta^2(a)\Vert
\end{equation}
with a similar estimate using $\epsilon\delta(a)h'(\epsilon|D|)$.
Indeed, using \eqref{dercom1}  with $|D|$ instead of $D$, one gets
$$
[h(\epsilon |D|),a]=(2\pi)^{-1}\int \hat h(s)[e^{is\epsilon
|D|},a]ds$$ so that by \eqref{dercom}
$$[h(\epsilon |D|),a]=(2\pi)^{-1}\int \hat
h(s)is\epsilon \int_0^1 e^{ius\epsilon |D|}[|D|,a]e^{i(1-u)s\epsilon
|D|}du ds
$$
and since  by \eqref{dercom} one has
$$
\Vert[[|D|,a],e^{i(1-u)s\epsilon |D|}]\Vert \leq |s\epsilon|\Vert
\delta^2(a)\Vert
$$
one gets
$$
\Vert [h(\epsilon|D|),a]-\epsilon h'(\epsilon|D|)\delta(a)\Vert \leq
C_2\, \epsilon^2\Vert \delta^2(a)\Vert\,, \ \ C_2=(2\pi)^{-1}\int
s^2|\hat h(s)|ds\,.
$$

We follow the proof of Lemma 10.29 in \cite{FGV}. One has
$$
[f(\epsilon |D|),a]-\frac 12 \epsilon (f'(\epsilon |D|)\delta(a)+
\delta(a)f'(\epsilon |D|))=A_\epsilon B_\epsilon+ C_\epsilon
A_\epsilon
$$
where $A_\epsilon =h(\epsilon|D|)$,
$B_\epsilon=[h(\epsilon|D|),a]-\epsilon h'(\epsilon|D|)\delta(a)$
and
$C_\epsilon=[h(\epsilon|D|),a]-\epsilon\delta(a)h'(\epsilon|D|)$. By
\eqref{derivlembis} one has $\Vert B_\epsilon\Vert =O(\epsilon^2)$,
$\Vert C_\epsilon\Vert =O(\epsilon^2)$, while $A_\epsilon$ is
uniformly bounded with   $\Rank A_\epsilon=O(\epsilon^{-p})$. Thus
by \eqref{esti1} one has $\Vert A_\epsilon
\Vert_{(p,1)}=O(\epsilon^{-1})$. Thus we get the required estimate
using \eqref{bimine}.
\endproof

We then let $K\subset X$ be a compact subset, as above, and consider
the operators
\begin{equation}\label{repsilon}
R_\epsilon=1_K\,f(\epsilon |D|)\,1_K
\end{equation}
We let $b\in \cA$ be equal to $1$ on $K$ \ie such that $b\,1_K=1_K$.
One then has:

\begin{lem} \label{derivlem1}
\begin{equation}\label{esti5}
   \Vert [R_\epsilon,a]-\frac 12 \epsilon (1_K\,f'(\epsilon |D|)b\delta(a)1_K\,+
1_K\,\delta(a)b f'(\epsilon |D|)1_K\,)\Vert_{(p,1)}=O(\epsilon)
\end{equation}
\end{lem}

\proof One has
$$
[R_\epsilon,a]=1_K\,[f(\epsilon |D|),a]\,1_K
$$
since $a$ commutes with $1_K$. Thus multiplying on both sides by
$1_K$ in \eqref{esti4}, one gets (using \eqref{bimine})
\begin{equation}\label{esti6}
\Vert [R_\epsilon,a]-\frac 12 \epsilon (1_K\,f'(\epsilon |D|)
\delta(a)1_K\,+ 1_K\,\delta(a) f'(\epsilon
|D|)1_K\,)\Vert_{(p,1)}=O(\epsilon)
\end{equation}
 Lemma \ref{lemnonzerodix} and \eqref{upperboundalpha} show, using
 \eqref{esti1},
  that one has a uniform upper
bound
$$
\Vert [f'(\epsilon |D|),b]\Vert_{(p,1)}\leq C \Vert [|D|,b]\Vert
$$
since $f'$ has compact support. Thus in \eqref{esti6} one can
replace $ 1_K\,f'(\epsilon |D|)=1_K\,b\,f'(\epsilon |D|)$ by $
1_K\,f'(\epsilon |D|)b $, without affecting the behavior in
$O(\epsilon)$. The same applies to the other term.\endproof

We recall the interpolation inequality used in \cite{Co-book} \S
IV.2.$\delta$, but stated without proof there.

\begin{lem} \label{interpol}
There exists for $1\leq p<\infty$, a constant $c_p$ such that, for
$S\in \cL^1$,
\begin{equation}\label{inter}
    \Vert S\Vert_{(p,1)}\leq c_p \Vert S\Vert_1^{1/p} \Vert S\Vert_\infty^{1-1/p}
\end{equation}
\end{lem}

\proof The inequality holds as an equality for $p=1$ with $c_1=1$,
thus we can assume that $p>1$. We use the fact that $\cL^{(p,1)}$ is
obtained by real interpolation of index $(\theta,1)$ for $\theta
=\frac 1p$ from the Banach spaces $Y_0=\cK$ and $Y_1=\cL^1$. The
functoriality of the interpolation gives an inequality of the form
$$
\Vert T(x)\Vert_{(Y_0,Y_1)_{(\theta,q)}}\leq
M_0^{1-\theta}M_1^\theta \Vert x\Vert_{(X_0,X_1)_{(\theta,q)}}
$$
for any linear operator from $X_0+X_1$ to $Y_0+Y_1$ such that
$$
\Vert Tx \Vert_{Y_i}\leq M_i \Vert x \Vert_{X_i}\qqq x\in X_i, i=0,1
$$
We can take $X_0=X_1=\C$ and let $T$ be such that $T(1)=S$. Then
$M_0=\Vert S\Vert_\infty$, $M_1=\Vert S\Vert_1$ and the norm $\Vert
x\Vert_{(X_0,X_1)_{(\theta,1)}}$ is finite and non-zero.\endproof

\begin{rem}\label{doublecheck} {\rm In order not to depend on
interpolation theory we give a direct proof of \eqref{inter}. We
assume that $p>1$. First, for $p>1$ an equivalent norm on
$\cL^{(p,1)}$ is
\begin{equation}\label{lp1norm}
    \Vert T\Vert_{(p,1)'}=(1-\theta)\;\sum N^{\theta -2}\sigma_N(T)  \,, \ \theta =\frac 1p\,.
\end{equation}
where $\sigma_N(T)$ is the sum of the first $N$ characteristic
values. The equivalence of the norms \eqref{lp1norm} and
\eqref{p1norm} follows from $\mu_N\leq \sigma_N/N$ one way. For the
other way, one applies Fubini to  the double sum
$$
\sum_n \sum_{m\geq n} \mu_n m^{\theta -2}=\sum_m \sum_{n\leq m}
\mu_n m^{\theta -2}\,.
$$
Now to estimate \eqref{lp1norm} assuming $\Vert T\Vert_\infty\leq 1$
and $\Vert T\Vert_1=\rho\geq 1$, one splits the sum as follows
$$
\sum_1^\infty N^{\theta -2}\sigma_N(T)=\sum_{N<\rho} N^{\theta
-2}\sigma_N(T)+\sum_{N\geq \rho} N^{\theta -2}\sigma_N(T)\,.
$$
Using $\Vert T\Vert_\infty\leq 1$ gives $\sigma_N(T)\leq N$ and one
bounds the first sum as
$$
\sum_{N<\rho} N^{\theta -2}N\sim C_\theta\,\rho^\theta
$$
Using $\Vert T\Vert_1=\rho\geq 1$ gives $\sigma_N(T)\leq \rho$ and
one bounds the second sum by
$$
\sum_{N\geq \rho} N^{\theta -2}\rho\sim C'_\theta\, \rho^\theta
$$
which gives the required  inequality \eqref{inter}. }\end{rem}

\begin{lem} \label{loclem} There exists a constant $C_f\leq \infty$
such that for  $b=b^*\in \cA$, $b\geq 0$,
\begin{equation}\label{locinequ}
  \liminf \epsilon^p   \Vert bf'(\epsilon |D|)b\Vert_1\leq
  C_f\cutint b^2 |D|^{-p}
\end{equation}
\end{lem}

\proof Note that by construction of $f$ as a cutoff function, its
derivative $f'\leq 0$ on $[0,\infty[$. Let $h=-f'\in C_c^\infty(\R)$
so that $h\geq 0$. One then has $bh(\epsilon |D|)b\geq 0$ and
$$
\Vert bh(\epsilon |D|)b\Vert_1=\Tr(bh(\epsilon |D|)b)
=\Tr(b^2h(\epsilon |D|))
$$
and the result follows from Corollary \ref{corheatequdix} (Appendix
2) which gives, $$\liminf \epsilon^p \Tr(b^2h(\epsilon |D|))\leq \mu
\cutint b^2 |D|^{-p}
$$
where $\mu= p\int_0^\infty u^{p-1}h(u)du$.
\endproof

\begin{lem} \label{finelem}There exists a constant $C'_f\leq \infty$
such that, for $b=b^*\in \cA$, $0\leq b\leq 1$,
\begin{equation}\label{locinequ1}
  \liminf   \Vert \epsilon  bf'(\epsilon |D|)b\Vert_{(p,1)}\leq
  C'_f(\cutint b^2 |D|^{-p})^{1/p}
\end{equation}
\end{lem}

\proof By Lemma \ref{loclem} one has, once $b$ is fixed, a sequence
$\epsilon_q\to 0$ such that
$$
\Vert bf'(\epsilon_q |D|)b\Vert_1\leq 2C_f \epsilon_q^{-p}\cutint
b^2 |D|^{-p}
$$
Also since $f'$ is bounded one has
$$
\Vert bf'(\epsilon_q |D|)b\Vert_\infty\leq B=\Vert
f'\Vert_\infty<\infty
$$
Thus it follows from \eqref{inter}, that
$$
\Vert bf'(\epsilon_q |D|)b\Vert_{(p,1)}\leq c_p (2C_f
\epsilon_q^{-p}\cutint b^2 |D|^{-p})^{1/p}\,B^{1-1/p}
$$
After multiplication by $\epsilon_q$ one gets the required
estimate.\endproof

\begin{thm} \label{main} There exits a finite constant $\kappa_p$ such that
for any operators $a_j\in \cA$ and compact subset $K\subset X$ one
has, with $J=\cL^{(p,1)}$, the inequality
\begin{equation}\label{voicprep}
   k_J(\{a_j\, 1_K\})\leq \kappa_p \max \Vert
\delta(a_j)\Vert_\infty(\lambda(K))^{1/p}
\end{equation}
where one lets\footnote{This is the natural extension of $\lambda$
given by the Riesz representation Theorem \cite{Rudin1}.}:
$$
\lambda(K)=\inf_{b\in \cA^+,\,b1_K=1_K} \cutint b |D|^{-p}\,.
$$
\end{thm}

\proof By definition one has
$$
k_J(\{a_j\, 1_K\})=\liminf_{A\in\cR_1^+,A\uparrow 1}\,\max \Vert
[A,a_j\, 1_K]\Vert_J
$$
where $\cR_1^+$ is the partially ordered set of positive, finite
rank operators of norm less than one, in $1_K\,\cH$. We take
$R_\epsilon=1_K\,f(\epsilon |D|)\,1_K$ as in \eqref{repsilon}. It is
by construction an element of $\cR_1^+$. Moreover since $f(\epsilon
|D|)\to 1$ strongly in $\cH$ one gets that $R_\epsilon \to 1$
strongly in $1_K\,\cH$. By Lemma \ref{derivlem1} one has
$$
 \Vert [R_\epsilon,a]-\frac 12 \epsilon (1_K\,f'(\epsilon |D|)b\delta(a)1_K\,+
1_K\,\delta(a)b f'(\epsilon |D|)1_K\,)\Vert_{(p,1)}=O(\epsilon)
$$
Using
$$
1_K\,f'(\epsilon |D|)b\delta(a)1_K=1_K\,b\,f'(\epsilon
|D|)b\delta(a)1_K
$$
and \eqref{bimine} for $A=1_K$, $T=b\,f'(\epsilon |D|)b$,
$B=\delta(a)1_K$ and similarly for the other term, one gets an
estimate of the form
$$
 \Vert [R_\epsilon,a]\Vert_{(p,1)}\leq O(\epsilon)+\Vert \epsilon  bf'(\epsilon |D|)b\Vert_{(p,1)}
 \Vert
\delta(a)\Vert_\infty
$$
Thus, by Lemma \ref{finelem} one gets  that, for any $b\in \cA^+$
equal to $1$ on $K$, there exists a sequence $\epsilon_q\to 0$ such
that
$$
\Vert \epsilon_q bf'(\epsilon_q |D|)b\Vert_{(p,1)}\leq
  2C'_f(\cutint b^2 |D|^{-p})^{1/p}
$$
which gives, for $q$ large enough,
$$
\Vert [R_{\epsilon_q},a]\Vert_{(p,1)}\leq
  2C'_f(\cutint b^2 |D|^{-p})^{1/p}\Vert
\delta(a)\Vert_\infty
$$
for any $a\in \cA$ and hence:
$$
\liminf \max \Vert [R_\epsilon,a_j\, 1_K]\Vert_J\leq 2 C'_f \max
\Vert \delta(a_j)\Vert_\infty (\cutint b^2 |D|^{-p})^{1/p}
$$
After varying $b$ one obtains the required estimate.
\endproof

\begin{rem} {\rm a) One may worry that Voiculescu's definition of $k_J$
involves the ordered set $\cR_1^+$ while all we got was $R_\epsilon
\to 1$ strongly in $1_K\,\cH$. Thus let us briefly mention how to
get the $A\uparrow 1$ from $R_\epsilon$ by a small modification.
Given a finite dimensional subspace of $1_K\,\cH$, one lets $P$ be
the corresponding finite rank projection, with fixed rank $N$. One
needs to construct $A\in\cR_1^+$, $A\geq P$ with a control on $\max
\Vert [A,a_j\, 1_K]\Vert_J$. One takes
$$
A_\epsilon=P+ (1-P)R_\epsilon (1-P)\,,
$$
so that $A\geq P$ by construction. Moreover one has
\begin{equation}\label{app}
   R_\epsilon-A_\epsilon=P(R_\epsilon-1) P+ P R_\epsilon(1-P)+
(1-P)R_\epsilon P\,.
\end{equation}
 Moreover by the strong convergence $R_\epsilon \to 1$, one has $$
\Vert P(R_\epsilon-1)\Vert_\infty=\Vert (R_\epsilon-1)P\Vert_\infty
\to 0 $$ so that all three terms in the rhs of \eqref{app} converge
to $0$ in norm and hence in the $J$ norm since their rank is less
than $N$ so that one can use \eqref{esti1}. Thus one has $$ \Vert
R_\epsilon-A_\epsilon \Vert_J\to 0
$$ and one controls $\max
\Vert [A,a_j\, 1_K]\Vert_J$ from  $\max \Vert [R_\epsilon,a_j\,
1_K]\Vert_J$.

b) It might seem possible at first sight to tensor the spectral
triple $(\cA,\cH,D)$ by $(\C,\cH',D')$, with the spectrum of $D'$
growing fast enough so that the product triple\footnote{in the even
case}
\begin{equation}\label{proucttriple}
    (\cA\otimes \C,\cH\otimes \cH',D''= D\otimes 1+\gamma\otimes D')
\end{equation}
would still be of dimension $p$, \ie such that the characteristic
values of the inverse of $D''$ are $O(n^{-1/p})$. Let us show that
this is only possible if the dimension of $\cH'$ is {\em finite}.
Indeed the eigenvalues of $(D\otimes 1+\gamma\otimes
D')^2=D^2\otimes 1+1\otimes D^{'2}$ are the independent sums of the
eigenvalues of $D^2$ and of $D^{'2}$. Thus having infinitely many
eigenvalues of $D^{'2}$ contradicts the two inequalities
$$
\alpha(\lambda)\geq c\lambda^p \,, \ \ \alpha''(\lambda)\leq C''
\lambda^p
$$
for the counting functions
$\alpha(\lambda)=\Tr(1_{[0,\lambda]}(|D|)$,
$\alpha''(\lambda)=\Tr(1_{[0,\lambda]}(|D''|)$ since they yield
$$
\dim(\cH')\leq C''/c\,.
$$

c) The constant $C'_f$ in \eqref{locinequ1} is given, up to a
function of $p$ alone, by
$$
C'_f=\left(\int_0^\infty u^{p-1}h(u)du\right)^{1/p}\Vert
h\Vert_\infty^{1-1/p}\,, \ \ h=-f'\geq 0\,.
$$
and one needs to check that there is a lower bound to $C'_f$
independently of the choice of the cutoff function $f$. Since
$f(0)=1$ the only information is about $\int_0^\infty h(u)du=1$ and
thus one needs to show a general inequality of the form:
\begin{equation}\label{checkineq}
    \int_0^\infty h(u)du\leq c(p)\left(p\int_0^\infty u^{p-1}h(u)du\right)^{1/p}\Vert
h\Vert_\infty^{1-1/p}\,.
\end{equation}
To prove this one lets $g(u)=h(u^{1/p})$, so that
$$
p\int_0^\infty u^{p-1}h(u)du=\int_0^\infty g(v)dv\,, \ \
\int_0^\infty h(u)du=\frac 1p\,\int_0^\infty v^{1/p-1}g(v)dv
$$
and one uses the same argument as in Remark \ref{doublecheck}.
First, with $k(u)=\int_0^ug(v)dv$,
$$
\int_0^\infty v^{1/p-1}g(v)dv=(1-\frac 1p)\int_0^\infty
v^{1/p-2}k(v)dv
$$
Next, assuming $\Vert h\Vert_\infty=1$, one has $g(v)\leq 1$ for all
$v>0$ and thus, with $\rho=\int_0^\infty g(v)dv$,
$$
\int_0^\infty v^{1/p-2}k(v)dv\leq \int_0^\rho
v^{1/p-2}k(v)dv+\int_\rho^\infty v^{1/p-2}k(v)dv
$$
so that, since $k(v)\leq v$ and $k(v)\leq \rho$ one gets
$$
\int_0^\infty v^{1/p-2}k(v)dv\leq \int_0^\rho
v^{1/p-1}dv+\int_\rho^\infty v^{1/p-2}\rho dv=c_p \rho^{1/p}
$$
with $c_p=p+(1-\frac 1p)^{-1}$ which gives \eqref{checkineq} with
$c(p)=1$.}
\end{rem}

\medskip

We can now combine this with Theorem 4.5 of \cite{Voic1} which gives
the equality, valid for $p$ self-adjoint operators $h_j$,
\begin{equation}\label{voicequ}
   k_J(\{h_j\})^p=\gamma_p\int_{\R^p} m(y)d^py
\end{equation}

\begin{cor} \label{inequbase} Let $a_j=a_j^*\in \cA$ be $p$ self-adjoint elements,
then for any compact subset $K\subset X$ one has
\begin{equation}\label{inequbase1}
    \int_{\R^p} m_{\rm ac}^K(y)d^py \leq \kappa'_p \max \Vert
\delta(a_j)\Vert_\infty^p \lambda(K)
\end{equation}
where the constant $\kappa'_p$ only depends on $p$, and the function
$m_{\rm ac}^K(y)$ is the multiplicity of the Lebesgue spectrum of
the restriction of the $a_j$ to $1_K(\cH)$.
\end{cor}

\medskip \section{Local bound on $\# (s_\alpha^{-1}(x)\cap
V)$}\label{sectlocalbound}

Let $V\subset U_\alpha$ be an open set with $\overline V\subset
U_\alpha$.

\begin{lem} \label{bound0}
There exists $C<\infty$ such that the spectral multiplicity
$m^V_{\rm ac}(x)$, on the absolutely continuous joint spectrum of
the restriction $a_\alpha^j|_V$ of the $a_\alpha^j$, to $1_V\cH$
fulfills:
\begin{equation}\label{finite1}
m^V_{\rm ac}(x)\leq C\,, \ \ \aee\  {\rm on}\; W=s_\alpha(V)
\end{equation}
\end{lem}

\proof By Theorem \ref{thmspecmult}, the joint spectral measure of
the $a_\alpha^j|_V$ is the Lebesgue measure on $s_\alpha(V)$.
 Let $E\subset W$ be a compact subset, and
$K=s_\alpha^{-1}(E)\cap \overline V$. Then Corollary \ref{inequbase}
gives an inequality of the form
\begin{equation}\label{voic}
  \int_{\R^p} m_{\rm ac}^K(y)d^py \leq \kappa' \lambda(K)
\end{equation}
One has
\begin{equation}\label{restri}
    m_{\rm ac}^V(y)\leq m_{\rm ac}^K(y)\qqq y\in E
\end{equation}
 since one has a direct sum decomposition
 $$
1_V\cH=1_{s_\alpha^{-1}(E)\cap V}\cH\oplus 1_{s_\alpha^{-1}(E^c)\cap
V}\cH
 $$
 where the representation in the second term in the right hand-side does not contribute
 to the multiplicity in $E$. Indeed, with $E^c=\cup E_n$ and $E_n$ compact disjoint from
 $E$, the joint spectrum of $a_\alpha^j|_{s_\alpha^{-1}(E_n)\cap
V}$ is contained in $E_n$ and disjoint from $E$. Moreover the
representation in the first term is dominated by the
 representation in $1_K\cH$ since $s_\alpha^{-1}(E)\cap V\subset K=s_\alpha^{-1}(E)\cap \overline
 V$.

 By
\eqref{specmeasmaj} one has an inequality
$$
 \int_{\overline V} f\circ s_\alpha\,d\lambda_{\overline V}\leq c
 \, \int_{s_\alpha({\overline V})} f(x) dx^p
 \leq c
 \, \int_{\R^p} f(x) dx^p\qqq f\in
 C_c^+(\R^p)\,,
$$
which shows, taking $1_E=\inf f_n$ as an infimum of continuous
functions $f_n\in C_c^+(\R^p)$, that
$$
\lambda(K)=\int_{\overline V} 1_E\circ s_\alpha\,d\lambda_{\overline
V}\leq \int_{\overline V} f_n\circ s_\alpha\,d\lambda_{\overline
V}\leq c
 \, \int_{\R^p} f_n(x) dx^p\to c
 \, \int_{E}  dx^p\,.
$$
Thus, using \eqref{restri} and \eqref{voic},
$$
\int_E m^V_{\rm ac}(x)d^px\leq\int_E m^K_{\rm ac}(x)d^px\leq \kappa'
\lambda(K)\leq c \kappa'
 \, \int_{E}  dx^p\,.
$$
 and there exists a constant $C=c \kappa'$ such that, for any compact
$E\subset W$,
\begin{equation}\label{finite}
   \int_E m^V_{\rm ac}(x)d^px\leq C\int_E d^px
\end{equation}
which gives the inequality, valid almost everywhere,
$$
m^V_{\rm ac}(x)\leq C\,.
$$
\endproof

\begin{lem} \label{bound}
Let $V$  be as above. Then there exists $m<\infty$ such that
\begin{equation}\label{finite2}
\#(s_\alpha^{-1}(x)\cap V)\leq m \qqq x\in W=s_\alpha(V)
\end{equation}
\end{lem}

\proof By Theorem \ref{thmspecmult}, one has, almost everywhere,
$$
 m^V_{ac}(y)\geq n(y)=\#\{s_\alpha^{-1}(y)\cap V\} \qqq y\in s_\alpha(V)\,,
$$
so that the result follows from Lemma \ref{bound0} and the
semi-continuity of $n(y)$ which shows that an almost everywhere
inequality remains valid everywhere.
\endproof

\begin{lem} \label{regular} Let $V\subset U_\alpha$ be an open set with
$\overline V\subset U_\alpha$. There exists a dense open subset
$Y\subset s_\alpha(V)$ such that every point of
$s_\alpha^{-1}(Y)\cap V$ has a neighborhood $N$ in $X$ such that the
restriction of $s_\alpha$ to $N$ is an homeomorphism with an open
set of $\R^p$.
\end{lem}

\proof Let $W= s_\alpha(V)$ and
$$
m_1=\sup_{x\in W} \#(s_\alpha^{-1}(x)\cap V)
$$
which is finite (and non-zero) by Lemma \ref{bound}. Let
$$
W_1=\{x\in W\,|\,\#(s_\alpha^{-1}(x)\cap V)=m_1\}
$$
It is by Lemma \ref{semicont} an open subset of $W$. Moreover for
$x\in W_1$ one can find $m_1$ disjoint open neighborhoods $V_j$ of
the preimages $x_j$ of $x$ such that all $V_j$ surject on the same
neighborhood $U$ of $x$ in $W$. It follows that the restriction of
$s_\alpha$ to each of the $V_j$ is a bijection onto $U$ and hence an
isomorphism of a neighborhood of $x_j$ with an open set in $\R^p$
given by the $a_\alpha^j$.

It can be that $W_1$ is not dense in $W$, but then we just take the
complement of its closure: $W^1=W\backslash \overline W_1$ and let
$$
m_2=\sup_{x\in W^1} \#(s_\alpha^{-1}(x)\cap V)
$$
which is $<m_1$ by construction. One then defines
$$
W_2=\{x\in W^1\,|\,\#(s_\alpha^{-1}(x)\cap V)=m_2\}
$$
which is by Lemma \ref{semicont} an open subset of $W^1$. The same
argument as above shows that the subset $Z=W_1\cup W_2$ fulfills the
condition of the Lemma.  One proceeds in the same way and gets, by
induction, a sequence $W_k$, with $Y=\cup  W_j$ fulfilling the
condition of the Lemma. Since the sequence $m_j$ is strictly
decreasing one gets that the process stops and $Y$ is dense in $W$.
\endproof

\medskip \section{Reconstruction Theorem}\label{sectloccoord}

We shall now use Lemma \ref{regular} together with the ability to
move around in $X$ by automorphisms of $\cA$ to prove the following
key Lemma:

\begin{lem} \label{reglem} For every point $\chi\in X$ there exists
$p$ real elements $x^\mu\in \cA$ and a smooth family $\tau_t\in
\Aut(\cA)$, $t\in \R^p$, $\tau_0=\id$, such that
\begin{itemize}
  \item The $x^\mu$ give an homeomorphism
of a neighborhood of $\chi$ with an open set in $\R^p$.
  \item The map $t\mapsto h(t)=\chi\circ \tau_t$ is an homeomorphism
of  a neighborhood of $0$ in $\R^p$  with a neighborhood of $\chi$.
\item The map $x\circ h$ is a local diffeomorphism.
\end{itemize}
\end{lem}

\proof Let $\chi\in X$. By Lemma \ref{coverlem},    there exists
$\alpha$ such that $\chi\in U_\alpha$.  By Lemma \ref{coropen},
there exists a smooth family $\sigma_t\in \Aut(\cA)$, $t\in \R^p$, a
neighborhood $Z$ of $\chi$ in $X=\Sp(\cA)$ and a neighborhood $W$ of
$0\in \R^p$ such that, for any $\kappa\in Z$, the map $t\mapsto
s_\alpha(\kappa\circ\sigma_t)$ is a diffeomorphism, depending
continuously on $\kappa$, of $W$ with a neighborhood of
$s_\alpha(\kappa)$ in $\R^p$.

 We now take for $V$ a ball
$$
V=B_r=\{y\in X\,| \, d(\chi,y)<r\}\subset U_\alpha\,, \ \ \overline
V \subset U_\alpha\,.
$$
 We apply Lemma \ref{regular} to $V=B_r$ and let $Y$ be a
dense open subset $Y\subset s_\alpha(V)$ such that every point of
$s_\alpha^{-1}(Y)\cap V$ has a neighborhood $N$ in $X$ such that the
restriction of $s_\alpha$ to $N$ is an homeomorphism with an open
set of $\R^p$. Since $Y$ is dense in $s_\alpha(V)$, and by Lemma
\ref{coropen} the image of $W$ by $t\mapsto
\psi(t)=s_\alpha(\chi\circ \sigma_t)$ is an open neighborhood of
$s_\alpha(\chi)$ one can choose a $u_0\in W$ such that $\chi\circ
\sigma_{u_0}\in V$ and $\psi(u_0)=s_\alpha(\chi\circ
\sigma_{u_0})\in Y$. One has $\kappa=\chi\circ \sigma_{u_0}\in
s_\alpha^{-1}(Y)\cap V$. Thus by Lemma \ref{regular} there exists a
neighborhood $N$ of $\kappa$ such that the restriction of $s_\alpha$
to $N$ is an isomorphism with an open set of $\R^p$. Thus the
$a_\alpha^\mu$ are good local coordinates near $\kappa$. The
automorphism $\sigma_{u_0}\in \Aut(\cA)$ is such that
\begin{equation}\label{move1}
    \kappa=\chi\circ \sigma_{u_0}\,,\ \  \chi=\sigma_{u_0}(\kappa)
\end{equation}
Recall that we use the covariant notation \eqref{covariant}.
 We take
\begin{equation}\label{move2}
x^\mu=\sigma_{u_0}(a_\alpha^\mu)
\end{equation}
 as local coordinates near $\chi$. The corresponding map $x$ from $X=\Sp(\cA)$ to
 $\R^p$ is given by
 $$
\zeta\in X\mapsto
\zeta(x^\mu)=\zeta(\sigma_{u_0}(a_\alpha^\mu))=s_\alpha(\zeta\circ
\sigma_{u_0})=s_\alpha\circ \sigma_{u_0}^{-1}(\zeta)
 $$
Thus $x=s_\alpha\circ \sigma_{u_0}^{-1}$ and, since $\sigma_{u_0}$
is an homeomorphism of $X$, $x=s_\alpha\circ \sigma_{u_0}^{-1}$ is
an homeomorphism of the neighborhood $\sigma_{u_0}(N)$ of $\chi$
with an open set of $\R^p$. Thus the $x^\mu$ are good local
coordinates at $\chi$. Then let $\tau_t\in \Aut(\cA)$ be given by
$$
\tau_t  = \sigma_{u_0+t}\circ\sigma_{u_0}^{-1}$$ so that $\tau_t
\circ\sigma_{u_0} = \sigma_{u_0+t} $. One has
$$
\chi\circ \tau_t(x^\mu)=\chi\circ
\tau_t(\sigma_{u_0}(a_\alpha^\mu))=\chi\circ
\sigma_{u_0+t}(a_\alpha^\mu)=s^\mu_\alpha(\chi\circ
\sigma_{u_0+t})=\psi^\mu(u_0+t)
$$
Now the map $h$ is given by $t\mapsto h(t)=\chi\circ \tau_t$ thus
one has
\begin{equation}\label{hmap}
    x\circ h(t)=\psi(u_0+t)
\end{equation}
This shows that the map $x\circ h$ is a  diffeomorphism from
$W_1=W-u_0$ (which is a neighborhood of  $t=0\in \R^p$  since
$u_0\in W$)  with an open set of $\R^p$.  On $W_1$, the map $h$ is
injective since $x\circ h$ is injective. Thus $h$ is an
homeomorphism with its range. One has $h(0)=\chi$, $h$ is continuous
thus $W_2=h^{-1}(\sigma_{u_0}(N))\cap W_1$ is an open set containing
$0$ and $W'_2=x\circ h(W_2)$ is an open set in $\R^p$. The map $x$
is an homeomorphism of $\sigma_{u_0}(N)$ with an open set in $\R^p$
and $x\circ h$ is an homeomorphism of $W_1$ with an open set in
$\R^p$. One has $h(W_2)\subset \sigma_{u_0}(N)$. Thus
$h(W_2)=x^{-1}(W'_2)\cap \sigma_{u_0}(N)$ is open in
$\sigma_{u_0}(N)$ and since it contains $h(0)=\chi$ we get that $h$
is an homeomorphism of  a neighborhood of $0$ in $\R^p$ with a
neighborhood of $\chi$. Moreover, as we have seen above,  the map
$x\circ h$ is a  diffeomorphism.
\endproof

\begin{lem}\label{localalg} The algebra $\cA$ is locally the algebra of restrictions
of smooth functions on $\R^p$ to a bounded open set of $\R^p$.
\end{lem}

\proof Let $\chi\in X$. By Lemma \ref{reglem}, we can assume that
some $x^\mu\in \cA$ give an homeomorphism of a neighborhood $U$ of
$\chi$ with a bounded open set $x(U)\subset \R^p$. By the smooth
functional calculus the algebra $C_c^\infty(x(U))$ is contained in
$\cA$ using the morphism $ f\in C_c^\infty(x(U))\mapsto f(x^\mu)\in
\cA\,. $ Moreover for any $\kappa\in U$ one has
$\kappa(f(x^\mu))=f(\kappa(x^\mu)$ so that the function $f\circ x$
coincides on $U$ with the element $f(x^\mu)\in \cA$.
 Taking a smaller neighborhood $V$ of $\chi$
  with compact closure in $U$ one
gets that the algebra $C^\infty(\R^p)|_{x(V)}$ of restrictions to
$x(V)$ of smooth functions on $\R^p$ is contained in the algebra of
restrictions to $V$ of elements of $\cA$, using $x$ to identify $V$
with the open set $x(V)\subset\R^p$. We need to show that any
element of $\cA$ restricts to a smooth function on $V$, using the
local coordinates $x$ to define smoothness. For this we use (Lemma
\ref{reglem}) the existence of a smooth family $ \tau\,:\,\R^p\to
\Aut(\cA) $ such that $x\circ \tau$ is a local diffeomorphism around
$\chi$. Thus given $b\in \cA$, to show that the restriction of $b$
to $V$ is smooth, it is enough to show that $\tau_t(b)$ evaluated at
$\chi$ is a smooth function of $t$. This follows from the smoothness
of the family $\tau_t$.
\endproof

\begin{thm} \label{regthm} Let $(\cA,\cH,D)$ be a strongly regular spectral triple
fulfilling the   five \axioms of   \S \ref{prelem} (\cf \cite{CoSM})
with $c$ antisymmetric  Then there exists an oriented   smooth
compact manifold $X$ such that $\cA=C^\infty(X)$.
\end{thm}

\proof We let $X=\Sp (\cA)$ be the spectrum of $\cA$ or equivalently
of the norm closure $A$. By construction it is a compact space. By
Lemma \ref{reglem}, for every point $\chi\in X$ there exists a
neighborhood $U$ of $\chi$ and $p$ real elements $x^\mu\in \cA$
which give a local homeomorphism $\phi$ of a neighborhood $V$ of $x$
with an open set in $\R^p$. Moreover by Lemma \ref{localalg} one has
$$
f\in \cA|_V \Leftrightarrow f\circ \phi^{-1}\in
C^\infty(\R^p)|_{\phi(V)}
$$
This shows that on the intersection of such domains of local charts,
the change of chart is of class $C^\infty$. We can thus, using
compactness, take a finite cover and this endows $X$ with a
structure of $p$-dimensional smooth manifold.  Lemma \ref{localalg}
shows that any $a\in \cA$ restricts to a smooth function in each
local chart and thus $\cA\subset C^\infty(X)$. Moreover given $f\in
C^\infty(X)$ there exists for each $V_j$ in the finite open cover of
$X$ an $a_j\in \cA$ with $f|_{V_j}=a_j|_{V_j}$. Then the existence
of partitions of unity (Lemma 2.10 of \cite{ReVa}):
$$
\psi_j\in \cA\,, \ \ \sum \psi_j=1 \,, \  \  {\rm Support} \,
\psi_j\subset V_j
$$
shows that $f$ agrees with $\sum \psi_ja_j\in \cA$. We have shown
that there exists a  smooth compact manifold $X$ such that
$\cA=C^\infty(X)$. The cycle $c$ gives a nowhere vanishing section
of the real exterior power $\wedge^p(TX)$ and hence shows that the
manifold $X$ is oriented.
\endproof

We thus obtain the following characterization of the algebras
$C^\infty(X)$:

\begin{thm} \label{regthm1} An involutive algebra $\cA$ is the algebra
of smooth functions on an oriented   smooth compact manifold if and
only if
 it admits a faithful \footnote{\ie with trivial kernel} representation in a pair $(\cH,D)$ fulfilling
  the five \axioms of
\S \ref{prelem} (\cf \cite{CoSM}) with the cycle antisymmetric and
the strong regularity.
\end{thm}

\proof The direct implication follows from Theorem \ref{regthm}. Conversely,
given an oriented   smooth compact manifold $X$ of dimension $p$, one can take
the representation in $\cH=L^2(X,\wedge^*_\C)$ the Hilbert space of square
integrable differential forms with complex coefficients, and use the choice of
a Riemannian metric to get the signature operator $D=d+d^*$ with the
$\Z/2$-grading $\gamma$ in the even case coming from the Clifford
multiplication by the volume form as in \cite{lawmich} Chapter 5. In the odd
case one uses the Clifford multiplication $\gamma$ by the volume form to reduce
the Hilbert space $\cH$ to the subspace given by $\gamma\xi=\xi$. More
specifically we consider the faithful representation of the Clifford algebra
${\rm Cliff}\,T^*_x(X)$ in $\wedge^*T^*_x(X)$ given by the symbol of $D$, \ie
\begin{equation}\label{cliffreg}
    v . \,\xi=v\wedge \xi-i_v \,\xi \qqq v\in T^*_x(X)\,, \ \xi\in \wedge^*T^*_x(X)
\end{equation}
where $i_v $ is the contraction by $v$. This gives (\cf \cite{lawmich}
Proposition 3.9) a canonical isomorphism of vector spaces ${\rm
Cliff}\,T^*_x(X)\sim \wedge^*T^*_x(X)$.
 We
let $\omega$ be the section of $\wedge^p T^*X$ given at each point by
$\omega=e_1\wedge \cdots \wedge e_p$ where $e_1,\ldots e_p$ is any positively
oriented orthonormal basis. In the Clifford algebra ${\rm Cliff}\,T^*_x(X)$ one
has $\omega^2=(-1)^{\frac{p(p+1)}{2}}$  (\cf \cite{lawmich} (5.26)) and one
defines
\begin{equation}\label{defngamma}
    \gamma\, \xi= i^{-\frac{p(p+1)}{2}}\, \omega\, \xi \qqq \xi\in {\rm
Cliff}\,T^*_x(X)\otimes \C
\end{equation}
where the product $\omega\, \xi$ is the left Clifford multiplication by
$\omega$. By \cite{lawmich} (5.35) this left multiplication is related to the
Hodge star operation by
\begin{equation}\label{hodge}
    \omega\, \xi=(-1)^{k(p-k)+\frac{k(k+1)}{2}}\star \xi   \qqq
    \xi\in\wedge^k\,.
\end{equation}
With these notations one has $d^*=(-1)^{(p+1)}\gamma d\gamma$ (\cf
\cite{lawmich} (5.10)) which shows that $D$ commutes with $\gamma$ when $p$ is
odd and anticommutes with $\gamma$ when $p$ is even.

 To check the orientability \axiom 4), one uses (in both cases of
the Dirac operator or the signature operator) local coordinates $x^\mu$ and the
equalities
\begin{equation}\label{cliffgrad}
    [D,f]= \sum \gamma^\mu\partial_\mu f\,, \ \
    \{\gamma^\mu,\gamma^\nu\}=-2g^{\mu\nu}\,,
\end{equation}
where the $\gamma^\mu$ correspond to the action of $dx^\mu$ through the
representation of the Clifford algebra (given by \eqref{cliffreg} for the
signature operator). One then has, for the multiple commutator,
\begin{equation}\label{multigamma}
    [\gamma^1,\gamma^2,\ldots,\gamma^p]=p!\, i^{\frac{p(p+1)}{2}}(\sqrt g)^{-1}\gamma
\end{equation}
where $\sqrt g$ is the square root of the determinant of the matrix
$g_{\mu\nu}$ and $\gamma=\gamma^*$, $\gamma^2=1$ is the grading in the even
case and  is just $1$ in the odd case\footnote{Since we reduced the Hilbert
space $\cH$ to the subspace given by $\gamma\xi=\xi$}. Thus in these local
coordinates $x^\mu$ the cycle associated to the volume form:
$$
c=\frac{1}{p!}\sum_\sigma \epsilon(\sigma) \sqrt g\otimes
x^{\sigma(1)}\otimes \cdots\otimes x^{\sigma(p)}
$$
fulfills, locally, \axiom 4), up to the power of $i$,
$i^{\frac{p(p+1)}{2}}$. Using a partition of unity gives the global
form of $c$ which is just the Hochschild cycle representing the
global volume form.

 The \axiom of strong regularity is checked using \S
\ref{regularity}. One applies Lemma \ref{delta1delta} to obtain the strong
regularity since we take for $D$ an elliptic differential operator of order one
on a smooth compact manifold and the principal symbol of $D^2$ is a scalar
multiple of the identity. This ensures that for any differential operator $P$
of order $m$ the symbol of order $m+2$ of $[D^2,P]$ vanishes as it is given by
the commutator of the principal symbols of order $m$ and $2$. Thus one gets
that for any differential operator $T$ of order $0$, the operators
$\delta_1^m(T)$ are of the form $P(1+D^2)^{-m/2}$ where $P$ is a differential
operator of order $m$. Thus the theory of elliptic operators (\cf \cite{Gilkey}
Lemma 1.3.4 and 1.3.5) shows that they are bounded. This applies for $D$ the
Dirac operator or the signature operator, thus one gets the strong regularity
in this case.
\endproof

\begin{thm} \label{fundthmspinc} Let $(\cA,\cH,D)$ be a  spectral triple with $\cA$
commutative, fulfilling the five \axioms of \S \ref{prelem} with the
cycle $c$ antisymmetric. Assume that the multiplicity\footnote{We
restrict to the even case} of the action of $\cA''$ in $\cH$ is
$2^{p/2}$. Then there exists a smooth oriented compact (spin$^c$)
manifold $X$ such that $\cA=C^\infty(X)$.
\end{thm}

\proof We need to show that we can dispense with the hypothesis of
strong regularity in Theorem \ref{regthm}. Indeed by the first part
of Remark \ref{modda}, we get
\begin{equation}\label{cliffrelcenter}
    [([D,a][D,b]+[D,b][D,a]),[D,c]]=0\qqq a, b, c \in \cA
\end{equation}
since this is implied by the commutation \eqref{commsymbd} of
$|[D,h]|$ with $[D,c]$. Thus if we work at a point $\chi\in
\Sp(\cA)$ and let $S_\chi$ be the fiber at $\chi\in X$ of the finite
projective module $\cH_\infty$ and $M_\chi\subset \End S_\chi$ be
the subalgebra generated by the $[D,a]$ for $a\in \cA$, it follows
from \eqref{cliffrelcenter} that
\begin{equation}\label{cliffrelcenterbis}
     [D,a][D,b]+[D,b][D,a]\in Z(M_\chi)\qqq a,  b\in  \cA
\end{equation}
where $Z(M_\chi)$ is the center of $M_\chi$. Let $e$ be a minimal
projection in the center of $M_\chi$. The equality $\pi_D(c)=\gamma$
shows that, at the point $\chi$,
$$
\gamma e= \sum_\alpha e
a_\alpha^0[[D,a_\alpha^1],[D,a_\alpha^2],\ldots,[D,a_\alpha^p]]\neq
0
$$
so that the dimension of the space $T^*_e(\chi)=\{e[D,a]\,|\, a\in
\cA\}_\chi$ is at least equal to $p$. In fact, more precisely, for
some $\alpha$, the multiple commutator
$$
[e[D,a_\alpha^1],e[D,a_\alpha^2],\ldots,e[D,a_\alpha^p]]=e
[[D,a_\alpha^1],[D,a_\alpha^2],\ldots,[D,a_\alpha^p]]\neq
0
$$
which can hold only if the $e[D,a_\alpha^j]$ are linearly
independent. By \eqref{cliffrelcenterbis} and the minimality of $e$,
the following equality defines a positive quadratic form $Q$ on the
self-adjoint part of $T^*_e(\chi)$:
\begin{equation}\label{positiveform}
    Q(e[D,a])e=(e[D,a])^2 \qqq a\in \cA\,.
\end{equation}
It is non-degenerate since when $e[D,a]$ is self-adjoint,
$(e[D,a])^2=0$ implies $e[D,a]=0$. Let then $C_Q$ be the Clifford
algebra associated to the quadratic form $Q$ on the self-adjoint
part of $T^*_e(\chi)$. The latter has real dimension $\geq p$ and
the relations \eqref{cliffrelcenterbis} show that the map $
e[D,a]\mapsto e[D,a]$ gives a representation of $C_Q$ in the complex
vector space $eS_\chi$.  Thus this shows that the dimension of
$eS_\chi$ is then at least equal to $2^{p/2}$. The hypothesis of the
Theorem on the multiplicity of  the action of $\cA''$ in $\cH$
shows, using the \axiom of absolute continuity, that the fiber
dimension of $S$ is $2^{p/2}$. This shows that $e=1$ and also, since
the complexification of the algebra $C_Q$ is an $N\times N$ matrix
algebra for $N\geq 2^{p/2}$, that $M_\chi= \End S_\chi$ for every
$\chi \in X$. It also shows that the dimension of $T^*(\chi)$ is
equal to $p$ and that on $U_\alpha$ the $[D,a_\alpha^j]$ form a
basis of $T^*(\chi)$. Consider then the following monomials
$$
\mu_F=[[D,a_\alpha^{j_1}],[D,a_\alpha^{j_2}],\ldots,[D,a_\alpha^{j_k}]]
$$
where $F=\{j_1<j_2<\ldots<j_k\}$ is a subset with $k$ elements of
$\{1,2,\ldots,p\}$. For every $\chi\in U_\alpha$ the $\mu_F$ form a
basis of $M_\chi= \End S_\chi$. Thus any element $T$ of $M_\chi=
\End S_\chi$ can be uniquely written in the form
\begin{equation}\label{tendo}
T=\sum a_F\, \mu_F
\end{equation}
The coefficients $a_F$ can be computed using the normalized trace on
$\End S_\chi$, the $\mu_F$ and the element $T$. Thus using the
conditional expectation $E_\cA$ of \eqref{tracemap} one gets, for
any
 endomorphism $T$ of $\cH_\infty$ with support in $U_\alpha$, that
 \eqref{tendo} holds with coefficients $a_F\in \cA$. This shows that
 any
 endomorphism $T$ of $\cH_\infty$ is a  polynomial in the $[D,a]$ with
coefficients in $\cA$ and it follows that it is automatically
regular. Thus the strong regularity holds and we can apply Theorem
\ref{regthm}. To see that $X$ is a spin$^c$ manifold one uses
\cite{CoSM} (see \cite{FGV} for the detailed proof).
\endproof

\medskip \section{Final remarks}\label{sectfinal}

\subsection{The role of $D$}

By Lemma \ref{insideweak}, the spectral triple $(\cA,\cH,D)$ is
entirely determined by $(\cA'',\cH,D)$ where $M=\cA''$ is the
commutative von Neumann algebra weak closure of $\cA$. It follows in
particular that, except for the dimension $N$ of the bundle $S$
which we may assume, for simplicity, to be constant and equal to
$2^{p/2}$, there is no information in the pair $(\cA'',\cH)$: they
are all pairwise isomorphic. Similarly the only invariant of the
pair $(\cH,D)$ is the spectrum of $D$ \ie a list of real numbers
with multiplicity. By \cite{Milnor} this spectrum does not suffice
to reconstruct the geometry, and it is natural to wonder what
additional invariant is required to do so. As we shall briefly
explain it is the relative position of $M$ and of the self-adjoint
operator $D$ which selects one geometric space, and it is worthwhile
to look at the \axioms from this point of view. The analogue in our
context of the geodesic flow is the following one parameter group
\begin{equation}\label{geodesicflowdefn}
\gamma_t(T)=e^{it|D|} Te^{-it|D|} \qqq T\in \cL(\cH)\,.
\end{equation}
We say\footnote{\cf Lemma \ref{geodcinfty} of \S \ref{regularity}}
that an operator $T\in \cL(\cH)$ is of class $C^\infty$ when the map
from $\R$ to $ \cL(\cH)$ given by $t\mapsto \gamma_t(T)$ is of class
$C^\infty$ (for the norm topology of $\cL(\cH)$) and we denote by
$C^\infty(\cH,D)$ this subalgebra of $\cL(\cH)$. This algebra only
depends upon $(\cH,D)$ and does not yet measure the compatibility of
$(M,\cH)$ and $(\cH,D)$. This is measured by the {\em weak density}
in $M$ of
\begin{equation}\label{cinftyalg}
C^\infty(M,\cH,D)=\{ T\in M\cap C^\infty(\cH,D) \,| [D,T]\in M'\cap
C^\infty(\cH,D)\}
\end{equation}
where $M'$ is the commutant of $M$. One checks that
$\cA=C^\infty(M,\cH,D)$ is a subalgebra of $M$ and its size measures
the compatibility of $(M,\cH)$ and $(\cH,D)$.

 We now come to two
equations which assert that $\cA=C^\infty(M,\cH,D)$ is large enough,
so that we have maximal compatibility. One checks that
$\cH_\infty=\cap \Dom D^m$ is automatically a module over $\cA$ (for
the obvious action). The first equation requires that this module is
{\em finite and projective} and that it admits a hermitian structure
$(\;|\;)$ (necessarily unique) such that:
\begin{equation}\label{first}
\langle  \xi ,a \,\eta  \rangle = \cutint \, (\xi |\eta)\,a \,
|ds|^{p} \,, \ \forall   a \in \cA   ,   \forall  \xi ,\eta \in
\cH_{ \infty}
\end{equation}
where $\cutint$ is the noncommutative integral given by the Dixmier
trace.

The second equation  means that we can find an element $c$ of the
tensor power $\cA^{\otimes n}$, $n=p+1$, totally antisymmetric in
its last $p$-entries, and such that
\begin{equation}\label{second}
c(D) =1\,, \ {\rm where} \ \  (a_0\otimes \cdots \otimes
a_p)(D)=a_0[D,a_1]\cdots [D,a_p]\qqq a_j\in\cA\,.
\end{equation}
(This assumes $p$ odd, in  the even case one requires that for some
$c$ as above $c(D)=\gamma$ fulfills $\gamma = \gamma^* ,  \gamma^2 =
1, \gamma D =-D\gamma$). We can now restate  Theorem
\ref{fundthmspinc} as:

\begin{thm} \label{compdisccont} Let $(M,\cH,D)$ fulfill \eqref{cinftyalg}, \eqref{first} and
\eqref{second}, and $N=2^{p/2}$, then there exists a unique smooth
compact oriented spin$^c$ Riemannian manifold $(X,g)$ such that the
triple $(M,\cH,D)$ is given by
\begin{itemize}
  \item $M=L^\infty(X,dv)$ where $dv$ is the Riemannian volume form.
  \item $\cH=L^2(X,S)$ where $S$ is the spinor bundle.
  \item $D$ is a Dirac operator associated to the Riemannian metric $g$.
\end{itemize}
\end{thm}

\proof We let $\cA=C^\infty(M,\cH,D)$.  By the weak density in $M$
of \eqref{cinftyalg}, we know that the multiplicity of the action of
$\cA''=M$ in $\cH$ is $N=2^{p/2}$. By construction the triple
$(\cA,\cH,D)$ fulfills the first three \axioms. The fourth and fifth
follow from \eqref{first} and \eqref{second}. Thus by Theorem
\ref{fundthmspinc} we get that $\cA=C^\infty(X)$ for a smooth
oriented compact spin$^c$ manifold $X$. The conclusion then follows
from \cite{CoSM} (see \cite{FGV} for the detailed proof). Note that
there is no uniqueness of $D$ since we only know its principal
symbol. This is discussed in \cite{CoSM} and \cite{FGV}.
\endproof

A striking feature of the above formulation is that the full
information on the geometric space is subdivided in two pieces
\begin{enumerate}
  \item The list of eigenvalues of $D$.
  \item The unitary relation $F$ between the Hilbert space of the canonical
  pair $(M,\cH)$ and the Hilbert space of the canonical pair
  $(\cH,D)$.
\end{enumerate}
Of course the conceptual meaning of the unitary $F$ is the Fourier transform,
but this second piece of data is now playing a role entirely similar to that of
the CKM matrix in the Standard Model \cite{mc2}. Moreover, in the latter, the
information about the Yukawa coupling of the Higgs fields with the Fermions
(quarks and leptons) is organised in a completely similar manner, namely 1) The
masses of the particles 2) The CKM (and PMNS) matrix. At the conceptual level,
such matrices describe the relative position of two different bases in the same
Hilbert space. They are encoded by a double coset space closely related to
Shimura varieties (\cite{mc2}). These points deserve further investigations and
will be pursued in a forthcoming paper.

\smallskip

\subsection{Finite propagation}

One can use in the above context a result of Hilsum \cite{hilsum} to
obtain:

\begin{lem} The support of the kernel $k_t(x,y)$ of the operator
$e^{itD}$ is contained in
$$
\{(x,y)\in X^2\,|\, d(x,y)\leq |t|\}
$$
where the distance $d$ is defined as
$$
d(x,y)=\sup |h(x)-h(y)|\,,\Vert [D,h]\Vert\leq 1
$$
\end{lem}

\proof Let $(x,y)\in X^2$ with $d(x,y)> |t|$. There exists $h=h^*$
in $\cA$ such that $\Vert [D,h]\Vert\leq 1$ and $h(y)-h(x)> |t|$.
Also $h$ and $[D,h]$ commute by the order one condition. Thus by
Lemma 1.10 of \cite{hilsum}, one has $b<h(y)$, $a>h(x)$ such that:
$$
(h-b)_+ e^{-itD}(h-a)_-=0
$$
so that $k_t(x,y)=0$.
\endproof

\smallskip
\subsection{Immersion versus embedding}

The proof of Theorem \ref{regthm} shows that, with the cycle $c$
given by \eqref{cyclec}, the  map $\psi$ from $X$ to $\R^N$ given by
the components $a_\alpha^j$ for $j\geq 1$ is an {\em immersion}. It
is not however an embedding in general even if one includes the
components $a_\alpha^0$. To see this consider  open balls $B\subset
\R^p$ and $B_1\subset B$ such that for some translation $v$ the ball
$B_2=B_1+v$ is disjoint from $B_1$ and contained in $B$. Then let
$x^\mu$ be the coordinates in $\R^p$ and $a_1^j\in C_c^\infty(B)$ be
such that
$$
a_1^j(x)=x^j \qqq x\in B_1\,,\ \ a_1^j(x)=x^j-v^j  \qqq x\in B_2
$$
Let $N$ be a neighborhood of the complement of $B_1\cup B_2$ in $B$.
Let then $a_2^j(x)=b(x) x^j$ where $b(x)=1$ for all $x\in N$ and
vanishes in an open set of the form $B'_1\cup B'_2$ where the
$B'_j\subset B_j$ are smaller concentric balls. Let $a_\alpha^0$ be
a partition of unity in $B$ for the covering by $B_1\cup B_2$ and
$N$. Then let $c$ be the antisymmetrization of
$$
\sum_1^2 a_\alpha^0\otimes a_\alpha^1\otimes \cdots\otimes
a_\alpha^p
$$
For $x\in B'_1$ all the $a_2^j$ vanish, including $a_2^0$, and
$a_1^0=1$ so that the following equality shows that the map $\psi$
is not injective:
$$
a_1^j(x+v)=a_1^j(x)\qqq x\in B'_1\,.
$$

\smallskip
\subsection{The antisymmetry condition}

We have used throughout the stronger form of \axiom 4) where the
Hochschild cycle $c\in H_p(\cA,\cA)$ is assumed to be totally
antisymmetric in its last $p$-entries. It is unclear that one can
relax the antisymmetry condition on $c$. It is not true in general
for commutative algebras that any Hochschild class can be
represented in this way, but this is the case for $\cA=C^\infty(X)$.
In general, one has a natural projection on the antisymmetric
chains, given by the antisymmetrisation map $P$. It is defined by
the equality
$$
P(a_0\otimes a_1 \otimes \cdots \otimes a_p)=\frac{1}{p!} \sum_\beta
\epsilon(\beta) a_0\otimes a_{\beta (1)}\otimes \cdots \otimes
a_{\beta(p)}
$$
Its range is contained in $Z_p(\cA,\cA)$ since $\cA$ is commutative
and any antisymmetric chain is a cycle (\cite{Loday} Proposition
1.3.5). It is not obvious that $P$ maps Hochschild boundaries to
Hochschild boundaries. This follows from the equality
$$
P=\frac{1}{p!}\, \varepsilon_p\circ \pi_p
$$
where one lets $\Omega_K^p=\wedge^p_{\cA}\Omega_K^1$ be the
$\cA$-module of K\"{a}hler $p$-forms (\cf \cite{Loday} 1.3.11) and
\begin{equation}\label{epsiandpi}
   \varepsilon_p\,:\,\Omega_K^p\to H_p(\cA,\cA)\,, \ \ \pi_p \,:\,
H_p(\cA,\cA)\to \Omega_K^p
\end{equation}
are defined in \cite{Loday} Proposition 1.3.12 and Proposition
1.3.15. They are given by
\begin{equation}\label{pimap}
    \pi_k(a_0\otimes a_1\otimes \cdots \otimes a_k)=a_0da_1\wedge \cdots \wedge da_k
\end{equation}
and
\begin{equation}\label{epsilonmap}
    \varepsilon_k(a_0da_1\wedge \cdots \wedge
    da_k)=\sum\epsilon(\sigma) a_0\otimes a_{\sigma(1)}\otimes \cdots \otimes
    a_{\sigma(k)}
\end{equation}

\smallskip
\subsection{Strong regularity}

The hypothesis of strong regularity is, in general, stronger than
regularity. Indeed  the operation of  direct sum
$(\cA,\cH_1\oplus\cH_2,D_1\oplus D_2)$ of two spectral triples for
the same algebra $\cA$ preserves regularity but not, in general,
strong regularity.

\begin{prop} Assuming regularity  the subalgebra $Z_D(\cA)$
of $\End_\cA(\cH_\infty)$ generated by $\cA$ and the
$[D,b][D,c]+[D,c][D,b]$ for $b, c \in \cA$ is a commutative algebra
containing $\cA$ and commuting with $[D,a]$ for all $a\in \cA$.
\end{prop}

\proof This follows from Remark \ref{modda} since \eqref{commsymbd}
shows that $[D,b]^2$ commutes with $[D,a]$ for all $a\in \cA$.
\endproof

The understanding of the general situation when one does not assume
strong regularity should be an interesting problem since the
inclusion $\cA\subset Z_D(\cA)$ should correspond to a finite
``ramified cover" of the corresponding spectra, with $Y=\Sp
Z_D(\cA)$ covering $X=\Sp \cA$. It is easy to construct examples
where $Y$  has singularities. It is not clear that, assuming the
first five \axiomss, the space $X$ is always smooth. Similarly it is
unclear what happens if one relaxes the regularity condition to the
Lipschitz regularity, since we made heavy use of at least
$C^{1+\epsilon}$-regularity in the above proofs. Finally it would be
interesting also to investigate the meaning of {\em real
analyticity} of the space $X$ in terms of the real analyticity of
the geodesic flow \eqref{geodesicflowdefn}.

\smallskip
\subsection{The Noncommutative case}

 Among the five \axioms of \S \ref{prelem} the \axioms 1), 3) and 5)
make no use of the commutativity of the algebra $\cA$ and they
extend as such to the noncommutative case. We refer to \cite{CoSM}
for the extension of \axioms 2) and 4) to the noncommutative case.
The extension of the order one condition involves a new key
ingredient which is an antilinear unitary operator $J$ in $\cH$
which encodes the nuance between spin and spin$^c$. It turns out to
be an incarnation not only  of the charge conjugation in physics
terms and of the needed ``real structure" to refine the K-theoretic
meaning of the spectral triple from ordinary $K$-homology to
$KO$-homology but, at a deeper level, of the Tomita operator which
plays in the noncommutative case the role of a substitute for
commutativity. All this plays an important role in the
noncommutative geometry understanding of the standard model
\cite{mc2}, \cite{cc5}, \cite{cc6}. The extension of the
``orientability" \axiom 4) exists and it certainly holds \eg for
noncommutative tori (\cite{CoSM}) but it is not fully satisfactory
yet and its clarification should be considered as an open question.

\medskip

\section{Appendix 1: Regularity}\label{regularity}

The \axiom of regularity is not easy to check for smooth manifolds
since it involves the module of the operator $D$. We give below the
equivalent formulation in terms of $D^2$ (\cf \cite{cmindex}).

We deal with operators $T$ which act on $\cH_\infty=\cap \Dom D^n$.
We say that $T$ is {\em bounded} when
\begin{equation}\label{normT}
    \Vert T\Vert =\sup\{ \Vert T\xi\Vert\,|\,\xi \in \cH_\infty\,,\ \Vert
    \xi\Vert\leq 1\}
\end{equation}
is finite. We still denote by $T$ the unique continuous extension to
a bounded operator in $\cH$.
 By self-adjointness of $D$ the
domain $\cH_\infty$ is a core for powers of $D$ or of $|D|$. The
derivation $\delta(T)=[|D|,T]$ is defined algebraically as an
operator in $\cH_\infty$. The relation with the commutator in $\cH$
is given as follows.

\begin{lem}\label{delta1def} Assume that both $T$ and $[|D|,T]$ are
bounded (as in \eqref{normT}). Then $T$ preserves $\Dom |D|=\Dom D$
and the bounded extension of $[|D|,T]$ coincides with the commutator
$|D|T-T|D|$ on $\Dom |D|$.
\end{lem}

\proof Let $\xi \in \Dom |D|$. There exists a sequence
$\xi_n\in\cH_\infty$ with $\xi_n\to \xi$ and $|D|\xi_n\to |D|\xi$.
Since $T$ is bounded the sequences $T\xi_n$ and $T|D|\xi_n$ are
convergent and converge to $T\xi$ and $T|D|\xi$. Since $[|D|,T]$ is
bounded the sequence $(|D|T-T|D|)\xi_n$ converges. Thus $|D|T \xi_n$
converges, and as $|D|$ is closed one gets that $T\xi$ is in the
domain of $|D|$. Thus $\Dom |D|$ is invariant under $T$. Moreover
one has $ |D|T \xi=(|D|T-T|D|)\xi+T |D|\xi$.
\endproof

In other words, saying that   both $T$ and $[|D|,T]$ are bounded is
equivalent, for operators acting in $\cH_\infty$ to $T\in \Dom
\delta$ and moreover $\delta(T)$ is then the bounded extension of
$|D|T-T|D|$.

We introduce the following variant of $\delta$, defined on operators
$T$ acting in $\cH_\infty$,
\begin{equation}\label{delta1bis}
    \delta_1(T)=[D^2,T](1+D^2)^{-1/2}
\end{equation}
\begin{lem}\label{delta1delta} Let $T$ acting in $\cH_\infty$ be
bounded.
\begin{enumerate}
  \item If $\delta_1(T)$ and $\delta_1^2(T)$ are bounded so is
  $\delta(T)$.
  \item The $\delta_1^n(T)$  are bounded for all $n$ iff so are the $\delta^n(T)$.
\end{enumerate}
\end{lem}

\proof (1) The module $|D|$ is given by the following integral, which makes
sense when applied to any $\xi\in \Dom D$, which we omit for simplicity
\begin{equation}\label{modD}
    |D|=\frac 2\pi\, \int_0^\infty\frac{D^2}{D^2+u^2}\,du
\end{equation}
To avoid dealing with the kernel of $D$ we use the decomposition
$\delta=\delta'+\delta_0$ where the derivations $\delta'$ and
$\delta_0$ commute, and $\delta_0$ is bounded,
\begin{equation}\label{delta'}
    \delta'(T)=[Q,T]\,, \ \ Q=D^2(1+D^2)^{-1/2}
\end{equation}
\begin{equation}\label{delta0}
   \delta_0(T)=[f_0(D),T]\,, \ \ f_0(x)=|x|-x^2(1+x^2)^{-1/2}\qqq
   x\in \R
\end{equation}
One has $f\in C_0(\R)$ and the derivation $\delta_0$ is bounded, in fact $\Vert
\delta_0\Vert\leq 1$ since $\Vert f_0\Vert_\infty <1/2$. One has
\begin{equation}\label{modDbet}
    Q=\frac 2\pi\, \int_0^\infty\frac{D^2}{D^2+1+u^2}\,du
\end{equation}
Thus
\begin{equation}\label{modD1}
\delta'(T)=[Q,T]=\frac 2\pi\,
\int_0^\infty[\frac{D^2}{D^2+1+u^2},T]\,du
\end{equation}
$$
[\frac{D^2}{D^2+1+u^2},T]=-[\frac{1+u^2}{D^2+1+u^2},T]=(1+u^2)\frac{1}{D^2+1+u^2}[D^2,T]\frac{1}{D^2+1+u^2}
$$
$$
=[D^2,T]\frac{1+u^2}{(D^2+1+u^2)^2}\,
+(1+u^2)[\frac{1}{D^2+1+u^2},[D^2,T]]\frac{1}{D^2+1+u^2}
$$
Thus using
$$
[\frac{1}{D^2+1+u^2},[D^2,T]]=-\frac{1}{D^2+1+u^2}[D^2,[D^2,T]]\frac{1}{D^2+1+u^2}
$$
we get
$$
[\frac{D^2}{D^2+1+u^2},T]=[D^2,T]\frac{1+u^2}{(D^2+1+u^2)^2}\,
-\frac{1}{D^2+1+u^2}[D^2,[D^2,T]]\frac{1+u^2}{(D^2+1+u^2)^2}
$$
Thus combining with \eqref{modD1}  one gets
\begin{equation}\label{modD3}
\delta'(T)=\frac 12 [D^2,T](1+D^2)^{-1/2}+\frac 12
[D^2,T](1+D^2)^{-3/2}
\end{equation}
$$
-\frac 2\pi\, \int_0^\infty
\frac{1}{D^2+1+u^2}[D^2,[D^2,T]]\frac{1+u^2}{(D^2+1+u^2)^2}du
$$
where we used
$$
\frac 2\pi\, \int_0^\infty\frac{u^2}{(D^2+1+u^2)^2}du=\frac 12
(1+D^2)^{-1/2}\,,
$$
and
$$
\frac 2\pi\, \int_0^\infty\frac{1}{(D^2+1+u^2)^2}du=\frac 12
(1+D^2)^{-3/2}\,.
$$
 Now one has $[D^2,[D^2,T]]=\delta_1^2(T)(1+D^2)$ and
$$
\Vert \frac{(1+u^2)(1+D^2)}{(D^2+1+u^2)^2}\Vert \leq 1
$$
so that
$$
\Vert   \frac{1}{D^2+1+u^2}[D^2,[D^2,T]]\frac{1+u^2}{(D^2+1+u^2)^2}
\Vert \leq   \Vert   \frac{1}{D^2+1+u^2}\Vert \Vert
\delta_1^2(T)\Vert
$$
and one gets
\begin{equation}\label{esti1bis}
    \Vert \delta'(T)\Vert \leq \Vert \delta_1(T)\Vert
    +\Vert \delta_1^2(T)\Vert\,.
\end{equation}
Now if both $\delta_1(T)$ and $\delta_1^2(T)$ are bounded, we get
that $\delta'(T)$ is bounded and since $\delta=\delta'+\delta_0$
with $\delta_0$ bounded, we get that $\delta(T)$ is bounded, with
\begin{equation}\label{estiweak}
    \Vert \delta(T)\Vert \leq \sum_0^2\Vert \delta_1^j(T)\Vert
    \,.
\end{equation}

\smallskip

(2) The operations $\delta$ and $\delta_1$ commute since $|D|$
commutes with $D^2$. Let us assume that the $\delta_1^n(T)$  are
bounded. We have seen that $\delta(T)$ is bounded. To show that
$\delta^2(T)$ is bounded it is enough to show that
$\delta_1^m(\delta(T))$ are bounded for $m=1,2$. But
$\delta_1^m(\delta(T))=\delta(\delta_1^m(T))$ which is bounded since
the $\delta_1^n(\delta_1^m(T))$ are bounded for $n\leq 2$, $m\leq
2$. More generally let us show by induction on $n$ an inequality of
the form
\begin{equation}\label{indinequ}
    \Vert \delta^n(T)\Vert \leq \sum_0^{2n} c_{n,k}\Vert \delta_1^k(T)\Vert
\end{equation}
To get it for $n+1$, assuming it for $n$, one uses \eqref{estiweak}
which gives
$$
\Vert \delta( \delta^n(T))\Vert \leq
\sum_0^2\Vert\delta_1^j(\delta^n(T))\Vert=\sum_0^2\Vert\delta^n(\delta_1^j(T))\Vert
$$
$$
\leq \sum_0^2\sum_0^{2n} c_{n,k}\Vert \delta_1^k(\delta_1^j(T))\Vert
$$
Thus we obtain by induction that $\delta^n(T)$ is bounded.

Conversely, the boundedness of the $\delta^n(T)$ implies that of the
$\delta_1^n(T)$. Indeed  the boundedness of the $\delta^n(T)$ is
equivalent to  the boundedness of the ${\delta''}^n(T)$ where
$\delta''(T)=[(1+D^2)^{1/2},T]$ since $|D|-(1+D^2)^{1/2}$ is bounded
and commutes with $|D|$. Moreover the square of the operation
$$
T\mapsto (1+D^2)^{1/2}T(1+D^2)^{-1/2}=T+\delta''(T)(1+D^2)^{-1/2}
$$
is
$$
T\mapsto (1+D^2)T(1+D^2)^{-1}=T+[D^2,T](1+D^2)^{-1}
$$
which gives
$$
[D^2,T](1+D^2)^{-1}=2\delta''(T)(1+D^2)^{-1/2}+
{\delta''}^2(T)(1+D^2)^{-1}
$$
so that
$$
\delta_1(T)=2\delta''(T)+ {\delta''}^2(T)(1+D^2)^{-1/2}
$$
and one can proceed as above to get the boundedness of the
$\delta_1^n(T)$. \endproof

\medskip

Finally we relate the regularity condition with the smoothness of
the geodesic flow $t\to \gamma_t(T)=e^{it|D|}Te^{-it|D|}$ of
\eqref{geod}.

\begin{lem} \label{geodcinfty} Let $T\in \cL(\cH)$, then the following conditions are
equivalent:
\begin{enumerate}
  \item  $T\in \cap_m \Dom\,\delta^m$.
  \item $t\to \gamma_t(T)$ is of class $C^\infty$ in the norm
  topology.
\end{enumerate}
\end{lem}

\proof Let us show that $(1)\Rightarrow (2)$. By \eqref{commudel}
$T$ preserves $\cH_\infty$. We write the Taylor formula with
remainder
\begin{equation}\label{taylor}
    f(t)=f(0)+tf'(0)+\ldots+
    \frac{t^n}{n!}f^{(n)}(0)+\frac{t^{n+1}}{n!}\int_0^1(1-u)^n
    f^{(n+1)}(tu)du
\end{equation}
for the function $f(t)=e^{it|D|}Te^{-it|D|}\xi$ with $\xi\in
\cH_\infty$. Since $T$ preserves $\cH_\infty$ this function is of
class $C^\infty$. One gets
\begin{equation}\label{taylor1}
\gamma_t(T)\xi=T\xi +
it\delta(T)\xi+\ldots+\frac{i^nt^n}{n!}\delta^n(T)\xi+\frac{i^{n+1}t^{n+1}}{n!}\int_0^1(1-u)^n
    \gamma_{tu}(\delta^{(n+1)}(T))\xi du
\end{equation}
since $f^{(k)}(s)=\gamma_s(\delta^{(k)}(T))\xi$ by induction on $k$.
This shows that $t\to \gamma_t(T)$ is of class $C^\infty$ in the
norm    topology, since the norm of the remainder is $O(t^{n+1})$.
Let us show the converse $(2)\Rightarrow (1)$. It is enough to show
that if $T\in \cL(\cH))$ and the following limit exists in norm
$\lim_{t\to 0} \frac 1t(\gamma_t(T)-T)$, then $T\in \Dom \delta$ and
the limit is $i\delta(T)$. One has, for $\xi \in \cH$,
\begin{equation}\label{domainofabsD}
    \xi\in \Dom |D|\Leftrightarrow \exists \lim_{t\to 0}\,\frac 1t
    (e^{it|D|}\xi-\xi)
\end{equation}
where the limit is supposed to exist in norm. Assuming that for some
bounded operator $Y\in \cL(\cH)$ one has $\lim_{t\to 0}\Vert \frac
1t(\gamma_t(T)-T)-Y\Vert=0$, one gets for any $\xi\in \Dom |D|$,
that $\frac 1t
    (e^{it|D|}T\xi-T\xi)\to iT|D|\xi+Y\xi$. This shows that $T\xi\in
    \Dom |D|$ and that $i|D|T\xi=iT|D|\xi+Y\xi$ which gives the required equality.
\endproof

\medskip

\section{Appendix 2: The Dixmier trace and the heat
expansion}\hfill\medskip\label{appendix}

We first recall the basic properties of the Dixmier trace. Recall
that  the characteristic value $\mu_n(T)$ of a compact operator $T$
is the $n$-th eigenvalue of $|T|$ arranged in decreasing order and
is equal to
\begin{equation}\label{minimax}
   \inf\{\Vert T|_{E^\perp}\Vert\,|\, {\rm dim}\, E=n-1\}
\end{equation}

\begin{defn}\label{weylnormdef} We define the \weyl by
\begin{equation}\label{weylnorm}
    \sigma_N(T)=\sum_1^{N}\mu_n(T)
\end{equation}
\end{defn}

The fact that they  are norms and in particular fulfill
\begin{equation}\label{sigmannorm}
\sigma_N(T_1+T_2)\leq \sigma_N(T_1)+\sigma_N(T_1)
\end{equation}
 follows from the next statement in which we use the same notation
for a subspace $E\subset \cH$ and the orthogonal projection on that
subspace.

\begin{prop} One has
\begin{equation}\label{weylnorm1}
    \sigma_N(T)=\sup \{ \Vert TE\Vert_1\,|\, {\rm dim}\, E=N\}
\end{equation}
Let $T$ be a positive operator, then
\begin{equation}\label{weylnorm2}
    \sigma_N(T)=\sup \{ \Tr(TE)\,|\, {\rm dim}\, E=N\}
\end{equation}
\end{prop}

We use the following notation for refined limiting processes,

\begin{defn} \label{defomegalim}
With the Ces\`aro mean $M$ defined by
\begin{equation}\label{cesaromean1}
M(f)(\lambda )=\frac{1}{{\log
\lambda}}\int_1^{\lambda}f(u)\frac{{\hbox{d} u}}{u}.
\end{equation}
and $h(\lambda)$ a bounded function of $\lambda>0$, $\omega$ a
linear form on $C_b(\R_+^*)$ which is positive, $\omega(1)=1$, and
vanishes on $C_0(\R_+^*)$, and $\phi$ an homeomorphism of $\R_+^*$,
we define:
\begin{equation}\label{limtwo}
    {\rm Lim}^k_{\phi(\lambda)\to\omega}h(\lambda)=\omega(M^k(g))\,,\ \
    g(\lambda)=h(\phi^{-1}(\lambda))
\end{equation}
where the upper index $k$ indicates that we iterate the Ces\`aro
mean $k$-times.
\end{defn}

\medskip

 We write ${\rm Lim}_{\omega}$ as an abbreviation for ${\rm
 Lim}^1_{\omega}$, and when we apply it to a sequence $(\alpha_N)_{N\in \N}$ we mean
 that the sequence has been extended to a function using
\begin{equation}\label{seqtofunc}
f_{\alpha}(\lambda )=\alpha_N\hbox{~for }\lambda\in \ ]N-1 ,N].
\end{equation}
Also we consider the two-sided ideal containing compact operators of
order one,
\begin{equation}\label{loneinfty}
{\mathcal L}^{(1,\infty )}(\cH)=\{T\in {\mathcal K};\hbox{~}\sigma_
N(T)=O(\log N)\}\,.
\end{equation}

\bigskip

\begin{defn} \label{defdixtrace}   For $T\geq 0$, $T\in {\mathcal
L}^{(1,\infty )}(\cH)$, we set
\begin{equation}\label{dixtracedef}
{\Tr}_{\omega}(T)={\rm Lim}_{\omega}\,\,\frac{1}{{\log N}}\sum_{
n=1}^{N}\mu_n(T).
\end{equation}
\end{defn}

The basic properties of the Dixmier trace ${\Tr}_{\omega}$ are
summarized in the following (\cite{Co-book} Proposition 3,
IV.2.$\beta$):

\begin{prop}\label{propbasicdixtrace} ${\Tr}_{\omega}$ extends uniquely by linearity to the entire
ideal $ {\mathcal L}^{(1,\infty )}(\cH)$ and has the following
properties:

\begin{itemize}
\item[(a)]   If $T\geq 0$ then $\Tr_{\omega}(T)\geq 0$.

\smallskip

\item[(b)]    If $S$ is any bounded operator and $\ T\in {\mathcal L}^{
(1,\infty )}(\cH)$\/,   then $\Tr_{\omega}(ST)=\Tr_{\omega}(TS)$.

\smallskip

\item[(c)] ${\Tr}_{\omega}(T)$   is independent of the choice of the inner product on~$
\cH$,   \ie   it depends only on the Hilbert space~$\cH$ as a
topological vector space.

\smallskip

\item[(d)] ${\Tr}_{\omega}$  vanishes on the ideal ${\mathcal L}^{
(1,\infty )}_0(\cH)$,   which is the closure, for the
$\Vert\hbox{~}\Vert_{1,\infty}$-norm,   of the ideal of finite-rank
operators.
\end{itemize}
\end{prop}

We fix $p\in [1,\infty]$. Let $D$ be a self-adjoint unbounded
operator such that its resolvent is an infinitesimal of order $1/p$,
\ie such\footnote{We replace $D$ by a non-zero constant on its
kernel so that $D^{-1}$ makes sense.} that
$\mu_n(D^{-1})=O(n^{-1/p})$. We shall compare ${\Tr}_{\omega}
(T|D|^{-p})$ and $\lim \epsilon^p\,\Tr (f(\epsilon D)T)$. We let
$E_N$ be the spectral projection\footnote{This is ambiguous when
there is spectral multiplicity.} on the first $N$-eigenvectors of
$|D|$ so that $\dim E_N=N$, $E_N< E_{N+1}$ and
\begin{equation}\label{endefn}
\Tr (E_N |D|^{-p})=\sigma_N( |D|^{-p})\,.
\end{equation}

\begin{lem} \label{lemheatequdix} For any bounded operator $T\in \cL(\cH)$ one has
\begin{equation}\label{heatequdix0}
    \limom\frac{1}{\log N}\,\Tr (E_N |D|^{-p}T)
    = {\Tr}_{\omega} (T|D|^{-p})
\end{equation}
\end{lem}

\proof The hypothesis on $D$ shows that $\Tr (E_N |D|^{-p})=O(\log
N)$. Moreover, by construction of the Dixmier trace, one has
\begin{equation}\label{heatequtotal}
\limom\frac{1}{\log N}\,\Tr (E_N |D|^{-p})=\limom\frac{1}{\log
N}\,\sigma_N( |D|^{-p})
    = {\Tr}_{\omega} (|D|^{-p})
\end{equation}
Let $\phi(T)$ be the left-hand side of \eqref{heatequdix0}. It makes
sense since $$|\Tr (E_N |D|^{-p}T)|\leq \Tr (E_N |D|^{-p})\Vert
T\Vert=O(\log N)$$ so that the sequence $\frac{1}{\log N}\,\Tr (E_N
|D|^{-p}T)$ is bounded. The functional $\phi$ on $\cL(\cH)$ is
linear and positive (the trace of the product of the two positive
operators $E_N |D|^{-p}$ and $T$ is positive). Let $\psi(T)$ be the
right-hand side of \eqref{heatequdix0}. Proposition
\ref{propbasicdixtrace} shows that, since $|D|^{-p}\in{\mathcal L}^{
(1,\infty )}(\cH)$, the functional $\psi$ is a positive linear
functional on $\cL(\cH)$. One uses Proposition
\ref{propbasicdixtrace}, (b) to check the positivity, using for
$T\geq 0$,
$$
{\Tr}_{\omega} (T|D|^{-p})={\Tr}_{\omega}
(T^{1/2}|D|^{-p}T^{1/2})\geq 0\,.
$$
Let us show that for any $T\geq 0$ one has $\phi(T)\leq \psi(T)$.
One has
$$
\sigma_N(T^{1/2}|D|^{-p}T^{1/2})=\sigma_N(|D|^{-p/2}T|D|^{-p/2})
$$
using $A=|D|^{-p/2}T^{1/2}$ in
\begin{equation}\label{inequmun1}
    \mu_n(A^*A)=\mu_n(AA^*) \qqq A\in \cK\,,\ n\in \N\,.
\end{equation}
 Thus one gets
\begin{equation}\label{dixfine}
\psi(T)={\Tr}_{\omega} (T|D|^{-p})=\limom \frac{1}{\log
N}\,\sigma_N( |D|^{-p/2}T|D|^{-p/2})
\end{equation}
 By
\eqref{weylnorm2}, one has
$$
\sigma_N(|D|^{-p/2}T|D|^{-p/2})=\sup \{
\Tr(|D|^{-p/2}T|D|^{-p/2}E)\,|\, {\rm dim}\, E=N\}$$ $$\geq
\Tr(|D|^{-p/2}T|D|^{-p/2}E_N)= \Tr (E_N |D|^{-p}T)
$$
since $E_N$ and $|D|^{-p/2}$ commute. Thus
$\sigma_N(|D|^{-p/2}T|D|^{-p/2})\geq \Tr (E_N |D|^{-p}T)$ and after
dividing by $\log N$ and applying $\limom$ to both sides one gets
the inequality $\phi(T)\leq \psi(T)$. But, by \eqref{heatequtotal},
$\phi(1)={\Tr}_{\omega} (|D|^{-p})=\psi(1)$ and thus the positive
functional $\theta=\psi-\phi$ is equal to $0$, by the Schwartz
inequality $|\theta(T)|^2\leq \theta(T^*T)\theta(1)$.
\endproof

With $|D|$ as above, we let as in \eqref{spectralprojbis}, for any
$\lambda>0$,
\begin{equation}\label{spectralproj}
    P(\lambda)={\bf 1}_{[0,\lambda]}(|D|)\,, \ \ \alpha(\lambda)=
    \Tr P(\lambda)\,.
\end{equation}

\begin{lem} \label{lemheatequdix2} Assume that
\begin{equation}\label{nonzerodix}
\liminf \lambda^{-p} \alpha(\lambda)>0
\end{equation}
Then, for any bounded operator $T\in \cL(\cH)$ one has
\begin{equation}\label{heatequdix0bis}
  p\;  \limom\frac{1}{\log N}\,\Tr (E_N |D|^{-p}T)
    = \limo_{\lambda^p\to\omega}\frac{1}{\log \lambda}\,\Tr (P(\lambda) |D|^{-p}T)
\end{equation}
\end{lem}

\proof We can assume by linearity that $T\geq 0$. We have (using
\eqref{upperboundalpha}) constants $c_1>0$ and $c_2<\infty$ such
that:
\begin{equation}\label{lambdaandn}
c_1\lambda^p\leq \alpha(\lambda)\leq c_2 \lambda^p
\end{equation}
We let:
\begin{equation}\label{lambdaandn1}
f(N)=\Tr (E_N |D|^{-p}T)\,,\ \ g(\lambda)=\Tr (P(\lambda)
   |D|^{-p}T)
\end{equation}
 Since $\dim P(\lambda)\leq N$ implies $P(\lambda)\leq E_N$ we
get, using $P(\lambda)
   |D|^{-p}\leq E_N |D|^{-p}$,
\begin{equation}\label{lambdalessn}
   f(N)\geq g(\lambda) \qqq \lambda\,, c_2 \lambda^p\leq N \,.
\end{equation}
Similarly, since $\dim P(\lambda)\geq N$ implies $P(\lambda)\geq
E_N$ we get
\begin{equation}\label{lambdalargen}
f(N)\leq g(\lambda)\qqq \lambda\,, c_1 \lambda^p\geq N \,.
\end{equation}
We extend $f(N)$ to positive real values of $N$ as a non-decreasing
step function. The arbitrariness of the extension is irrelevant
since $f(N+1)-f(N)\to 0$ when $N\to \infty$  and we are interested
in $\limom\frac{1}{\log N}f(N)$ which is insensitive to bounded
perturbations of $f$. By construction, the Ces\`aro mean satisfies
the following scale invariance, for bounded functions~$ f$,
\begin{equation}\label{cesaromean2}
\vert M(\theta_{\mu} (f))(\lambda )-M(f)(\lambda )\vert\to
0\hbox{~~as } \lambda\to\infty,
\end{equation}
 where $\mu >0$ and $\theta_{\mu} (f)(\lambda )=f(\mu^{-1}\lambda )$ $
\forall \lambda\in \R^{*}_{+}$. It follows from \eqref{lambdalessn}
and \eqref{lambdalargen} that $f(c_1N)\leq g(N^{1/p})\leq f(c_2N)$
up to $o(N)$ and for any positive real $N$. Thus the scale
invariance of the Ces\`aro mean \eqref{cesaromean2}, together with
$\log N/\log cN\to 1$ gives:
\begin{equation}\label{cesequiv}
    M(\frac{1}{\log N}f(N))- M(\frac{1}{\log N}g(N^{1/p}))\to 0
\end{equation}
so that
$$
\limom \frac{1}{\log N}f(N)=\frac 1p\limom \frac{1}{\log
N^{1/p}}g(N^{1/p})
$$
and the required equality \eqref{heatequdix0bis} follows from
Definition \ref{defomegalim}.
\endproof

\begin{cor} Assuming \eqref{nonzerodix}, one has
\begin{equation}\label{heatequdix0bet}
  p\; {\Tr}_{\omega} (T|D|^{-p})
    = \limo_{\lambda^p\to\omega}\frac{1}{\log \lambda}\,\Tr (P(\lambda)
    |D|^{-p}T)\qqq T\in \cL(\cH)\,.
\end{equation}
\end{cor}

\proof This follows from Lemmas \ref{lemheatequdix} and
\ref{lemheatequdix2}.\endproof

\medskip

\begin{thm} \label{thmheatequdix} Assume \eqref{nonzerodix}.
Let $f\in C_c([0,\infty[)$. Let $\rho=p\int_0^\infty u^{p-1}f(u)du$.
Then for any bounded operator $T\in \cL(\cH)$ one has
\begin{equation}\label{heatequdix}
    \limtwo_{\epsilon^{-p}\to\omega}\epsilon^p\,\Tr (f(\epsilon |D|)T)
    =\rho\,{\Tr}_{\omega} (T|D|^{-p})
\end{equation}
\end{thm}

\proof Let $g(u)=u^pf(u)$ viewed as an integrable function on the
multiplicative group $\R_+^*$, endowed with its normalized Haar
measure $d^*u=\frac{du}{u}$. We can assume that $T\geq 0$. We
consider the positive measure on $\R_+^*$ given by $d\beta(\lambda)$
where
\begin{equation}\label{posmeas}
    \beta(\lambda)=\Tr (P(\lambda)
    |D|^{-p}T)
\end{equation}
which is a non-decreasing step function of $\lambda$. The measure
$d\beta$ is a positive linear combination of Dirac masses, $
d\beta=\sum\alpha_n\delta_{\lambda_n} $.
 One has
$$
d\beta(\lambda)=\Tr (dP(\lambda)
    |D|^{-p}T)=\lambda^{-p}\Tr (dP(\lambda)
    T)
    $$
    $$
    \epsilon^p\,\Tr (f(\epsilon |D|)T)=\epsilon^p\,\int f(\epsilon
    \lambda)\Tr (dP(\lambda)
    T)=\int \epsilon^p\lambda^pf(\epsilon
    \lambda)d\beta(\lambda)
    $$
    so that:
\begin{equation}\label{posmeas1}
    \epsilon^p\,\Tr (f(\epsilon |D|)T)=\int g(\epsilon
    \lambda)d\beta(\lambda)
\end{equation}
The  convolution of the measure $d\beta$ with the function $\tilde
g(u)=g(u^{-1})$ makes sense, since both have support in an interval
$[u_0,\infty[$ with $u_0>0$,  and gives the function
\begin{equation}\label{posmeas3}
   (\tilde g\star d\beta)(u)=\int g(u^{-1}
    \lambda)d\beta(\lambda)\,.
\end{equation}
Thus, with $h(\epsilon)=\epsilon^p\,\Tr (f(\epsilon |D|)T)$, one
gets using \eqref{posmeas1},
\begin{equation}\label{posmeas4}
 h(u^{-1})=(\tilde g\star d\beta)(u)
 \end{equation}
The convolution of the measures $\tilde g (u)d^*u$ and $d\beta$ is
absolutely continuous with respect to $d^*u$ and is given, with
$\theta_u(v)=uv$ for all $u,v>0$, by
\begin{equation}\label{posmeas2}
   (\tilde g\star d\beta)d^*u= \int \tilde
   g(u)\theta_u(d\beta)d^*u\,.
\end{equation}
We extend the definition of the Ces\`aro mean \eqref{cesaromean1} to
measures $\mu$ by:
\begin{equation}\label{cesaromean1meas}
M(\mu)(\lambda )=\frac{1}{{\log \lambda}}\int_1^{\lambda}d\mu\,,
\end{equation}
so that
\begin{equation}\label{cesaromean1measbis}
M(\mu)(\lambda )=M(h)(\lambda )\  {\rm for}\ \mu=h \,d^*u\,.
\end{equation}
One has $\beta(\lambda)=O(\log\lambda)$ since $$
\beta(\lambda)\leq \Tr(P(\lambda)|D|^{-p})\Vert T\Vert\leq
\int_0^\lambda u^{-p}d\alpha(u)\Vert T\Vert $$ while $\alpha(u)=0$
near $0$, and $u^{-p}\alpha(u)$ is bounded by \eqref{lambdaandn}.
This gives after integrating by parts $$\int_0^\lambda
u^{-p}d\alpha(u)=\lambda^{-p}\alpha(\lambda)+\int_0^\lambda
p\,u^{-p-1}\alpha(u)du\leq c_2(1+p\log\lambda) +c'
$$
Moreover  for $v>1$ one gets, by the above integration by parts,
\begin{equation}\label{lambv}
\beta(v\lambda)-\beta(\lambda)\leq \Vert
T\Vert\int_\lambda^{v\lambda}u^{-p}d\alpha(u)\leq \Vert
T\Vert\,c_2(1+\log v)\,.
\end{equation}
One has
\begin{equation}\label{posmeas6}
 M(d\beta)(\lambda)=\frac{1}{\log \lambda}\int_1^{\lambda}d\beta=\frac{1}{\log
 \lambda}(\beta(\lambda)-\beta(1))
 \end{equation}
 Thus one has constants $a$ and $b$ such that for any $u$,
\begin{equation}\label{posmeas5prel}
 |M(\theta_u(d\beta))(\lambda)-M(d\beta)(\lambda)|=\leq (a+b|\log u|)(\log\lambda)^{-1}
 \end{equation}
Thus since $\tilde g(u)$ and $|\log u|\tilde g(u)$ are integrable,
\begin{equation}\label{posmeas5}
 M(\int \tilde g(u)\theta_u(d\beta)d^*u)-M(d\beta)\,\int \tilde g(u)d^*u\to 0
 \end{equation}
 Equivalently, using \eqref{posmeas4}, \eqref{posmeas2}, \eqref{cesaromean1measbis}
  and $\int \tilde g(u)d^*u=\int
 g(u)d^*u$,
 \begin{equation}\label{posmeas5.5}
 M(\tilde h)-M(d\beta)\,\int g(u)d^*u\to 0 \,, \ \ \tilde
 h(u)=h(u^{-1})\,.
 \end{equation}
 Now by \eqref{posmeas6} and \eqref{heatequdix0bet} one has
 \begin{equation}\label{firstM}
 \limo_{\lambda^p\to\omega} M(d\beta)(\lambda)= p\; {\Tr}_{\omega} (T|D|^{-p})
 \end{equation}
 Thus we finally get
 \begin{equation}\label{secondM}
 (p\int g(u)d^*u) {\Tr}_{\omega} (T|D|^{-p})=\limo_{\lambda^p\to\omega} M(\tilde h)(\lambda)
 \end{equation}
The right hand side is given, by definition, by
$$
\limo_{\lambda^p\to\omega} M(\tilde h)(\lambda)=\omega(M(k)(u))\,, \
\ k(u)=M(\tilde h)(u^{1/p})
$$
Thus we still need to compare $k(u)=M(\tilde h)(u^{1/p})$ with
$k_1(u)=M(\tilde h(\lambda^{1/p}))(u)$, but a simple computation
shows that $k(u)=k_1(u)$.\endproof

\medskip

\begin{cor} \label{corheatequdix} Assume \eqref{nonzerodix}.
Let $f\in C_c([0,\infty[)^+$ be a positive function. Let
$\rho=p\int_0^\infty u^{p-1}f(u)du$. One has, when $\epsilon \to 0$,
\begin{equation}\label{heatequdix1}
    \liminf\epsilon^p\,\Tr (f(\epsilon |D|)T)
    \leq \rho\,{\Tr}_{\omega} (T|D|^{-p})\,.
\end{equation}
\end{cor}

\proof Let $\delta=\liminf\epsilon^p\,\Tr (f(\epsilon |D|)T)$. Then
for any $c<1$ one has $h(\epsilon)=\epsilon^p\,\Tr (f(\epsilon
|D|)T)\geq c\delta$ for $\epsilon\leq \epsilon_c>0$. It follows that
$\limtwo_{\epsilon^{-p}\to\omega}h(\epsilon)\geq c\delta$. Thus by
\eqref{heatequdix} one has $c\delta\leq \rho\,{\Tr}_{\omega}
(T|D|^{-p})$ and one gets \eqref{heatequdix1}.
\endproof

\end{document}